\documentclass[12pt,a4paper]{amsart}
\usepackage{fontenc,amsmath,mathrsfs,amsfonts,amsthm, verbatim,amssymb, MnSymbol, enumerate, multicol, hyperref, geroENG, graphicx, color}

\usepackage[english]{babel}
\usepackage{boondox-cal}
\usepackage{setspace}
\usepackage[left=2.5 cm, bottom=2.5 cm, right=2.0 cm, top = 2.0 cm]{geometry}
\usepackage{times}
\allowdisplaybreaks
\begin{document}
\title{Metrical properties of Hurwitz Continued Fractions}

\author{Yann Bugeaud}
\address{I.R.M.A., UMR 7501, Universit\'e de Strasbourg
et CNRS, 7 rue Ren\'e Descartes, 67084 Strasbourg Cedex, France}
\address{Institut universitaire de France}
\email{bugeaud@math.unistra.fr}

\author{Gerardo Gonz\'alez Robert}
\address{Department of Mathematical and Physical Sciences,  La Trobe University, Bendigo 3552, Australia.}
\email[Corresponding Author]{G.Robert@latrobe.edu.au, matematicas.gero@gmail.com}

\author{Mumtaz Hussain}
\address{Department of Mathematical and Physical Sciences,  La Trobe University, Bendigo 3552, Australia. }
\email{m.hussain@latrobe.edu.au}

\date{}
\maketitle

\begin{abstract}
We develop the geometry of Hurwitz continued fractions, a major tool in understanding the approximation properties of complex numbers by ratios of Gaussian integers. Based on a thorough study of the geometric properties of Hurwitz continued fractions, among other things, we determine that the space of valid sequences is not a closed set of sequences. 
Additionally, we establish a comprehensive metrical theory for Hurwitz continued fractions.

Let $\Phi:\Na\to \RE_{>0}$ be any function. For any complex number $z$ and $n\in\Na$, let $a_n(z)$ denote the $n$th partial quotient in the Hurwitz continued fraction of $z$. One of the main results of this paper is the computation of the Hausdorff dimension of the set
\[E(\Phi)
\colon=
\left\{ z\in \mathbb C: |a_n(z)|\geq \Phi(n) \text{ for infinitely many }n\in\mathbb{N} \right\}.
\]
This study is a complex analog of a well-known result of Wang and Wu [Adv. Math. 218 (2008), no. 5, 1319--1339].
\end{abstract}
\tableofcontents

\section{Introduction}
A classical result states that for each real number $x$ there is a sequence $(A_n(x))_{n\geq 0}$ of integers, where $A_n(x)\geq 1$ for $n\geq 1$, such that
\[
x=
A_0(x) + \cfrac{1}{A_1(x) + \cfrac{1}{A_2(x)  + \cfrac{1}{\ddots}}}.
\]
Such expansion, which is infinite if and only if $x$ is an irrational number, is known as the regular (real) continued fraction of $x$. 
Regular continued fractions are a fundamental tool for understanding irrational numbers, partly due to their ability to provide the best rational approximations of real numbers \cite{Cas1957}. However, their applications go beyond this property. They enable us to draw connections with various branches of mathematics such as dynamical systems, ergodic theory, and fractal geometry, thus unveiling shared characteristics among different classes of irrational numbers 
\cite{MR4557871, MR3418243, FaLiWaWu2009, MR3713001, MR46583, MR2263718}.
We refer the reader to \cite{Khi1997} for the classical theory of regular continued fractions and to \cite{Bug2004} for a modern account.

The Borel-Bernstein Theorem (see, for example, \cite[Theorem 1.11]{Bug2004}) and similar results provide us with families of null sets (Lebesgue measure zero) determined by their continued fraction expansions and, hence, by their approximation properties. A classic example is the set of badly approximable numbers; this is precisely the set of real numbers with bounded partial quotients. Hausdorff and packing dimensions have been used to determine the size of such sets. The research in this vein dates back to Jarník \cite{Jar1928}, who proved that the set of badly approximable numbers has full Hausdorff dimension. Good's celebrated work \cite{Goo1941} concerns the size of small null sets obtained by imposing asymptotic conditions on the partial quotients (for refinements and extensions, see \cite{Luc1997, WanWu2008}).  The Hausdorff dimension of the null sets appearing in the Borel-Bernstein Theorem was calculated by Wang and Wu in their seminal paper \cite{WanWu2008}.
Other sets related to Diophantine approximation, such as the Markov and Lagrange spectrum, have also been studied with the aid of regular continued fractions \cite{Mor2018}.

Regular continued fractions link analytic and algebraic properties of irrational numbers. Among the classical results in this direction, we have Lagrange's Theorem on quadratic surds \cite[Theorem 28]{Khi1997} or Liouville's construction of transcendental numbers, where transcendence is obtained by ensuring good rational approximation \cite[Section 9]{Khi1997}. More recently, Adamczewski and Bugeaud used Schmidt's Subspace Theorem to build families of badly approximable transcendental numbers \cite{AdaBug2005,Bug2013}.  

Extending the theory of continued fraction expansions to higher dimensions while preserving most of the aforementioned properties is, thus, highly desirable. However, such extensions are far from straightforward and may not always be possible. In this paper, we focus on the complex plane $\Cx$. In $\Cx$, several continued fraction algorithms have been suggested depending on the preferred properties  \cite{Che21, Hoc2020, Hur1887, Lev1952, SchmidtA1975, Sha1979}. In all of these expansions, the resulting space of sequences is a major concern. However, Shallit's work \cite{Sha1979} completely describes the sequence space of his algorithm using automata theory. 
In \cite{Hur1887}, Adolf Hurwitz introduced a complex continued fraction algorithm. 
Hurwitz continued fractions (defined in Section \ref{Section:HCF_Basics}) give a continued fraction expansion for complex numbers that share several helpful properties with regular continued fractions (see, for example, propositions \ref{Propo2.1} and \ref{Propo3.2}). 
In particular, they provide good Gaussian rational approximations to a given complex number. {In \cite{Lak1974}, Lakein relied on these approximation properties to obtain a new proof of Ford's Theorem on approximation to complex numbers by elements of $\QU(i)$.} A word of caution: while drawing an analogy between Hurwitz and the nearest-integer continued fraction algorithm may seem more natural, we choose to focus on the more widely studied regular continued fractions.


The theory of Hurwitz continued fractions has recently been used in complex Diophantine approximation, yet it is far from being as well understood as their real counterpart. 
In \cite[Section 7]{DanNog2014}, Dani and Nogueira extended their work \cite{DanNog2002} on real binary quadratic forms to complex numbers using Hurwitz continued fractions. In \cite{Sim2016}, Simmons studied Diophantine approximation properties of Hurwitz continued fractions when restricted to real numbers by comparing them with regular continued fractions.
Hines \cite{Hin2019} and Gonz{\'a}lez Robert \cite{Gero2020-01} showed independently that badly approximable complex numbers have bounded Hurwitz continued fraction partial quotients. Hines' work, however, includes other quadratic imaginary fields. 
There are some (rather surprising) results on the partial quotients of complex algebraic numbers. Hensley \cite{Hen2006} was the first to find examples of algebraic complex numbers of degree four over $\QU(i)$ and with bounded complex partial quotients. This was subsequently generalized by Bosma and Gruenewald \cite{BosGru2012}, who showed that, for every even integer $d$, there exist algebraic complex non-real numbers of degree $d$ over $\QU$ for which the Hurwitz continued fraction expansion has bounded partial quotients. In \cite{Gero2020-01}, Gonz{\'a}lez Robert constructed a family of badly approximable transcendental complex numbers giving a partial analog of \cite{Bug2013}. 
In \cite{HeXio2021-02}, He and Xiong calculated the Hausdorff dimension of sets of complex numbers that can be approximated by quotients of Gaussian integers at a given order $\psi$ but no better. This is a complex version of a problem solved by Bugeaud \cite{Bug2003}. Lastly, the theory developed in this paper was applied by Garc{\'i}a-Ramos, Gonz{\'a}lez Robert, and Hussain to study topological, dynamical, and descriptive set-theoretic properties of Hurwitz continued fractions \cite{GRGRH24}. Some of the applications in their paper include determining the rank in the Borel hierarchy of the set of Hurwitz normal numbers with respect to the complex Gauss measure and the construction of a family of
complex transcendental numbers with bounded partial quotients. 

In this paper, we compute the Hausdorff dimension of sets of complex numbers whose Hurwitz continued fraction has a subsequence of partial quotients of fast growth. Consequently, we obtain a complex analog of the Borel-Bernstein refinement by Wang-Wu \cite{WanWu2008}.
A significant part of our work is dedicated to studying two intimately related aspects of Hurwitz continued fractions: their geometric properties and their symbolic space. In particular, we show that the space of sequences,  
the partial quotients in the Hurwitz continued fraction of a complex number, is not a closed subset  of the space of sequences of Gaussian integers $\Za[i]^{\Na}$ (see the proof of Lemma \ref{PROPO:GC01:Rho:01}). However, we can replace it with the closure of a subset $\sfR$ of $\Za[i]^{\Na}$. The set $\sfR$ will be defined in simple topological terms in Section \ref{Section:HCF_Basics}. 

Let us state our main results. 
For any $t\in\RE$, we let $\lfloor t\rfloor$ be the integer part of $t$, that is
\[
\lfloor t\rfloor \colon= \max \{m\in\Za: m\leq t\}.
\]
For each $z\in \Cx$, define the nearest Gaussian integer $[z]$ to $z$ as 
\[
[z]\colon= 
\left\lfloor \real(z) + \frac{1}{2}\right\rfloor + i \left\lfloor \imag(z) + \frac{1}{2}\right\rfloor.
\]
Define the set
\[
\mfF\colon= \{z\in \Cx: [z]=0\}.
\]
We will restrict our attention to $\mfF$.
Note that $\mfF$ is the unit square centered at the origin whose sides are parallel to the coordinate axes, it includes the left and the bottom sides and excludes the remaining two. Also, the Gaussian integral translates of $\mfF$ form a partition of $\Cx$.  
The rigorous definitions of Hurwitz continued fractions and some related terms are fully explained in Section \ref{Section:HCF_Basics}. For now, we just need to know that the Hurwitz continued fraction process associates to each complex number $z$ a (finite or infinite) sequence $(a_n(z))_{n\geq 0}$ in $\Za[i]$ such that
\begin{equation}\label{Eq:HCF:01}
z
=
[a_0(z);a_1(z),a_2(z),\ldots]
=
a_0(z) + \cfrac{1}{a_1(z) + \cfrac{1}{a_2(z) + \cfrac{1}{\ddots}}}.
\end{equation}
If $z\in\mfF$, then $a_0(z)=0$ and we denote the continued fraction in \eqref{Eq:HCF:01} by $[a_1(z),a_2(z),\ldots]$ or simply $[a_1,a_2,\ldots]$. 
For any $B>1$, we write
\begin{equation}\label{Eq:Defsb}
s_B
\colon=
\lim_{n\to\infty} \inf\left\{\rho\geq 0: \sum_{(a_1, \ldots, a_n)\in\sfR(n)} \left(\frac{1}{B^n |q_n(a_1,\ldots,a_n) |^2}\right)^\rho  <1\right\}.
\end{equation}

Let $\dimh(A)$ denote the Hausdorff dimension of a given subset $A\subseteq \mathbb{C}$, we refer the reader to \cite{Falconer2014} for an account of the theory of Hausdorff dimension.
Write $\|z\|\colon=\max\{ |\real(z)|, |\imag(z)|\}$ for any $z\in\Cx$.

\begin{teo01}\label{TEO:MAIN}
For any function $\Phi:\Na\to \RE_{>0}$, define the set 
\[
E_{\infty}(\Phi)
\colon=
\left\{ z\in \mathfrak{F}: \|a_n(z)\|\geq \Phi(n) \text{ \rm for infinitely many }n\in\mathbb{N} \right\}
\]
and $B\geq 1$ by
\[
\log B = \liminf_{n\to\infty} \frac{\log \Phi(n) }{n}.
\]
The Hausdorff dimension of $E_{\infty}(\Phi)$ is as follows:
\begin{enumerate}[\rm i.]
\item If $B=1$, then $\dimh(E_{\infty}(\Phi))=2$.
\item If $1<B<\infty$, then $\dimh(E_{\infty}(\Phi))=s_B$.
\item When $B=\infty$, let $b$ be given by $$\log b = \displaystyle\liminf_{n\to\infty} \frac{\log\log \Phi(n)}{n}.$$
Then, we have $$\dimh(E_{\infty}(\Phi))=\frac{2}{1+b}.$$
\end{enumerate}
\end{teo01}
\begin{rema01}
Theorem \ref{TEO:MAIN} holds if we replace $\|\cdot\|$ with the complex absolute value $|\cdot|$. Indeed, for every $\Phi:\Na\to  \RE_{>0}$ and $c>0$, we have
\begin{align*}
\liminf_{n\to\infty} \frac{\log \Phi(n) }{n}
&=
\liminf_{n\to\infty} \frac{\log c\,\Phi(n) }{n}, \\
\liminf_{n\to\infty} \frac{\log\log \Phi(n)}{n}
&=
\liminf_{n\to\infty} \frac{\log\log  c\,\Phi(n) }{n}.
\end{align*}
Then, by Theorem \ref{TEO:MAIN}, $\dimh(E_{\infty}(\Phi))=\dimh(E_{\infty}(c\Phi))$. Putting
\[
E(\Phi) 
\colon=
\left\{ z\in \mathbb C: |a_n(z)|\geq \Phi(n) \text{ \rm for infinitely many }n\in\mathbb{N} \right\},
\]
we have
\[
E_{\infty} \left( 2^{-\frac{1}{2}} \Phi \right)
\subseteq 
E(\Phi)
\subseteq 
E_{\infty} (\Phi),
\]
and, thus, $$\dimh(E(\Phi)) =\dimh(E_{\infty}(\Phi)).$$
\end{rema01}
A crucial step towards proving Theorem \ref{TEO:MAIN} is the computation of the Hausdorff dimension of the sets $F(B)$ for $B>1$, defined as
\[
F(B)\colon= \{ z =[a_1,a_2,\ldots]\in\mathfrak{F}: \|a_n(z)\|\geq B^n \text{ for infinitely many } n\in\mathbb{N}\}.
\]
\begin{teo01}\label{TEO:DIMHB}
For each $1<B<\infty$, we have $\dimh(F(B))=s_B$.
\end{teo01}

In \cite{Luc1997}, {\L}uczak calculated the Hausdorff dimension of sets of real numbers with fastly growing partial quotients. Given $x\in (0,1)\setminus \QU$, let $[A_{1}(x),A_{2}(x),\ldots]$ denote its regular continued fraction. For $b,c>1$, define the sets
\begin{align*}
\widetilde{\Xi}_{\RE}(b,c) &\colon=  \{ x\in(0,1): c^{b^n}\leq A_n(x) \text{ for all } n\in\Na\}, \\
\Xi_{\RE}(b,c) &\colon=  \{ x\in(0,1): c^{b^n}\leq A_n(x) \text{ for infinitely many } n\in\Na\}.
\end{align*}
\begin{teo01}[{\L}uczak, \cite{Luc1997}]\label{TEO:LUCZAK:R}
For every pair of real numbers $b,c>1$, we have
\[
\dimh  \left(\widetilde{\Xi}_{\RE}(b,c)\right)
=
\dimh \left(\Xi_{\RE}(b,c)\right)
=
\frac{1}{1+b}.
\]
\end{teo01}

The case $B=\infty$ in Theorem \ref{TEO:MAIN} requires a complex version of Theorem \ref{TEO:LUCZAK:R}. 
\begin{teo01}\label{TEO:LUC:USUALABS}
If $b,c>1$, then
\begin{align*}
\frac{2}{1+b}
&=
\left\{ z\in\mfF: c^{b^n} \leq |a_{n}(z) |\, \text{ for all } n\in\Na\right\} \\
&=
\left\{ z\in\mfF: c^{b^n} \leq |a_{n}(z) | \, \text{ for infinitely many } n\in\Na\right\}.
\end{align*}
\end{teo01}
Theorem \ref{TEO:LUC:USUALABS} provides an alternative proof of the lower bound in Good's theorem for Hurwitz continued fractions \cite[Theorem 1.3]{Gero2020-02}; that is,
\[
\dimh\left(\left\{ z\in\mfF: \lim_{n\to\infty} |a_n(z)|=\infty \right\}\right) = 1.
\]
Certainly, note that for any pair of real numbers $b>1$ and $c>1$, we have
\[
\left\{ z\in\mfF: c^{b^n} \leq |a_n(z) | \text{ for all } n\in\Na \right\}
\subseteq 
\left\{ z\in\mfF: \lim_{n\to\infty} |a_n(z)|=\infty \right\}.
\]
Theorem \ref{TEO:LUC:USUALABS} implies that
\[
\frac{2}{1+b}\leq \dimh\left(\left\{ z\in\mfF: \lim_{n\to\infty} |a_n(z)|=\infty \right\}\right)
\]
and letting $b\to 1$ we conclude that
\[
1\leq \dimh\left(\left\{ z\in\mfF: \lim_{n\to\infty} |a_n(z)|=\infty \right\}\right).
\]
The upper bound can be shown with a natural covering argument.

The paper is organized as follows. Section \ref{SC:Notation} contains a summary of notation. 
In Section \ref{Section:TransAndInver}, we establish some terminology concerning Gaussian-integral translations and the complex inversion. 
In Section \ref{Section:HCF_Basics}, we define Hurwitz continued fractions and recall some of their basic properties; in particular, we define various types of cylinders (regular, irregular, extremely irregular, full, and almost full). 
In Sections  \ref{SEC:CYL} and \ref{Sec:GC}, we study the geometric properties of Hurwitz continued fractions.
More precisely, in Section \ref{SEC:CYL}, we discuss the symmetries and the size of regular cylinders. 
In Section \ref{Sec:GC}, we estimate how many regular cylinders of a given level are intersected by small discs. These estimate will be used for the lower bound in Theorem \ref{TEO:DIMHB}. The detailed explanation of the geometric aspects of Hurwitz continued fractions is of independent interest.
In Section \ref{SEC:IRRREG}, we discuss how we can approximate irregular numbers by regular numbers (defined in Section \ref{Section:HCF_Basics}). 
Sections \ref{Sec:TeodimBprf} and \ref{Sec:CxLuczak} are dedicated to the proofs of Theorems \ref{TEO:DIMHB} and \ref{TEO:LUC:USUALABS}, respectively. 
Section \ref{Sec:TeoMain} contains the proof of Theorem \ref{TEO:MAIN}.
Lastly, in \ref{SECTION:ConcludingRemarks}, we give some final remarks and open problems.

\medskip

\noindent{\bf Acknowledgements.} Yann Bugeaud is partially supported by project ANR-18-CE40-0018 funded by the French National Research Agency. The research of Mumtaz Hussain and Gerardo Gonz\'alez Robert is supported by the Australian Research Council Discovery Project (200100994). Gerardo Gonz\'alez Robert was partially funded by the program \textit{Estancias posdoctorales por M\'exico} of the {CONACyT}, Mexico. We thank the anonymous referees for useful comments and suggestions that led to significant improvements.

\section{Frequently used notation}\label{SC:Notation}
For convenience, we collect some of the frequently used notation. The vast majority of the symbols, however, are defined throughout the text.
\begin{enumerate}[1.]
\item If $(x_n)_{n\geq 1}$ and $(y_n)_{n\geq 1}$ are two sequences of positive numbers, we write $x_n\ll y_n$ if there is some constant $c>0$ such that $x_n\leq c y_n$ for all large $n\in\Na$. We write $x_n\asymp y_n$ if $x_n\ll y_n$ and $y_n\ll x_n$.

\item $\overline{\Cx}$ stands for the Riemann sphere.

\item We define $\iota:\overline{\Cx}\to \overline{\Cx}$ by $\iota(z)=z^{-1}$ and, for $a\in\Za[i]$,  $\tau_a:\overline{\Cx}\to\overline{\Cx}$ by $\tau_{a}(z)=z+a$. We consider the usual conventions for $z=0$ and $z=\infty$. The composition of these functions is denoted by concatenation.

\item Put $\alpha:=\frac{2-\sqrt{3}}{2}$ and $\zeta_1 = -\frac{1}{2} + i \alpha$.

\item \label{Not:Intervals} For all $z,w\in\Cx$,  write 
\[
[z,w)
:=
\{z + t(w-z): t\in [0,1)\}, 
\quad
[z,w]
:=
\{z + t(w-z): t\in [0,1]\}, 
\]
\[
(z,w)
:=
\{z + t(w-z): t\in (0,1)\}, 
\quad
(z,w]
:=
\{z + t(w-z): t\in (0,1]\}.
\]

\item For $n\in\Na$ and $\bfa\in \scD^n$, the functions $Q_n(\,\cdot\,;\bfa)$, $T_n(\,\cdot\,;\bfa):\overline{\Cx}\to\overline{\Cx}$ are given by
\[
T_n(\,\cdot\,;\bfa) \colon= \tau_{-a_n}\iota\cdots\tau_{-a_1}\iota,
\quad
Q_n(\,\cdot\,;\bfa) \colon= \iota\tau_{a_1}\cdots \iota\tau_{a_n}.
\]
\item For every $A\subseteq \Cx$, we denote the closure of $A$ by $\Cl(A)$, the interior of $A$ by $\inte(A)$ or $A^{\circ}$, and the boundary of $A$ by $\Fr(A)$.


\item We denote the empty word by $\epsilon$.

\item If $X$ is a non-empty finite set, $\#X$ denotes the number of elements contained in $X$.
\item If $X$, $Y$ are non-empty sets, $A\subseteq X$, and $f:X\to Y$ is any function, the image of $A$ under $f$ is $f[A]$; that is, $f[A]\colon=\{f(a):a\in A\}$.

\item For $z\in\Cx$, write $\|z\|\colon=\max\{|\real(z)|,|\imag(z)|\}$.

\item For $z\in\Cx$, write $\Pm(z)\colon=\min\{|\real(z)|,|\imag(z)|\}$.

\item $\scD\colon=\{a \in\Za[i]: |a|\geq \sqrt{2}\}$.

\item $\scE\colon=\{a\in \scD: |a|\geq \sqrt{8}\}$. 

\item For $M>0$, we define $\scD(M)\colon=\{a\in\scD:\|a\|\leq M\}$.

\item The functions $\Rota,\Mir_1,\Mir_2:\overline{\Cx}\to\overline{\Cx}$ are given by
\[
\Rota(z) = iz, 
\quad
\Mir_1(z) = \overline{z}, 
\quad\text{and }\quad
\Mir_2(z) = -\overline{z}
\quad\text{ for all } z\in \Cx
\]
and $\Rota(\infty)=\Mir_1(\infty)=\Mir_2(\infty)=\infty$.
\item The dihedral group of order $8$, which is generated by $\Rota$ and $\Mir_1$, is denoted by $\Di_8$. The subgroup of $\Di_8$ generated by $\Rota$ is denoted by $\RotaG$.

\item The orbit of $\zeta_1$ under the action of $\RotaG$ is denoted by $\RotaG(\zeta_1)$ and similarly for $\Di_8$; that is, $\RotaG(\zeta_1)=\{f(\zeta_1):f\in \RotaG\}$ and $\Di_8(\zeta_1)=\{f(\zeta_1):f\in \Di_8\}$.

\item For any two sets $A,B\subseteq \Cx$, we write $A\equiv B \mod \RotaG$ (resp. $A\equiv B \mod \Di_8$) if and only if there is some $f\in\RotaG$ (resp. $f\in\Di_8$) such that $A=f[B]$.

\item The distance between two non-empty sets $A,B\subseteq \Cx$ is
\[
\dist(A,B)
\colon=
\inf\{ |a-b|:a\in A,b\in B\}.
\]
\item \label{Not:DxCirc} For any $z\in \Cx$ and any $\rho>0$, we write
\[
\Dx(z;\rho) \colon= \{w\in \Cx: |z-w|<\rho\},
\quad
\overline{\Dx}(z;\rho) \colon= \{w\in \Cx: |z-w|\leq \rho\},
\]
and
\[
C(z,\rho)
\colon=
\{w\in \Cx: |z-w|=\rho\}.
\]
Following \cite{Lak1973}, we write 
\[
\Dx(z)\colon= \Dx(z;1),\quad
\overline{\Dx}(z)\colon=\overline{\Dx}(z;1).
\]
\item Given $k\in\Na$ and $\bfa=\sanu\in\scD^{\Na}$, the $k$-th \textbf{prefix} of $\bfa$ is $\prefi(\bfa;k)\colon= (a_1,\ldots, a_k)$.

\item \label{Not:ConcatWords} If $m,n\in\Na$, $\bfa=(a_1,\ldots,a_n)\in\scD^n$, and $\bfb=(b_1,\ldots,b_m)\in\scD^m$, we write
\[
\bfa\bfb \colon= (a_1,\ldots,a_n,b_1,\ldots,b_m)\in\scD^{n+m}.
\]

\item If $\bfa=\sanu$ is a sequence and $u,v$ are natural numbers such that $u<v$, we write $\bfa[u,v]=(a_{u}, a_{u+1}, \ldots, a_{v} )$.
\end{enumerate}

\section{Translations and the complex inversion}\label{Section:TransAndInver}
In this section, we set up the notation for the complex inversion and a family of translations which will be used throughout this paper. 
Let $\overline{\Cx}:=\Cx\cup\{\infty\}$ denote the Riemann sphere. The \textbf{inversion} $\iota:\overline{\Cx}\to \overline{\Cx}$ is the map given by 
\[
\forall z\in \overline{\Cx} \quad
\quad
\iota(z)
\colon=
\begin{cases}
0, & \text{ if } z=\infty, \\
\infty, & \text{ if } z=0, \\
\frac{1}{z}, & \text{ if } z\in\overline{\Cx}\setminus\{0, \infty\}.
\end{cases}
\]

For any $a\in \Za[i]$, the \textbf{translation} $\tau_{a}:\overline{\Cx}\to \overline{\Cx}$ is given by
\[
\forall z\in \overline{\Cx} \quad
\quad
\tau_a(z)
\colon=
\begin{cases}
z+a, & \text{ if } z\in\Cx, \\
\infty &\text{ if }z=\infty.
\end{cases}
\]
We represent the composition of these functions by juxtaposing them; for instance, $\tau_{a}\iota:=\tau_{a}\circ \iota$ for all $a\in \Za[i]$.
\begin{rema01}\label{Rema:FunDom}
The set $\mfF$ is a fundamental domain for the action of $\{\tau_a:a\in\Za[i]\}$ on $\Cx$. Then, the nearest Gaussian integer $[z]$ of a given $z\in \Cx$ is the only number $a\in\Za[i]$ satisfying one of the equivalent conditions:
\begin{enumerate}[\rm i.]
\item $z=\tau_{a}(w)$ for some $w\in \mfF$,
\item  $\tau_{-a}(z)\in\mfF$.
\end{enumerate}
\end{rema01}

Given $n\in\Na$ and $\bfa=(a_1,\ldots, a_n)\in\Za[i]^n$, define $T_n(\cdot;\bfa), Q_n(\cdot;\bfa):\overline{\Cx}\to \overline{\Cx}$ by 
\[
\forall z\in \overline{\Cx} 
\quad
T_n(z;\bfa)\colon= \tau_{-a_n}\iota \cdots \tau_{-a_1}\iota(z)
\;\text{ and }\;
Q_n(\cdot;\bfa)\colon=\iota\tau_{a_1}\cdots \iota \tau_{a_n}(z).
\]
Clearly, $T_n(\cdot,\bfa)$ and $Q_n(\cdot,\bfa)$ are mutually inverse homeomorphisms from $\overline{\Cx}$ onto itself. For any set $A\subseteq \overline{\Cx}$, we denote its closure by $\Cl(A)$, its interior by $\inte(A)$ and its image under $T_n(\cdot;\bfa)$ and $Q_n(\cdot;\bfa)$, respectively, by $T_n(A;\bfa)$ and $Q_n(A;\bfa)$; that is,
\[
T_n(A;\bfa)
:=
\left\{ T_n(z;\bfa): z\in A\right\}
\;\text{ and }\;
Q_n(A;\bfa)
:=
\left\{ Q_n(z;\bfa): z\in A\right\}.
\]
Observe that
\begin{align*}
\Cl(T_n(A;\bfa))
=
T_n(\Cl(A);\bfa), 
&\quad
\inte(T_n(A;\bfa))
=
T_n(\inte(A);\bfa),  \\
\Cl(Q_n(A;\bfa))
=
Q_n(\Cl(A);\bfa), 
&\quad
\inte(Q_n(A;\bfa))
=
Q_n(\inte(A);\bfa). 
\end{align*}
We will use these four identities without further reference. 

Lastly, we recall some well-known formulas concerning the action of $\iota$ on circles in the Riemann sphere. 
For $z\in \Cx$ and $\rho>0$, write
\[
C(z;\rho)
:=
\left\{w\in \Cx: |z-w|=\rho\right\}.
\]
First, if $k\in\RE\setminus\{0\}$, then
\begin{equation}\label{Eq:FormRecInv}
\iota\left[\{z\in \Cx: \real(z)=k\}\right]
=
C\left( \frac{1}{2k}; \frac{1}{2|k|}\right), 
\quad
\iota\left[ \{z\in \Cx: \imag(z)=k\} \right]
=
C\left( \frac{-i}{2k}; \frac{1}{2|k|}\right).
\end{equation}
These equations are implicitly used throughout the paper. 
If $z_0\in\Cx$ and $\rho>0$ are such that $\rho\neq |z_0|$, then
\begin{equation}\label{Eq:FormCircInv}
\iota\left[ C(z_0;\rho)\right]
=
C\left( \frac{\overline{z_0}}{|z|^2-\rho^2}; \frac{\rho}{|\rho^2 - |z_0|^2|}\right).
\end{equation}

\section{Basics on Hurwitz continued fractions}\label{Section:HCF_Basics}
In this section, we recall the definition of Hurwitz continued fractions and some of their basic properties. We define Hurwitz continued fractions in two steps: first, we restrict ourselves to numbers in $\mfF\setminus\{0\}$ and, then, consider numbers in $\Cx\setminus\{0\}$. 

Let the \textbf{Hurwitz-Gauss map} $T:\mfF\to \mfF$ be given by
\[
\forall z\in \mfF \quad
\quad
T(z) = 
\begin{cases}
z^{-1}- [z^{-1}], &\text{ if } z\neq 0,\\
0, &\text{ if } z=0.
\end{cases}
\]
For any $z\in \mfF\setminus \{0\}$, define $a_0(z) \colon= 0$, $a_1(z) \colon= [z^{-1}]$ and $a_n(z)\colon= a_1(T^{n-1}(z))$ for any $n\in\Na$ for which $T^{n-1}(z)\neq 0$. Naturally, the exponent in $T$ denotes iteration. For $z\in \Cx\setminus \mfF$, we define $a_0(z)=[z]$ and $a_n(z)=a_n(z - [z])$ whenever the right-hand side is defined. The \textbf{Hurwitz continued fraction} (Hcf) of $z\in\Cx\setminus\{0\}$ is the expression
\[
[a_0(z);a_1(z),a_2(z),\ldots]\colon= 
a_0(z) + \cfrac{1}{a_1(z) + \cfrac{1}{a_2(z) + \cfrac{1}{\ddots}}}.
\]
As mentioned in the introduction, when $z\in \mfF$, we have $a_0(z)=0$ and we write $[a_1(z),a_2(z),\ldots]$ (resp. $[a_1,a_2,\ldots]$) rather than $[0;a_1(z),a_2(z),\ldots]$ (resp. $[0;a_1,a_2,\ldots]$). As usual, we can interpret the previous expression as a limit (see Proposition \ref{Propo2.1}.\ref{Propo2.1_iii}). The terms of the sequence $(a_n(z))_{n\geq 1}$ are called the \textbf{Hurwitz partial quotients} of $z$ or simply the \textbf{partial quotients} of $z$. Let $\sepn$, $\seqn$ be the possibly finite sequences of Gaussian integers given by
\begin{align*}
p_0=0, \quad p_1 &= 1, \quad p_{n+1}=a_{n+1}p_n + p_{n-1} \text{ for }n\geq 1,  \\ 
q_0=1, \quad q_1 &= a_1, \quad q_{n+1}=a_{n+1}q_n + q_{n-1} \text{ for }n\geq 1.
\end{align*}
(see \cite[Section 2]{Khi1997}). Following \cite{DanNog2014}, we refer to $\sepn,\seqn$ as the $\clQ$\textbf{-pair} of $z$.

Hcfs produce an actual continued fraction representation of complex numbers. This is made precise in Proposition \ref{Propo2.1}. For the rest of the paper, we write
\[
\scD:=\Za[i]\setminus \{0,1,i,-1,-i\}.
\]
Sometimes, we refer to an element of $\scD^n$, $n\in\Na$, as a word of length $n$ over the alphabet $\scD$.
\begin{propo01}\label{Propo2.1}
Take any $z\in\Cx\setminus\{0\}$. Let $(a_n)_{n\geq 0}$ be the possibly finite sequence of partial quotients of $z$ and $\sepn$, $\seqn$ its $\clQ$-pair. 
\begin{enumerate}[\rm i.]
\item \label{Propo2.1_i} \cite[p. 73]{Hen2006} For any $n\in\Na$, if $a_n$ exists, then $a_n\in\scD$.
\item \label{Propo2.1_ii} \cite[Theorem 1]{Khi1997} For any $n\in\Na$, if $a_n$ exists, then
\[
\frac{p_n}{q_n} = [a_0;a_1,\ldots,a_n]
:=
a_0 + \cfrac{1}{a_1+ \cfrac{1}{ \ddots + \cfrac{1}{a_n}}}.
\]
\item \label{Propo2.1_iii} \cite[Theorem 6.1]{DanNog2014} The sequence $(a_n)_{n\geq 0}$ is infinite if and only if $z\in\Cx\setminus\QU(i)$; in this case,
\[
z
=
\lim_{n\to\infty} [a_0;a_1,a_2,\ldots,a_n].
\]
\end{enumerate}
\end{propo01}

\begin{figure}[h!]
\begin{center}
\includegraphics[scale=0.85,  trim={5cm 17.25cm 10cm 4cm},clip]{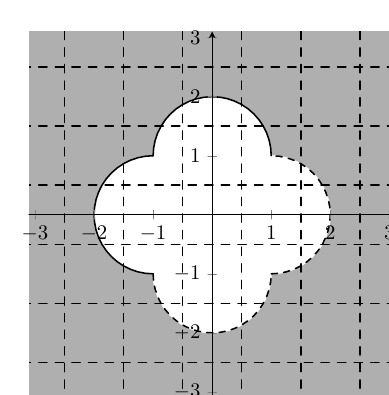}
\caption{ The set $\iota[\mfF]$. \label{Fig:02}}
\end{center}
\end{figure}

For readability and brevity, many of our arguments rely on pictures. Let us discuss Proposition \ref{Propo2.1}.\ref{Propo2.1_i} as an example. Figure \ref{Fig:02} depicts the set $\iota[\mfF]$. An explicit description of $\iota[\mfF]$ can certainly be obtained using \eqref{Eq:FormRecInv}. However, it is immediate from Figure \ref{Fig:02} that, if $a\in \Za[i]$, then
\[
 \tau_{a}[\mfF]\cap \iota[\mfF]\neq \vac
\quad \text{ if and only if }\quad
 a\in \scD.
\]
As a consequence, for any non-zero $z\in \mfF$ we have $a_1(z)=[z^{-1}]\in \scD$ (see Remark \ref{Rema:FunDom}) and, therefore, $a_n(z)\in\scD$ whenever it exists.

Proposition \ref{Propo3.2} collects some properties of Hcfs. Although these properties resemble well-established facts on regular continued fractions--except for part \ref{Propo3.2_iii}--the known proofs are much more involved. 

\begin{propo01}\label{Propo3.2}
Consider an arbitrary $z=[a_1,a_2,\ldots] \in\mfF\setminus\QU(i)$ with $\clQ$-pair $\sepn$, $\seqn$. The next assertions hold:
\begin{enumerate}[\rm i.]
\item \label{Propo3.2_i} \cite[p. 195]{Hur1887} The sequence $(|q_n|)_{n\geq 0}$ is strictly increasing. 
\item \label{Propo3.2_ii} \cite[Corollary 5.3]{DanNog2014} If $\psi\colon=\left( \frac{1+\sqrt{5}}{2} \right)^{1/2}$, then 
\[
|q_n|\geq \psi^{n-1}
\;\text{ for all }\;n\in\Na.
\]
\item \label{Propo3.2_iii} \cite[Proposition 3.3]{DanNog2014} We have
\[
z= \frac{(a_{n+1}+[a_{n+2},a_{n+3},\ldots])p_{n}+p_{n-1}}{(a_{n+1}+[a_{n+2},a_{n+3},\ldots])q_{n}+q_{n-1}}
\;\text{ for all } n\in\Na.
\]
\item \label{Propo3.2_iv} \cite[Theorem 1]{Lak1973}  We have
\[
\left| z - \frac{p_n}{q_n}\right| < \frac{1}{|q_n|^2}
\;\text{ for all } n\in\Na.
\]
\end{enumerate}
\end{propo01}
We have restricted  Proposition \ref{Propo3.2} to numbers in $\mfF\setminus\QU(i)$, but it can be easily adapted to $\Cx\setminus\{0\}$ (including finite sequences of partial quotients).

We may express the correspondence rule of the functions $Q_n(\,\cdot\,;\bfa)$ in  terms of the sequences $\sepn$ and $\seqn$. 

\begin{propo01}[{\cite[Theorem 5]{Khi1997}}]\label{Prop:QnCorrespRule}
Let $n\in\Na$ and $\bfa\in\sfR(n)$ be arbitrary. For every $z\in \Cx$ we have
\[
Q_n(z;\bfa)
= \frac{zp_{n-1} + p_{n}}{zq_{n-1} + q_{n}} 
= \frac{p_n}{q_n} + \frac{(-1)^nz}{q_n^2 \left( 1+ \frac{q_{n-1}}{q_n} z\right)}.
\]
Moreover, if $z\in\mfF_{n}(\bfa)$ (see Definition \ref{Def:CylProtSets} below), then
\[
Q_n(z;\bfa)
= 
[a_1,a_2,\ldots, a_{n - 1}, a_n + z ].
\]
\end{propo01}

\section{Cylinders and prototype sets}\label{SEC:CYL}
In this section, we introduce the building blocks of the metrical theory of Hcfs: cylinders and prototype sets. After their definitions, we classify them into six categories and provide examples justifying the need for this classification. Later, we study the symmetries of the family of one of these categories: the family of regular cylinders. 

\subsection{Definitions and basic properties}
In the definition below,  the term prototype set is borrowed from \cite{HeXio2021-01}.
\begin{def01}\label{Def:CylProtSets}
Given $n\in\Na$ and $\bfa=(a_1,\ldots,a_n)\in\scD^n$, the \textbf{cylinder} of level $n$ based on $\bfa$ is the set
\[
\clC_n(\bfa)\colon=
\{z \in\mfF: a_1(z)=a_1, \ldots, a_n(z)=a_n\}.
\] 
The \textbf{prototype set} $\mfF_n(\bfa)$ is $\mfF_n(\bfa)\colon= T^n[\clC_n(\bfa)]$.  
\end{def01}

\begin{figure}[h!]
\begin{center}
\includegraphics[scale=0.65,  trim={5cm 17cm 5cm 1.5cm},clip]{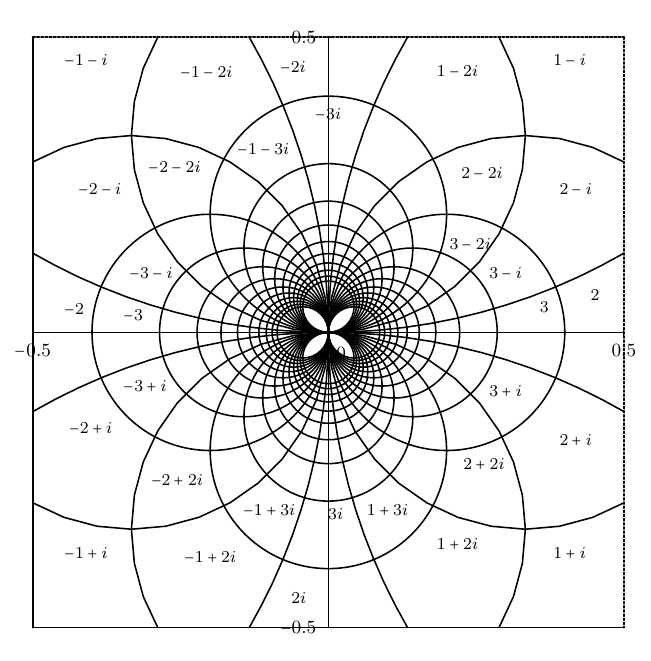}
\caption{ Cylinders of level $1$. \label{Fig:mfFPart}}
\end{center}
\end{figure}
Figure \ref{Fig:mfFPart} depicts cylinders of level $1$.
The next proposition gathers some properties of cylinders and prototype sets stemming directly from their definitions.

\begin{propo01}\label{Prop:ObsCyls}
Let $n\in\Na$ and $\bfa=(a_1,\ldots, a_n)\in\scD^{n}$ be arbitrary. 
\begin{enumerate}[\rm i.]
\item\label{Prop:ObsCyls:00} The restriction of $T_n(\,\cdot\,;\bfa)$ to $\clC_n(\bfa)$ is a homeomorphism from $\clC_n(\bfa)$ onto $\mfF_n(\bfa)$. Similarly, the restriction of $Q_n(\,\cdot\,;\bfa)$ to $\mfF_n(\bfa)$ is a homeomorphism from $\mfF_n(\bfa)$ onto $\clC_n(\bfa)$.
\item\label{Prop:ObsCyls:01} We have $\clC_1(a_1) = \iota\tau_{a_1}[\mfF]\cap \mfF$.
\item\label{Prop:ObsCyls:02} If $a_{n+1}\in\scD$, then 
\begin{equation*}
\clC_{n+1}(a_1, \ldots, a_n, a_{n+1} )
=
\clC_1(a_1)\cap \iota\tau_{a_1}[\clC_{n}(a_2,\ldots, a_{n+1})].
\end{equation*}
In general, for any $j\in\{1,\ldots,n\}$ we have
\[
\clC_{n+1}(a_1, \ldots, a_{n+1})
=
\clC_{j}(a_1,\ldots, a_j)\cap Q_{j}(\clC_{n+1-j}(a_{j+1},\ldots, a_{n+1});(a_1,\ldots, a_j) ).
\]
The assertion remains true if we replace $\clC$ with $\clCc$ (resp. with $\oclC$).
\end{enumerate}
\end{propo01}

\subsection{Taxonomy}
\begin{figure}[h]
\begin{center}
\includegraphics[scale=0.65,  trim={5cm 13cm 9cm 4.00cm},clip]{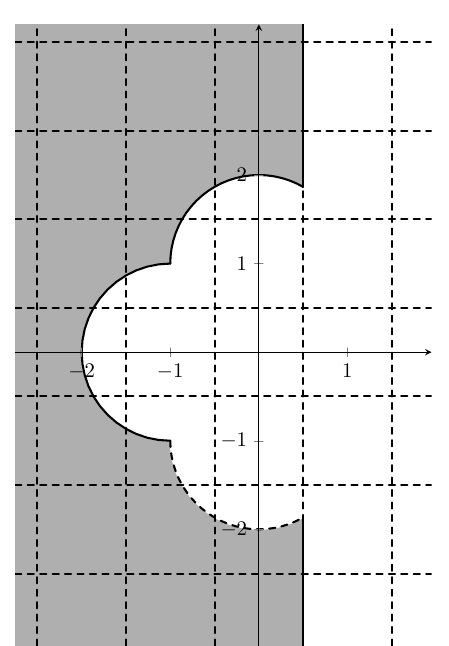}
\caption{ The set $\iota[\mfF_1(-2)]$. \label{Fig:iotaF1-2}}
\end{center}
\end{figure}

We introduce a taxonomy of cylinders, prototype sets, sequences, and numbers in $\mfF$. 
To motivate it, we discuss a particular examples of prototype sets using the notation in items \ref{Not:Intervals} and \ref{Not:DxCirc} of Section \ref{SC:Notation}.
The definition of prototype sets yields
\[
\tau_{-2}[\mfF_1(-2)]
=
\iota[\mfF]\cap \tau_{-2}[\mfF].
\] 
Using \eqref{Eq:FormCircInv}, we can verify $\iota[\mfF]\cap \tau_{-2}[\mfF]=\tau_{-2}[\mfF\setminus\Dx(1;1)]$, then, $\mfF_1(-2)=\mfF\setminus \Dx(1;1)$.
Applying equations \eqref{Eq:FormRecInv} and \eqref{Eq:FormCircInv}, we obtain the set
\[
\iota[\mfF_1(-2)]
=
\iota[\mfF]\setminus \left\{ z\in \Cx: \real(z) > \frac{1}{2}\right\},
\]
which is depicted in Figure \ref{Fig:iotaF1-2}. Define
\[
\alpha\colon= \frac{2-\sqrt{3}}{2}
\;\text{ and }\;
\zeta_1 \colon= -\frac{1}{2} + i\alpha.
\]

Then, for all $m\in\Za$, $|m|\geq 2$, we have 
\[
\iota[\mfF_1(-2)]\cap \tau_{1+im}[\mfF]
=
\begin{cases}
\left[ \frac{-1-i}{2} + 1 + im,  \frac{-1+i}{2} +1 + im\right), &\text{ if } |m|\geq 3, \\[2ex]
\left[\, \overline{\zeta_1} + 1 + 2i,  \frac{-1+i}{2} +1 + 2i\right), &\text{ if } m=2, \\[2ex]
\left[ \frac{-1-i}{2} + 1 + 2i,  \zeta_1 + 1 - 2i \right), &\text{ if } m=-2,
\end{cases}
\]
which implies that
\[
\mfF_1(-2, 1+im)
=
\begin{cases}
\left[ \frac{-1-i}{2} ,  \frac{-1+i}{2}\right), &\text{ if } |m|\geq 3, \\[2ex]
\left[ \, \overline{\zeta_1} ,  \frac{-1+i}{2} \right], &\text{ if } m=2, \\[2ex]
\left[ \frac{-1-i}{2} ,  \zeta_1  \right), &\text{ if } m=-2.
\end{cases}
\]
Note that $\iota[\mfF_1(-2, 1+2i)]$ is the arch of the circle $C(-1;1)$ going from $-2+i+\alpha - \frac{1}{2}$ (included) to $-1-i$ (excluded), so
\[
\tau_{-2+i}[\mfF]
\cap
\iota[\mfF_2(-2, -1+2i)]
=
\left\{ \tau_{-2+i}\left(\alpha - \frac{i}{2}\right)\right\};
\]
hence,
\[
\mfF_2(-2, -1+2i, -2+i)
=
\left\{\alpha - \frac{i}{2}\right\}.
\]
To sum up, we have shown the next points:
\begin{enumerate}[i.]
\item $\clCc_1(-2)\neq \vac$;
\item If $m\in\Za$ and $|m|\geq 2$, then $\clC_2(-2, 1+im)\neq \vac$ but $\clCc_2(-2,1+im)=\vac$;
\item $\#\clC_3(-2, -1+2i, -2+i)=1$.
\end{enumerate}
\begin{rema01}
From the previous discussion, we also conclude that if $z\in\mfF\setminus\QU(i)$ satisfies $a_1(z)=-2$, then $\real(a_2(z))\leq 1$. Similar observations hold for other small digits; for example, $a_1(z)=1+i$ implies $\imag(z)\leq 0\leq \real(z)$.
\end{rema01}
In what follows, we denote
\[
\omfF:=\Cl(\mfF) 
\quad\text{ and }\quad
\mfFc:=\inte(\mfF),
\]
and, for all $n\in\Na$ and $\bfa\in\Omega(n)$, we write
\begin{align*}
\oclC_n(\bfa)\colon= \Cl(\clC_n(\bfa)), &\quad
\clCc_n(\bfa)\colon= \inte(\clC_n(\bfa)),\\
\omfF_n(\bfa)\colon= \Cl(\mfF_n(\bfa)), &\quad
\mfFc_n(\bfa)\colon= \inte(\mfF_n(\bfa)).
\end{align*}

\begin{def01}
For any $n\in\Na$, we say that $\bfa \in\scD^n$ is
\begin{enumerate}[\rm i.]
\item \textbf{valid} if $\clC_n(\bfa)\neq\vac$;
\item \textbf{invalid} if $\clC_n(\bfa)=\vac$;
\item \textbf{regular} if $\clCc_n(\bfa)\neq\vac$;
\item \textbf{irregular} if $\clC_n(\bfa)\neq\vac$ and $\clCc_n(\bfa)=\vac$;
\item \textbf{extremely irregular} if $\#\clC_n(\bfa)=1$;
\item \textbf{full} if $\mfF_n(\bfa)=\mfF$;
\item \textbf{almost full} if for every $b\in\scD$ the word $\bfa b$ is regular.
\end{enumerate}
We denote by $\Omega(n)$, $\sfR(n)$, $\sfIr(n)$, $\sfEI(n)$, $\sfF(n)$, and $\sfAF(n)$ to be the sets of valid, regular, irregular, extremely irregular, full, almost full words of length $n$ respectively. 

We extend these classifications to infinite sequences. We say that $\bfa=\sanu\in\scD^{\Na}$ is
\begin{enumerate}[\rm i.]
\item \textbf{valid} if it is the Hcf of some $z\in\mfF$;
\item \textbf{invalid} if it is not valid;
\item \textbf{regular} if it is valid and $(a_1,\ldots, a_n) \in\sfR(n)$ for all $n\in\Na$;
\item \textbf{irregular} if it is valid and $(a_1,\ldots, a_n)\in\sfIr(n)$ for some $n\in\Na$;
\item \textbf{extremely irregular} if it is valid and $(a_1,\ldots, a_n)\in\sfEI(n)$ for some $n\in\Na$;
\item \textbf{full} if it is valid and $(a_1,\ldots, a_n)\in\sfF(n)$ for all $n\in\Na$;
\item \textbf{almost full} if it is valid and $(a_1,\ldots, a_n)\in\sfAF(n)$ for all $n\in\Na$.
\end{enumerate}
We denote by $\Omega$, $\sfR$, $\sfEI$, $\sfF$, and $\sfAF$ the set of valid, regular, extremely irregular, full, and almost full words, respectively. We extend in an obvious way the notions of being full, almost full, regular, irregular, and extremely irregular to prototype sets.
\end{def01}

\begin{rema01}\label{Rema:AltDescriptFull}
This notation allows us to give an alternative description of prototype sets. For $n\in\Na$ and $\bfb=(b_1,\ldots, b_n)\in\scD^n$, we have
\[
\mfF_n(\bfb)\setminus \QU(i)
=
\{z \in\mfF\setminus\QU(i): (b_1,\ldots, b_n,a_1(z),a_2(z),\ldots)\in\Omega\}.
\]
\end{rema01}
The next proposition is a direct consequence of the definitions of cylinders and full words.

\begin{propo01}\label{Prop:FullCyl}
Let $n\in\Na$ and $\bfa=(a_1,\ldots, a_n)\in\scD^{n}$ be arbitrary, then
\[
\bfa\in \sfF(n)
\;\text{ if and only if }\;
\clC_n(\bfa)
=
\iota\tau_{a_1}\cdots \iota\tau_{a_n}[\mfF].
\]
\end{propo01}
A complex number $z=[a_1,a_2,\ldots]\in\mfF$ is \text{full} (resp. \textbf{almost full}, \textbf{regular}, \textbf{irregular}, \textbf{extremely irregular}) if $\sanu$ is full (resp. almost full, regular, irregular, extremely irregular). The number $\zeta_1$ discussed above is an example of an extremely irregular number. It is well established that Lebesgue almost all $z\in\mfF$ are regular. This follows from the existence of a $T$-ergodic Borel measure equivalent to the Lebesgue measure (see, for example, \cite[Proposition 1]{EiItoNak2019}) and Birkhoff's ergodic theorem. 



\begin{def01}
An \textbf{open cylinder} (resp. \textbf{open prototype set}, \textbf{closed cylinder}, \textbf{closed prototype set}) is a set of the form $\clCc_n(\bfa)$ (resp. $\mfFc_n(\bfa)$, $\oclC_n(\bfa)$, $\omfF_n(\bfa)$) for some $n\in\Na$ and $\bfa\in\sfR(n)$. 
\end{def01}
\begin{rema01}
There are only 13 closed prototype sets (see, for example, \cite[Section 2]{EiItoNak2019}). In fact, any closed prototype set is one of the sets depicted in Figure \ref{Fig:ClosedPrototypeSets} or one of their right-angled rotations. Note that the boundaries of $\omfF$, $\omfF_1(2)$, and $\omfF_1(1+i)$ consist of four smooth curves whereas the boundary of $\omfF_1(2+i)$ consists of five smooth curves. We refer to each of these curves as a \textbf{side} of the prototype set. For $n\in\Na$ and $\bfa\in \sfR(n)$, a \textbf{side} of $\oclC_n(\bfa)$ is the image under $Q_n(\,\cdot\,;\bfa)$ of a side of $\omfF_n(\bfa)$.
Analogous observations hold for open prototype sets. 
\begin{figure}[h!]
\begin{center}
\includegraphics[scale=0.5,  trim={5cm 9.0cm 5.0cm 8.0cm},clip]{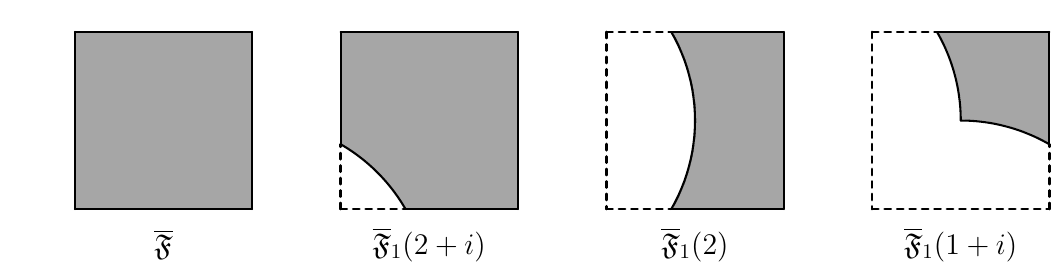}
\caption{Closed prototype sets \label{Fig:ClosedPrototypeSets}}
\end{center}
\end{figure}

\end{rema01}

There are sequences $\bfa\in\scD^{\Na}\setminus \Omega$ for which every prefix is valid. This is the case, for example, of the purely periodic sequence 
\[
\bfa
=
(-2,2i,2,-2i,-2,2i,2,-2i,\ldots).
\]
However, every sequence (finite or infinite) with sufficiently large terms is full. To be more precise, define 
\[
\clE
:=
\left\{ a\in \Za[i]: |a|\geq \sqrt{8}\right\}.
\]
{
The following proposition was mentioned (without proof) in \cite{GeroPhDThesis}. We provide its proof for the reader's convenience. 
}
\begin{propo01}\label{PROP:SQRT8-REG:NUVO}
We have $\clE^{n}\subseteq \sfF(n)$ for all $n\in\Na$ and $\clE^{\Na}\subseteq \sfF$.
\end{propo01}
\begin{proof}
First, we show inductively the assertion for finite words. Take any $a_1\in\clE$. We obtain $\tau_{a_1}[\mfF]\subseteq \iota[\mfF]$ from Figure \ref{Fig:02}. Proposition \ref{Prop:ObsCyls}.\ref{Prop:ObsCyls:01} then gives $\clC_1(a_1) = \iota\tau_{a_1}[\mfF]$ and, by Proposition \ref{Prop:FullCyl}, $a_1\in\sfF(1)$. Assume $\clE^{N}\subseteq \sfF(N)$ for some $N\in\Na$ and take $\bfa=(a_1,\ldots, a_N, a_{N+1})\in\clE^{N+1}$. In view of Proposition \ref{Prop:ObsCyls}.\ref{Prop:ObsCyls:02}, we have
\begin{align*}
\clC_{N+1}(a_1,\ldots, a_{N+1})
&=
\clC_1(a_1)\cap \iota\tau_{a_1}[\clC_{N}(a_2,\ldots, a_{N+1})] \\
&=
\iota\tau_{a_1}[\mfF\cap \clC_{N}(a_2,\ldots, a_{N+1})] &&(\text{by the case } N=1) \\
&=
\iota\tau_{a_1}[\clC_{N}(a_2,\ldots, a_{N+1})]\\ & =
\iota\tau_{a_1}\iota\tau_{a_2}\cdots \iota\tau_{a_{N+1}} [\mfF].
\end{align*}
The last equality follows from the induction hypothesis and Proposition  \ref{Prop:FullCyl}. Lastly, by Proposition \ref{Prop:FullCyl} again, we conclude that $\bfa\in \sfF(N+1)$ and $\clE^{N+1}\subseteq \sfF(N+1)$ follows.

Now we turn to the second assertion.  Consider any $\bfa=\sanu\in \clE^{\Na}$. By the previous part, we have $\clCc_n(a_1,\ldots, a_n)\neq \vac$ for all $n\in\Na$. Let $\epsilon$ denote the empty word and define $\clC_0(\epsilon)=\mfF$. We claim that
\[
\clCc_{n+1}(a_1,\ldots, a_{n+1})
\subseteq
\oclC_n(a_1,\ldots, a_n)
\;\text{ for all } n\in\Na_{0}.
\]
The case $n=0$ follows from Figure \ref{Fig:02} while for an arbitrary $n\in\Na$ we have
\begin{align*}
\oclC_{n+1}(a_1,\ldots, a_{n+1})
&= 
\iota \tau_{a_1} \cdots \iota\tau_{a_n} \iota\tau_{a_{n+1}} [\omfF]\\
&=
\iota \tau_{a_1}\cdots \iota\tau_{a_n} [\oclC_1(a_{n+1})] \\
&\subseteq
\iota \tau_{a_1}\cdots \iota\tau_{a_n} [\mfFc]
\\ &=
\clCc_{n}(a_1,\ldots, a_{n}),
\end{align*}
thus proving the claim.
Parts \ref{Propo3.2_ii} and \ref{Propo3.2_iv} of Proposition \ref{Propo3.2} imply $|\clC_n(a_1,\ldots, a_n)|\to 0$ as $n\to\infty$. Therefore, there is a single element $z$ in 
\[
\bigcap_{n\in \Na} \clC_n(a_1,\ldots, a_n),
\]
which clearly satisfies $z=[a_1,a_2,\ldots]$ and, thus, $\bfa\in \sfF$.
\end{proof}
 
\begin{propo01}\label{Prop:ClosureOfInterior}
If $n\in\Na$ and $\bfa\in\sfR(n)$, then $\omfF_n(\bfa)=\Cl\left( \mfFc_n(\bfa)\right)$ and $\oclC_n(\bfa)=\Cl\left( \clCc_n(\bfa)\right)$.
\end{propo01}
\begin{proof}
The result is obvious for $n=1$ (see Figure \ref{Fig:02}). For the general case, use induction on~$n$.
\end{proof}
For all $z\in\Cx$, write $\Pm(z) \colon= \min\{ | \real(z) |,  | \imag(z) |\}$.

\begin{propo01}\label{Prop:PmandRn}
Let $n$ be a natural number and $\bfa\in \sfR(n)$.
\begin{enumerate}[\rm i.]
\item \label{Prop:PmandRn:i} If $\Pm(a_n)\geq 3$, then $\bfa\in\sfF(n)$.
\item \label{Prop:PmandRn:ii} If $\|a_n\|\geq 3$, then $\bfa\in\sfAF(n)$.
\item \label{Prop:PmandRn:iii} There exists $d\in\scD$ such that $\Pm(d)=\|d\|=3$ and $\bfa d\in\sfF(n+1)$.\end{enumerate}
\end{propo01}
\begin{proof}
Parts \ref{Prop:PmandRn:i} and \ref{Prop:PmandRn:ii} are shown by checking case by case all the open prototype sets. Part \ref{Prop:PmandRn:iii} is essentially \cite[Lemma 2.6(b)]{HeXio2021-01}.
\end{proof}
\subsection{Symmetries of regular cylinders}\label{Subsection:Symmetries oF Regular Cylinders}
The symmetries of the Hcfs process are at the core of the geometric properties that
we study. For instance, we extend the symmetries exhibited by open cylinders of level 1 in Figure \ref{Fig:mfFPart} to open cylinders of arbitrary level. With this in mind, consider the maps $\Rota, \Mir_1:\overline{\Cx}\to\overline{\Cx}$ given by
\[
\forall z\in \overline{\Cx}
\quad
\Rota(z)\colon=iz
\quad\text{ and }\quad
\Mir_1(z)\colon= \overline{z}
\]
and call $\Di_8$ the group of isometries generated by $\Rota$ and $\Mir_1$ along with the composition of functions. Clearly, $\Di_8$ is the dihedral group of order $8$. In what follows, the juxtaposition of elements of $\Di_8$, $\iota$, and the maps $\tau_a$ stands for their composition; for example, $\Rota\,\Mir_1(z)=i\overline{z}$ for $z\in \Cx$. 

We will denote the cyclic group generated by $\Rota$ by $\RotaG$. If $A,B\subseteq \Cx$, we write $A\equiv B\mod \RotaG$ whenever there is some $f\in \RotaG$ such that $f[A]=B$. We likewise define $A\equiv B\mod\Di_8$.

For reference, we gather some straightforward observations in the next proposition. 
\begin{propo01}\label{PROPO:SIMET:CILINDROS:01}
The following assertions hold:
\begin{enumerate}[\rm i.]
\item \label{PROPO:SIMET:CILINDROS:01:01} Every $f\in \Di_8$ gives a bijection from $\scD$ onto itself and $\|f(a)\|=\|a\|$ for all $a\in\scD$. 

\item \label{PROPO:SIMET:CILINDROS:01:02} $\Rota\, \iota=\iota\,\Rota^{-1}$, $\Rota^{-1}\,\iota=\iota\,\Rota$.

\item \label{PROPO:SIMET:CILINDROS:01:03} For all $a\in\scD$ we have $\Rota\,\tau_a=\tau_{\Rota(a)}\,\Rota$, $\Rota^{-1} \tau_a=\tau_{\Rota^{-1}(a)}\, \Rota^{-1}$.
\item \label{PROPO:SIMET:CILINDROS:01:04} For all $a\in\scD$ we have $\Rota\,\iota\,\tau_{a} = \iota\,\tau_{\Rota^{-1}(a)}\Rota^{-1}$ and $\Rota^{-1}\,\iota\,\tau_{a} = \iota\,\tau_{\Rota(a)}\,\Rota$.
\end{enumerate}
\end{propo01}

\begin{propo01}\label{PROPO:SIMET:CILINDROS:02}
For all $n\in\Na$, $\bfa=(a_1,\ldots, a_n)\in\sfR(n)$, and $f\in\Di_8$, there are some $\bfb=(b_1,\ldots,b_n)\in \sfR(n)$ and $f_1,\ldots,f_n\in \Di_8$ such that 
\[
f\left[\clC^{\circ}_n(\bfa)\right]
= \clC^{\circ}_n(\bfb)
\quad\text{ and }\quad
b_j=f_j(a_j) \;\text{ for all }\; j\in\{1,\ldots,n\}.
\]
We can replace $\clC^{\circ}_n$ with $\oclC_n$ in the previous equality.
\end{propo01}
\begin{proof}
It is enough to show the result when $f$ is one of the generators $\Rota$ or $\Mir_1$. Assume that $f=\Rota$ and let us demonstrate inductively that for all $n\in\Na$ and $\bfa=(a_1,\ldots,a_n)\in\sfR(n)$ we have
\begin{equation}\label{EQ:ROTA:CIL}
\Rota\left[ \clC^{\circ}_n(\bfa)\right]
=
\clCc_n\left( \Rota^{-1}(a_1), \Rota(a_2), \Rota^{-1}(a_3),\ldots,\Rota^{(-1)^n}(a_n)\right).
\end{equation}
Take any $a_1\in\scD$. First, we show that 
\begin{equation}\label{EQ:PROPO:SIMET:CILINDROS:ind:01}
\Rota[\clCc_1(a_1)]
=
\clCc_1\left(\Rota^{-1}(a_1)\right)
\;\text{ and }\;
\Rota^{-1}\left[\clCc_1\left(a_1 \right)\right]
=
\clCc_1\left(\Rota(a_1)\right).
\end{equation}
To this end, note that a complex number $z\in\mfFc$ belongs to $\clCc_1(a_1)$ if and only if there is some $w\in\mfFc_1(a_1)$ such that
\[
z=\frac{1}{a_1 + w}.
\]
This is equivalent to
\[
\Rota(z)=\frac{1}{\Rota^{-1}(a_1) + \Rota^{-1}(w)},
\]
or, equivalently, $\iota\, \Rota(z) = \Rota^{-1}(a_1) + \Rota^{-1}(w)$, and since $\Rota(z)\in \mfFc$, $\Rota^{-1}(a_1)\in\scD$, and $\Rota^{-1}(w)\in\mfF$, we have
\[
a_1\left(\Rota(z)\right) = \Rota^{-1}(a_1). 
\]
We may, thus, conclude
\[
\Rota\left[ \clCc_1(a_1)\right]
\subseteq 
\clC_1\left(\Rota^{-1}(a_1)\right).
\]
Moreover, because $\Rota$ maps open sets onto open sets, we even have
\[
\Rota\left[ \clCc_1(a_1)\right]
\subseteq 
\clCc_1\left(\Rota^{-1}(a_1)\right).
\]
In a similar fashion, we may show $\Rota^{-1}[\clCc_1\left(\Rota^{-1}(a_1)\right)] \subseteq \clCc_1(a_1)$ and we get the first identity in \eqref{EQ:PROPO:SIMET:CILINDROS:ind:01}. The second equality follows from the first one. 

Assume \eqref{EQ:ROTA:CIL} is true for some $n=N\in\Na$. Take any $\bfa=(a_1,\ldots, a_N,a_{N+1}) \in\sfR(N+1)$ and let us show that
\begin{equation}\label{EQ:PROPO:SIMET:CILINDROS:02}
\Rota\left[ \clC^{\circ}_{N+1}(\bfa)\right]
=
\clC^{\circ}_{N+1} \left( \Rota^{-1}(a_1), \Rota(a_2), \Rota^{-1}(a_3),\ldots,\Rota^{(-1)^{N+1}}(a_{N+1})\right).
\end{equation}
Proposition \ref{Prop:ObsCyls}.\ref{Prop:ObsCyls:02} implies that
\begin{equation}\label{EQ:PROPO:SIMET:CILINDROS:03}
\Rota\left[ \clC_{N+1}^{\circ}(a_1, \ldots, a_{N+1})\right]
=
\Rota\left[ \clC_{N}^{\circ}(a_1,\ldots, a_N )\right]
\cap \Rota \left[ Q_N(\clCc_1(a_{N+1});(a_1,\ldots, a_N) \right].
\end{equation}
We apply the induction hypothesis on the first term to the right in order to get
\begin{equation}\label{EQ:PROPO:SIMET:CILINDROS:04}
\Rota\left[ \clC_{N}^{\circ}(a_1,\ldots, a_N)\right]
=
\clC_N^{\circ}\left( \Rota^{-1}(a_1),\ldots, \Rota^{(-1)^{N}}(a_N)\right).
\end{equation}
For the second term on the right-hand side of \eqref{EQ:PROPO:SIMET:CILINDROS:03}, observe that if $N=2$ then
\[
\Rota\,\iota\,\tau_{a_1}\,\iota\,\tau_{a_2}
=
(\iota \tau_{\Rota^{-1}(a_1)}\,\Rota^{-1})\iota\,\tau_{a_2}
=
\iota \,  \tau_{\Rota^{-1} (a_1)}\, \iota\, \tau_{\Rota(a_2)}\,\Rota
\]
and $\Rota\circ Q(\,\cdot\, ;  (a_1,a_2)) = Q(\,\cdot\, ; (\Rota^{-1}(a_1),\Rota(a_2)) \circ \Rota$. In general, from Proposition \ref{PROPO:SIMET:CILINDROS:01}.\ref{PROPO:SIMET:CILINDROS:01:04}, we obtain
\[
\Rota
\circ
Q(\,\cdot\, ;  \bfa)
=
Q(\,\cdot\, ; (\Rota^{-1}(a_1),\Rota(a_2),\ldots, \Rota^{(-1)^N}(a_2))
\circ
\Rota^{(-1)^N}.
\]

Thus, depending on the parity of $N$, we use one of the identities in \eqref{EQ:PROPO:SIMET:CILINDROS:ind:01} to conclude
\[
Q_N (\clCc_1(a_{N+1}) ;  (a_1,\ldots, a_N) ) 
=
Q_N\left( \clC_1^{\circ}\left(\Rota^{(-1)^{N+1}}(a_{N+1} )\right); \left( \Rota^{-1}(a_1),\ldots,\Rota^{(-1)^N}(a_N) \right)\right).
\]
We now apply Proposition \ref{PROPO:SIMET:CILINDROS:01}.\ref{PROPO:SIMET:CILINDROS:01:04} to obtain \eqref{EQ:PROPO:SIMET:CILINDROS:02}.

Using a similar yet simpler inductive argument, we may show that for all $n\in\Na$ and $ \bfa=(a_1,\ldots,a_n)\in\sfR(n)$, we have
\begin{equation}\label{EQ:MIR1:CIL}
\Mir_1\left[\clC^{\circ}_n(\bfa)\right]
=
\clC^{\circ}_n\left(\Mir_1(a_1), \Mir_1(a_2), \ldots, \Mir_1(a_n)\right).
\end{equation}

The assertion for closed cylinders follows from the proof above, Proposition \ref{Prop:ClosureOfInterior}, and the fact that $\Mir_1$ and $\Rota$ are homeomorphisms from $\overline{\Cx}$ onto itself.
\end{proof}

Define $\Mir_2\colon=\Rota^2\,\Mir_1$, so $\Mir_2(z)=-\overline{z}$ for $z\in \Cx$. When $f=\Mir_2$, we can use equations \eqref{EQ:MIR1:CIL} and \eqref{EQ:ROTA:CIL} to determine the functions $f_1,\ldots, f_n$ appearing in Proposition \ref{PROPO:SIMET:CILINDROS:02}.

\begin{coro01}\label{CORO:SIMET:CILINDROS:02}
For all $n\in\Na$ and $\bfa\in\sfR(n)$, we have
\begin{equation}\label{EQ:MIR2:CIL}
\Mir_2\left( \clC^{\circ}_n(\bfa)\right)
=
\clC^{\circ}_n\left(\Mir_2(a_1), \Mir_2(a_2),\ldots, \Mir_2(a_n)\right).
\end{equation}
\end{coro01}

\subsection{Size of regular cylinders}
In this subsection, we recall some estimates of the diameter and the Lebesgue measure $\leb$ of regular cylinders in terms of the sequence $\seqn$. See item \ref{Not:ConcatWords} in Section \ref{SC:Notation} for the definition of $\bfa\bfb$.
\begin{propo01}\label{Prop:SizeEstimates}
Let $m,n\in\Na$, $\bfa\in\Omega(n)$, and $\bfb\in \Omega(m)$ be arbitrary.
\begin{enumerate}[\rm i.]
\item\label{Prop:SizeEstimates:i}  \cite[Lemma 2.7]{HeXio2021-01}
There is an absolute constant $0<\kappa_1<1$ such that if $\bfa\in\sfR(n)$, then
\[
\frac{\kappa_1}{|q_n(\bfa)|^2} \leq |\clC_n(\bfa)| \leq \frac{2}{|q_n(\bfa)|^2} 
\]
and
\[
\frac{\pi\kappa_1}{|q_n(\bfa)|^4} \leq \leb(\clC_n(\bfa)) \leq \frac{\pi}{|q_n(\bfa)|^4}.
\]
\item\label{Prop:SizeEstimates:ii} If $\bfa\bfb\in \Omega(m+n)$, then
\[
\frac{1}{5\sqrt{2}}\|q_n(\bfa)\|\|q_m(\bfb)\|
\leq 
\|q_{n+m}(\bfa\bfb)\|
\leq 
6\|q_n(\bfa)\|\|q_m(\bfb)\|.
\]
\item\label{Prop:SizeEstimates:iii} There is an absolute constant $\kappa_2>1$ such that if $\bfa$, $\bfb$, and $\bfa\bfb$ are regular, then
\[
\frac{1}{\kappa_2}|\clC_{n}(\bfa)||\clC_{m}(\bfb)|
\leq 
|\clC_{n+m}(\bfa\bfb)|
\leq 
\kappa_2|\clC_{n}(\bfa)||\clC_{m}(\bfb)|.
\]
\item \label{Prop:SizeEstimates:iv} There is an absolute constant $\kappa_3>1$ such that for any $n\in\Na$, $\bfa\in\Omega(n)$, and $A \subseteq \omfF_n(\bfa)$ we have
\[
\frac{1}{\kappa_3}\left| \clC_n(\bfa)\right| |A|
<
 \left| Q_n(A;\bfa) \right| 
<
\kappa_3\left| \clC_n(\bfa)\right| |A|.
\]
\end{enumerate}
\end{propo01}
\begin{proof}
For part \ref{Prop:SizeEstimates:ii}, we know that $\bfa\bfb\in \Omega(n+m)$ implies
\begin{equation}\label{Eq:Lem513:Pf}
\frac{1}{5}|q_n(\bfa)||q_m(\bfb)|
\leq 
|q_{n+m}(\bfa\bfb)|
\leq 
3|q_n(\bfa)||q_m(\bfb)|
\end{equation}
(for a proof, see \cite[Lemma 2.2.(g)]{HeXio2021-01}). Since $\|z\|\leq |z|\leq \sqrt{2}\|z\|$ for all $z\in\Cx$, we have
\[
\|q_{n+m}(\bfa\bfb)\| 
\leq 
|q_{n+m}(\bfa\bfb)| 
\leq 
3|q_n(\bfa)||q_m(\bfb)| 
\leq 
6\|q_n(\bfa)\|\|q_m(\bfb)\|
\]
and, similarly,
\[
\|q_{n+m}(\bfa\bfb)\| \geq \frac{1}{\sqrt{2}} |q_{n+m}(\bfa\bfb)| 
\geq 
\frac{1}{5\sqrt{2}} \|q_n(\bfa)\|\|q_m(\bfb)\|.
\]
For part \ref{Prop:SizeEstimates:iii}, combine parts \ref{Prop:SizeEstimates:i} and \ref{Prop:SizeEstimates:ii}.

Part \ref{Prop:SizeEstimates:iv} is a consequence of Proposition \ref{Prop:QnCorrespRule} and the next inequality, which follows from Proposition \ref{Propo3.2}.\ref{Propo3.2_i} and $z \in \mfF$,
\begin{equation}\label{Eq:BndOnDenomQn}
1 - \frac{1}{\sqrt{2}}
\leq 
\left| 1 + \frac{q_{n-1}}{q_{n}}z\right|
\leq 
1+\frac{1}{\sqrt{2}}.
\end{equation}
A similar argument is developed in full detail in the proof of Proposition \ref{PROPO:GC01:QnDISCOS}.
\end{proof}

\section{Geometric constructions}\label{Sec:GC}

In this section, we count how many sets of a particular kind are intersected by a disc of small radius and center in $\omfF$. To this end, we consider two possible geometric configurations, which we call \textit{constructions}. The results in this section will be used for the proof of the lower bound in Theorem \ref{TEO:MAIN}.

We start off by proving that the corners of $\mfF$ cannot be attained using closed cylinders of level at least $2$ (Proposition \ref{GC:PROPO:02}). Let $\Esq$ be the set of corners\footnote{The notation $\Esq$ comes from the word \textit{esquina}, which is Spanish for \textit{corner}.}  of $\mfF$; that is,
\[
\Esq
\colon=
\left\{ \frac{1+i}{2}, \frac{-1+i}{2}, \frac{-1-i}{2}, \frac{1-i}{2}\right\}.
\]
Note that
\begin{equation}\label{Eq:invEsq}
\iota[\Esq]
\colon=
\left\{  1 - i , -1 - i, -1 + i,  1 + i \right\}
\subseteq
\scD.
\end{equation}
We will write
\[
\Di_8(\zeta_1)
\colon=
\{f(\zeta_1):f\in \Di_8\}
\;\text{ and }\;
\RotaG(\zeta_1)
\colon=
\{f(\zeta_1):f\in \RotaG\}.
\]
\begin{propo01}\label{GC:PROPO:01}
For all $n\in\Na$ and $\bfa=(a_1,\ldots, a_n)\in \sfR(n)$, we have
\[
\bigcup_{ \substack{d\in \scD \\ \bfa d\in \sfR(n+1)} } \oclC_{n+1}(\bfa d)
=
\oclC_n(\bfa) \setminus\{ [a_1,\ldots,a_n]\}.
\]
\end{propo01}
\begin{proof}
From Figure \ref{Fig:mfFPart}, we immediately conclude that
\begin{equation}\label{GC:EQ:01}
\bigcup_{ d\in \scD } \oclC_{1}(d)
=
\omfF\setminus\{ 0\}.
\end{equation}
We, however, provide a full proof to emphasize the clarity we gain by relying on pictures. First, for each $d\in\scD$ and $z\in \clC_1(d)$, there is some $w\in\mfF$ such that $z= \frac{1}{d+w}$, so
\[
|z| \geq \frac{1}{|d|+|w|} \geq \frac{1}{|d| + \frac{\sqrt{2}}{2}},
\;\text{ and }\;
0\not\in \oclC_1(d).
\]
This shows the inclusion 
\[
\bigcup_{ d\in \scD } \oclC_{1}(d)
\subseteq
\omfF\setminus\{ 0\}.
\]
For the other inclusion, take any $z\in \omfF\setminus\{0\}$. If $z\in \mfF$, then $z\in \clC_1(a_1(z))\subseteq \oclC_1(a_1(z))$. If $z\in\omfF\setminus\mfF$, then $z$ belongs to either the upper or the right side of $\mfF$. Hence, there is some $\Mir\in \{\Mir_1,\Mir_2\}$ such that $\Mir(z)\in \mfF$, so
\[
\Mir(z) \in \oclC_1(a_1(\Mir(z)).
\]
Then, in view of \eqref{EQ:MIR1:CIL} or \eqref{EQ:MIR2:CIL}, we have
\[
z \in  \Mir ( \oclC_1(a_1(\Mir(z)) 
=
\oclC_1( \Mir (a_1(\Mir(z))).
\]
For the general case, take any $n\in\Na$ and $\bfa\in\sfR(n)$. Equation \eqref{GC:EQ:01} implies that
\[
\omfF_n(\bfa) \cap \bigcup_{ d\in \scD } \oclC_{1}(d)
=
\omfF_n(\bfa)\setminus\{0\}.
\]
Lastly, apply $Q_n(\,\cdot\,;\bfa)$ on both sides of the previous equality and conclude.
\end{proof}

\begin{propo01}\label{GC:PROPO:02}
For each $n\in\Na_{\geq 2}$, we have
\[
\Esq
\cap
\bigcup_{\bfa\in\sfR(n)} \oclC_n(\bfa)
=\vac.
\]
\end{propo01}
\begin{proof}
Let $n\in\Na_{\geq 2}$ and $\bfa=(a_1,\ldots, a_n)\in\sfR(n)$ be arbitrary. By Proposition \ref{Prop:ObsCyls}.\ref{Prop:ObsCyls:02}, 
\begin{align*}
\iota[\oclC_n(\bfa)]
&=
\iota\left[\oclC_1(a_1)\cap \iota\tau_{a_1}[\oclC_{n-1}(a_2,\ldots, a_n)]\right] \\
&=
\iota[\oclC_1(a_1)]\cap \tau_{a_1}[\oclC_{n-1}(a_2,\ldots, a_n)] \\
&=
\tau_{a_1}[\omfF_1(a_1)]\cap \tau_{a_1}[\oclC_{n-1}(a_2,\ldots, a_n)] \\
&=
\tau_{a_1}[\omfF_1(a_1)\cap \oclC_{n-1}(a_2,\ldots, a_n)].
\end{align*}
Form \eqref{Eq:invEsq}, we obtain $\iota[\Esq]\cap \tau_{a_1}[\omfF]\subseteq \{a_1\}$, so
\begin{align*}
\iota[\Esq\cap \oclC_n(\bfa)]
&=
\iota[\Esq]\cap \iota[\oclC_n(\bfa)] \\
&\subseteq 
\iota[\Esq]\cap \tau_{a_1}[\omfF_1(a_1)\cap \oclC_{n-1}(a_2,\ldots, a_n)] \\
&\subseteq 
\iota[\Esq]\cap \tau_{a_1}[\omfF \cap \oclC_{n-1}(a_2,\ldots, a_n)] \\
&\subseteq
\{a_1\} \cap \tau_{a_1}[\omfF_1(a_1)\cap \oclC_{n-1}(a_2,\ldots, a_n)] \\
&=
\tau_{a_1}\left[ \{0\} \cap \omfF_1(a_1)\cap \oclC_{n-1}(a_2,\ldots, a_n) \right] \\
&\subseteq 
\tau_{a_1} \left[ \{ 0 \} \cap \oclC_{n-1}(a_2,\ldots, a_n) \right] \\
&=\vac. &&\text{(by Proposition \ref{GC:EQ:01})}
\end{align*}
\end{proof}

\subsection{First geometric construction}
The main result of this subsection is Theorem \ref{TEO:GC:01}. It tells us that any disc of sufficiently small radius and center in a closed cylinder of level $n$ intersects at most six regular cylinders of level $n$. More precisely, given $n\in \Na$ and $\bfa\in\sfR(n)$, the set of \textbf{neighbors} of $\bfa$ is\footnote{The notation $\Veci_n(\bfa)$ comes from the Spanish word \textit{vecino}, which means \textit{neighbor}.}
\[
\Veci_n(\bfa) 
\colon= 
\left\{ \bfb\in \sfR(n): \bfb\neq \bfa \text{ and } \oclC_n(\bfa) \cap \oclC_n(\bfb) \neq\vac\right\}.
\]
For $\rho>0$ and $z\in\omfF$, define 
\[
G_n(z,\bfa,\rho)
\colon=
\{ \bfb\in \sfR(n): \Dx(z;\rho|\clC_n(\bfa)|) \cap \oclC_n(\bfb)\neq\vac\}.
\]
\begin{teo01}\label{TEO:GC:01}
There is a constant $\rho>0$ such that for every $n\in\Na$, $\bfa\in \sfR(n)$ and $z\in \oclC_n(\bfa)$ we have
\[
\#G_n(z,\bfa, \rho )\leq 6, 
\quad
G_n(z,\bfa, \rho )\subseteq \{\bfa\}\cup \Veci_n(\bfa),
\]
and
\[
\Dx\left(z; \rho  |\clC_n(\bfa)|\right)\cap \omfF
\subseteq 
\bigcup_{\bfb\in G_n(z,\bfa, \rho )} \oclC_n(\bfb).
\]
\end{teo01}

We start with an auxiliary result. 
It is well known that M\"obius transformations (or fractional linear transformations) map discs in $\overline{\Cx}$ onto discs in $\overline{\Cx}$. Moreover, if $t$ is a M\"obius transformation and $D=\Dx(z;r)$ for some $z\in\Cx$ and $r>0$, then $t[D]$ is also a disc in $\Cx$ whenever the pole of $t$ is not contained in $D$. In this case, the center of $t[D]$ might not be $t(z)$, but there is some $\rho>0$ such that
\[
\Dx\left( t(z);\rho\right)
\subseteq
t\left[ \Dx(z;r)\right].
\]
Proposition \ref{PROPO:GC01:QnDISCOS} below provides some values for $\rho$ when $t=Q_n(\,\cdot \,;\bfa)$ for some $\bfa\in\sfR(n)$, $n\in\Na$.

Note that, by the triangle inequality, for all $w_1\in\Cx$, $\rho>0$, and $w_2\in \Dx(w_1;\rho)$ we have
\begin{equation}\label{Eq:DscInDsc}
\Dx\left( w_2;\rho - |w_1 - w_2|\right)
\subseteq 
\Dx\left( w_1; \rho\right).
\end{equation}
Furthermore, the next equality holds:
\[
\rho - |w_1-w_2|
=
\inf_{0\leq \theta \leq 2\pi} \left|w_2 - (w_1 + \rho e^{i\theta}) \right|.
\]
\begin{propo01}\label{PROPO:GC01:QnDISCOS}
Let $n\in\Na$ and $\bfa\in\sfR(n)$ be arbitrary. If $w\in\omfF$ and $0<r<1$ are such that $\Dx(w;r)\subseteq \Dx (0; 1)$, then $Q_n(\Dx(w;r);\bfa)$ is a disc in the complex plane and
\[
\Dx\left( Q_n(w;\bfa) ; \frac{r}{2\left( 1 + \frac{\sqrt{2}}{2}\right)\left( 2 + \frac{\sqrt{2}}{2}\right)} |\clC_n(\bfa)|\right)
\subseteq
Q_n(\Dx(w;r);\bfa).
\]
\end{propo01}
\begin{proof}
The set $Q_n(\Dx(w;r);\bfa)$ is a disc on the Riemann sphere. 
Writing $q_n=q_n(a_1,\ldots, a_n)$ and $q_{n-1} = q_{n-1}(a_1,\ldots, a_{n-1})$, Proposition \ref{Propo3.2}.\ref{Propo3.2_i} and Proposition \ref{Prop:QnCorrespRule} tell us that $\overline{\Dx}(w;r)$ does not contain the pole of $Q_n(\,\cdot\,;\bfa)$. Then, there are some $w'\in\Cx$ and $\rho>0$ such that 
\[
Q_n(\Dx(w;r);\bfa)
=
\Dx(w';\rho).
\]
Next, we show the inclusion in the statement. For any $\theta\in [0,2\pi)$, we have
\begin{align*}
\left|  Q_n(w ; \bfa) - Q_n(w + re^{i\theta}; \bfa) \right| 
&=
\left|
\frac{w}{q_n^2 \left( 1+ \frac{q_{n-1}}{q_n} w\right)}
 -
\frac{w + re^{i\theta}}{q_n^2 \left( 1+ \frac{q_{n-1}}{q_n} (w + re^{i\theta})\right)}
\right| 
\quad \text{(by Proposition \ref{Prop:QnCorrespRule})}
\\
&=
\frac{r}{|q_n|^2\left| 1+ \frac{q_{n-1}}{q_n}w\right|\left| 1+ \frac{q_{n-1}}{q_n} (w + re^{i\theta}) \right| } \\
&\geq 
\frac{r}{2\left( 1 + \frac{\sqrt{2}}{2} \right)\left(2 + \frac{\sqrt{2}}{2} \right)}|\clC_n(\bfa)|. 
\qquad \text{(by Proposition \ref{Prop:SizeEstimates}.\ref{Prop:SizeEstimates:i} and }w\in \omfF)
\end{align*}
As a consequence, since $Q_n(C(w;r);\bfa)=C(w';\rho)$, we have
\[
\inf_{0\leq \theta < 2\pi}
\left|Q_n(w;\bfa) -  (w' + \rho e^{i\theta}) \right|
\geq 
\frac{r}{2\left( 1 + \frac{\sqrt{2}}{2} \right)\left(2 + \frac{\sqrt{2}}{2} \right)}|\clC_n(\bfa)|.
\]
Lastly, in view of $Q_n(w;\bfa)\in \Dx(w';\rho)$, we appeal to \eqref{Eq:DscInDsc} to conclude the result.
\end{proof}
\subsubsection{Proof of Theorem \ref{TEO:GC:01}}
We split the proof of Theorem \ref{TEO:GC:01} into five parts (lemmas \ref{PROPO:GC01:Rho:01} through \ref{PROPO:GC01:Rho:06}) depending on the position of $z$ with respect to the boundary of $\clC_n(\bfa)$. We denote the boundary of a given $A\subseteq\Cx$ by $\Fr(A)$. The constants $\rho_1, \ldots, \rho_5$ in the lemmas below are absolute. 

\begin{lem01}\label{PROPO:GC01:Rho:01}
There is a constant $\rho_{1}>0$ such that for all $f\in\Di_8$, $n\in\Na$, and $\bfa\in\sfR(n)$ for which $\xi:=f(\zeta_1) \in \oclC_n(\bfa)\cap \Fr(\mfF)$, we have 
\[
G_n(\xi,\bfa,\rho_1)\subseteq \Veci_n(\bfa)\cup \{\bfa\},
\quad
\#G_n(\xi,\bfa,\rho_1)\leq 3, 
\]
and
\[
\Dx(\xi; \rho_{1}|\clC_n(\bfa)|)\cap \omfF
\subseteq 
\bigcup_{\bfb\in G_n(\xi,\bfa,\rho_1)} \oclC_n(\bfb).
\]
\end{lem01}

\begin{lem01}\label{PROPO:GC01:Rho:02}
Take any $n\in\Na_{\geq 2}$, $\bfa=(a_1,\ldots, a_n) \in \sfR(n)$ such that $\oclC_n(\bfa)\cap \Fr(\mfF)\neq\vac$, and  $z\in \Fr(\mfF)$ such that
\begin{enumerate}[\rm i.]
\item There exists $\bfc\in \sfR(n)$ satisfying $\bfa\neq\bfc$ and $z\in \oclC_n(\bfa)\cap \oclC_n(\bfc)$;
\item If $z\in \oclC_1(a)$, then $a=a_1$.
\end{enumerate}
There is a constant $\rho_2>0$ for which
\[
G_n(z,\bfa,\rho_2)\subseteq \Veci_n(\bfa)\cup \{\bfa\},
\quad
\#G_n(z,\bfa,\rho_2)\leq 3, 
\]
and
\[
\Dx(z; \rho_{2}|\clC_n(\bfa)|)\cap \omfF\subseteq \bigcup_{\bfb\in G_n(z,\bfa,\rho_2)} \oclC_n(\bfb).
\]
\end{lem01}

\begin{lem01}\label{PROPO:GC01:Rho:03}
There is a constant $\rho_{3}>0$ such that for any $n\in\Na$, $\bfa\in\sfR(n)$ with $\Fr(\mfF)\cap \oclC_n(\bfa)\neq\vac$, and any $z\in\Fr(\mfF)\cap \oclC_n(\bfa)$, we have
\[
G_n(z,\bfa,\rho_3)\subseteq \Veci_n(\bfa)\cup \{\bfa\},
\quad
\#G_n(z,\bfa,\rho_3)\leq 3, 
\]
and
\[
\Dx(z; \rho_{3}|\clC_n(\bfa)|)\cap \omfF\subseteq \bigcup_{\bfb\in G_n(z,\bfa,\rho_3)} \oclC_n(\bfb).
\]
\end{lem01}

\begin{lem01}\label{PROPO:GC01:Rho:05}
There is a constant $\rho_{4}>0$ such that for all $n\in\Na$, $\bfa\in\sfR(n)$, and $z\in \mfFc \cap \Fr(\clC_n(\bfa))$, we have
\[
G_n(z,\bfa,\rho_4)\subseteq \Veci_n(\bfa)\cup \{\bfa\},
\quad
\#G_n(z,\bfa,\rho_4)\leq 6, 
\]
and
\[
\Dx(z; \rho_{4}|\clC_n(\bfa)|)\cap \omfF\subseteq \bigcup_{\bfb\in G_n(z,\bfa,\rho_4)} \oclC_n(\bfb).
\]
\end{lem01}
\begin{lem01}\label{PROPO:GC01:Rho:06}
There is a constant $\rho_{5}>0$ such that if $n\in\Na$, $\bfa\in\sfR(n)$, and $z\in \clCc_n(\bfa)$, then
\[
G_n(z,\bfa,\rho_5)\subseteq \Veci_n(\bfa)\cup \{\bfa\},
\quad
\#G_n(z,\bfa,\rho_5)\leq 6, 
\]
and
\[
\Dx(z; \rho_{5}|\clC_n(\bfa)|)\cap \omfF\subseteq \bigcup_{\bfb\in G_n(z,\bfa,\rho_5)} \oclC_n(\bfb).
\]
\end{lem01}
\subsubsection{Proof of Lemma \ref{PROPO:GC01:Rho:01}}
Let us assume that $f=\Mir_1$, so $\xi=\Mir_1(\zeta_1)= -\frac{1}{2} - i\alpha$, the same argument will hold for any $f\in \Di_8$ by Subsection \ref{Subsection:Symmetries oF Regular Cylinders}. 
Define the sequences $\bfd^1 \colon= (d_{1,n})_{n\geq 1}$, $\bfd^2 \colon= (d_{2,n})_{n\geq 1}$, and $\bfd^3 \colon= (d_{3,n})_{n\geq 1}$ by  
\begin{align*}
\bfd^1 &\colon= (-2,-2i,2, 2i, -2,-2i,2, 2i, \ldots ), \\
\bfd^2 &\colon= (-2 + i, 2i, 2, -2i, -2, 2i, 2, -2i, -2,\ldots ), \\
\bfd^3 &\colon= (-2 + i, 1+2i, -2, -2i, 2, 2i, -2, -2i, 2, 2i,  \ldots );
\end{align*}
that is, $\bfd^1$,$(d_{2,n})_{n\geq 2}$, and  $(d_{3,n})_{n\geq 3}$ are purely periodic sequences of period length $4$.

\begin{figure}[h!]
\begin{center}
\includegraphics[scale=0.85,  trim={5cm 19cm 11cm 6cm},clip]{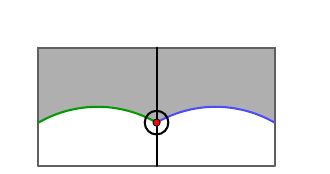}
\caption{  The set $\iota[\omfF_1(-2+i)]\cap(\tau_{2i}[\omfF]\cup\tau_{1+2i}[\omfF])$ and the disc $\Dx\left( \iota T_1(\xi;-2+i)  ;\tfrac{1}{10} \right) $. \label{Fig:GC02:02}}
\end{center}
\end{figure}

The proof technique of the next proposition will be used extensively in the forthcoming arguments.
\begin{propo01}\label{Prop:Rho:01:01}
Let $n\in\Na$. A word $\bfa\in \sfR(n)$ satisfies $\xi\in \oclC_n(\bfa)$ if and only if it is a prefix of $\bfd^1$, $\bfd^2$ or $\bfd^3$.
\end{propo01}
\begin{proof}
Consider $n=1$. By Figure \ref{Fig:02}, the only numbers $a_1\in \scD$ satisfying $\xi\in\oclC_n(a_1)$ are $a_1=-2$ and $a_1=-2+i$. Alternatively, from
\[
-\frac{5}{2}
<
\real(\xi^{-1})
<
-\frac{3}{2}
\;\text{ and }\;
\imag(\xi^{-1})=\frac{1}{2},
\]
we conclude $\xi^{-1}\in \tau_{a_1}[\omfF]$ if and only if $a_1\in \{-2,-2+i\}$.

Now consider $n=2$. Let $(a_1,a_2)\in\scD^2$ be such that 
\[
(a_1,a_2)\in \sfR(2)
\;\text{ and }\;
\xi\in\oclC_2(a_1,a_2).
\]
By the case $n=1$, we must have $a_1\in \{-2,-2+i\}$. First, assume that $a_1=-2$. We claim that $a_2=-2i$ is the single $a_2\in\scD$ satisfying 
\[
(-2,a_2)\in\sfR(2)
\;\text{ and }\;
T_1(\xi;-2)\in \omfF(-2)\cap \oclC_n(a_2).
\]
\begin{figure}[h!]
\begin{center}
\includegraphics[scale=0.5,  trim={5cm 10.5cm 6.0cm 6.5cm},clip]{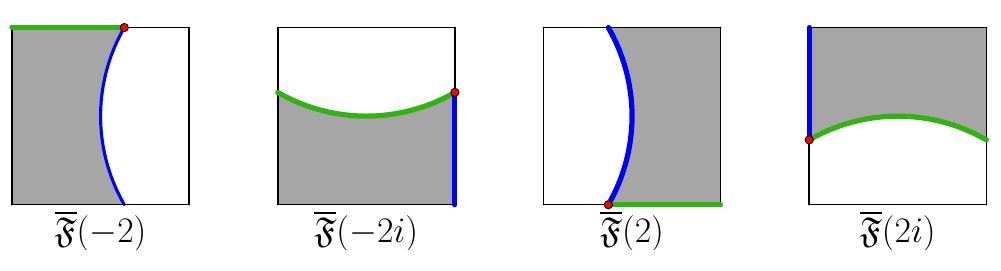}
\caption{Possible configurations \label{FigureB}}
\end{center}
\end{figure}
\normalsize
Indeed, consider the set $\omfF_1(-2) = \tau_2\iota[\oclC_1(-2)]\cap \omfF$, which is the first configuration depicted in Figure \ref{FigureB}. We have distinguished the point $T_2(\xi;-2)$ in red, the arch $T_2(\Fr(\mfF)\cap \oclC_1(2);-2)$ in blue and the line segment $T_2(\oclC_1(-2)\cap \oclC_1(-2+i);-2)$ in green. We apply $\iota$ on $\omfF_1(-2)$ and obtain the set $\iota[\omfF_1(-2)]$, which is depicted in Figure \ref{Fig:GC02:01-01}. 

\footnotesize
\begin{figure}[h]
\begin{center}
\includegraphics[scale=0.85,  trim={5cm 15cm 10cm 4cm},clip]{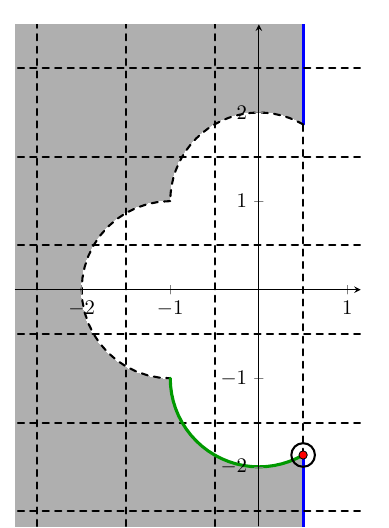}
\caption{ The set $\iota[\mfF_1(-2)]$ and the disc $\Dx\left(\iota T_1(\xi;-2) ;\tfrac{1}{10} \right)$. \label{Fig:GC02:01-01}}
\end{center}
\end{figure}
\normalsize

Observe that
\[
\iota \tau_2\iota (\xi)
\in 
\tau_{-2i}[\omfF]\cap \tau_{1-2i}[\omfF].
\]
However, $(-2,-2i)\in \sfR(2)$ and $(-2,1-2i)\not\in \sfR(2)$, so $(a_1,a_2)=(-2,-2i)$. This shows our claim. We also conclude that 
\[
\omfF_1(-2i) = T_2(\oclC_1(2);(-2,-2i))\cap \omfF,
\]
thus obtaining the second configuration in Figure \ref{FigureB}. 
In the second configuration in Figure \ref{FigureB}, the red point is $T_2(\xi;(-2,-2i))$; 
the green arch is the restriction to $\omfF$ of the image under $T_2(\,\cdot\,;(-2,-2i))$ of the green segment in the first configuration, and a similar interpretation for the blue segment.

We may repeat the argument to conclude that $a_3=2$ is the only solution $a_3\in \scD$ of
\[
(-2,-2i, a_3)\in\sfR(3)
\;\text{ and }\;
\xi\in \oclC_3(-2,-2i,a_3).
\]
Moreover, we get $\omfF_1(2) = T_3(\oclC_1(2);(-2,-2i, 2))\cap \omfF$, and obtain the third configuration in Figure \ref{FigureB}. Likewise, we can show that $a_4=2i$ is the only $a_4\in\scD$ for which
\[
(-2,-2i, 2,a_4)\in\sfR(4)
\;\text{ and }\;
\xi\in \oclC_3(-2,-2i,2, a_4).
\]
This argument will lead us to the fourth configuration in Figure \ref{FigureB}, which depicts the set $\omfF_1(2i)=T_4(\oclC_1(2);(-2,-2i,2,2i))$. Finally, the same strategy shows that $a_5=-2$ is the only solution of 
\[
(-2,-2i, 2, 2i, a_5)\in\sfR(5)
\;\text{ and }\;
\xi\in\oclC_5(-2,-2i, 2, 2i, a_5).
\]
and we obtain again the first configuration in Figure \ref{FigureB}. 
Therefore, we conclude that if $n\in\Na$ and $\bfa=(a_1,\ldots, a_n)\in\sfR(n)$ satisfy $a_1=-2$ and $\xi\in\oclC_n(\bfa)$, then $\bfa=\prefi(\bfd^1;n)$.

The proof is similar when $a_1=-2+i$. First, we show that $a_2\in \{2i, 1+2i\}$. Consider the set $\iota[\omfF_1(-2+i)]\cap(\tau_{2i}[\omfF]\cup\tau_{1+2i}[\omfF])$, depicted in Figure \ref{Fig:GC02:02}. We have distinguished the point $T(\xi;-2+i)$ in red, the arch $T(\oclC_1(-2)\cap \oclC_1(-2+i);-2+i)$ in green, and the arch $T(\Fr(\omfF)\cap \oclC_1(-2+i);-2+i)$ in blue. We conclude $a_2\in \{2i,1+2i\}$ arguing as above. Furthermore, we may proceed as in the case $a_1=-2$ to show that for any $n\in\Na_{\geq 3}$ and $\bfa\in\sfR(n)$ such that $\xi\in \oclC_n(\bfa)$ and $(a_1,a_2)=(-2 + i, 2i)$ (resp. $(a_1,a_2)=(-2 + i, 1+2i)$), we have $\bfa=\prefi(\bfd^2;n)$ (resp. $\bfa=\prefi(\bfd^3;n)$).
\end{proof}

\begin{coro01}\label{Corollary:GC1}
For some absolute constants, we have that if $n\in\Na$ and $\bfa,\bfb\in\sfR(n)$ are such that $\xi \in \oclC_n(\bfa)\cap \oclC_n(\bfb)$, then $|\clC_n(\bfa)| \asymp |\clC_n(\bfb)|$.
\end{coro01}
\begin{proof}
Let $n\in\Na$. The corollary is obvious if $n\in \{1,2\}$, for there are finitely many cylinders satisfying the hypothesis. Suppose that $n\geq 3$. Let $\bfa=(a_1, \ldots, a_n)$, $\bfb=(b_1,\ldots, b_n)\in\sfR(n)$ be such that $\xi\in \oclC_n(\bfa)\cap \oclC_n(\bfb)$. By Proposition \ref{Prop:Rho:01:01}, there is some $f$ in the group generated by $\Mir_1$ and $\Mir_2$ such that $(a_3,\ldots, a_n) 
= \left(f(b_3), \ldots, f(b_n)\right)$. Then, Proposition \ref{PROPO:SIMET:CILINDROS:02} and Corollary \ref{CORO:SIMET:CILINDROS:02} give
\[
\left|\clC_{n-2}(a_3,\ldots, a_n)\right|
=
|\clC_{n-2}(b_3,\ldots, b_n)|.
\]
Moreover, we also have
\[
\clC_n(\bfa)
=
Q_2\left( \clC_{n-2}(a_3,\ldots, a_n);(a_1,a_2)\right).
\]
Since there are only three options for $(a_1,a_2)$, for some absolute constants we have
\[
\left|\clC_n(\bfa) \right|
\asymp 
\left|\clC_{n-2}(a_3,\ldots, a_n)\right|.
\] 
Therefore, we conclude
\[
\left|\clC_n(\bfa) \right|
\asymp
\left|\clC_{n-2}(a_3,\ldots, a_n) \right|
=
\left|\clC_{n-2}(b_3,\ldots, b_n) \right|
\asymp
\left|\clC_n(\bfb) \right|.
\]
\end{proof}

\begin{proof}[Proof of Lemma \ref{PROPO:GC01:Rho:01}]
The lemma will be proven once we show the existence of some $\rho_1>0$ such that for any $j\in\{1,2,3\}$ and $n\in \Na$, we have
\begin{equation}\label{Eq:proof:rho1:00}
\Dx\left( \xi; \rho_1\left| \clC_n(\bfd^j)\right|\right)
\cap
\omfF
\subseteq \bigcup_{k=1}^3 \oclC_n\left(d_{k,1},\ldots, d_{k,n}\right).
\end{equation}
First, note that every small $r'>0$ verifies
\begin{equation}\label{Eq:proof:rho1:01:bis}
\omfF\cap \Dx(\xi;r')
\subseteq 
\oclC_1(-2)\cap \oclC_1(-2+i).
\end{equation}
Indeed, by Figure \ref{Fig:GC02:01}, for $0<r<\frac{1}{10}$ we have 
\[
\Dx\left( \iota(\xi);r\right)
\cap \iota\left[ \mfF \right]
\subseteq 
\left(\tau_{-2}\left[ \mfF \right]\cup \tau_{-2+i}\left[ \mfF \right]\right)\cap \iota\left[ \mfF \right].
\]
\begin{figure}[h!]
\begin{center}
\includegraphics[scale=0.85,  trim={5cm 19cm 13.5cm 4cm},clip]{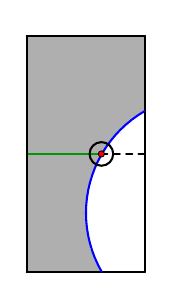}
\caption{ The set $\iota[\mfF]\cap(\tau_{-2}[\omfF]\cup\tau_{-2+i}[\omfF])$ and the disc $\Dx\left(\iota(\xi), 1/10 \right)$. \label{Fig:GC02:01}}
\end{center}
\end{figure}

Then, we apply $\iota$ to obtain \eqref{Eq:proof:rho1:01:bis} for every small $r'>0$. This proves \eqref{Eq:proof:rho1:00} for $n=1$ and $j\in \{1,2,3\}$. Next, we show the existence of some constant $\tilde{\rho}_1>0$ such that for $n=2$, we have
\begin{align}
\Dx\left( \xi; \tilde{\rho}_1\left| \oclC_2(-2,-2i)\right| \right)
\cap
\oclC_1(-2)
&\subseteq 
\oclC_2(-2, -2i), \label{Eq:proof:rho1:01} \\
\Dx\left( \xi; \tilde{\rho}_1\left| \oclC_2(d_{j,1} ,d_{j,2})\right| \right)
\cap
\oclC_1(d_{j,1} )
&\subseteq 
\oclC_2(-2+i, 2i) \cup \oclC_2(-2+i, 1+2i) \; \text{ for } j=1,2. \label{Eq:proof:rho1:02}
\end{align}
As in the proof of Proposition \ref{Prop:Rho:01:01}, we conclude from Figure \ref{Fig:GC02:01-01} that 
\[
\Dx\left( \iota\left(T_1(\xi;-2);\frac{1}{10}\right)\right) \cap \iota \left[\omfF_1(-2)\right]
\subseteq 
\tau_{-2i}[\omfF]\cap \iota \left[\omfF_1(-2)\right].
\]
Because of $\iota[\omfF_1(-2)]= \iota\tau_2\,\iota \left[\oclC_1(-2)\right]$, we apply $\tau_{2i}$ on both sides to get
\[
\Dx\left( T_2(\xi;-2,-2i);\frac{1}{10}\right)\cap T_2(\clC_1(-2);(-2,-2i))
\subseteq 
\omfF \cap  T_2(\oclC_1(-2);(-2,-2i)).
\]
We now apply $Q_2 ( \,\cdot\, , (-2,-2i) )$ and obtain
\begin{align*}
Q_2\left( \Dx\left( T_2(\xi;(-2,-2i));\frac{1}{10}\right), (-2,-2i) \right) \cap \oclC_1(-2)
&\subseteq 
Q_2(\omfF;(-2,-2i)) \cap  \oclC_1(-2) \\
&= \oclC_2(-2,-2i).
\end{align*}
By Proposition \ref{PROPO:GC01:QnDISCOS}, there is some $\tilde{\rho}_1>0$ for which \eqref{Eq:proof:rho1:01} is true. Similarly, we can show \eqref{Eq:proof:rho1:02} using Figure \ref{Fig:GC02:02}. Furthermore, essentially the same argument and the proof of Proposition \ref{Prop:Rho:01:01} allow us to conclude that, for $n\geq 3$ and $j\in\{1,2,3\}$, the same $\tilde{\rho}_1>0$ satisfies
\begin{equation}\label{Eq:proof:rho1:03}
\Dx\left( \xi; \tilde{\rho}_1\left| \clC_n\left(d_{j,1}, \ldots, d_{j,n}\right)\right|\right)
\cap 
\oclC_{n-1}\left(d_{j,1}, \ldots, d_{j,n-1}\right)
\subseteq
\oclC_n\left(d_{j,1}, \ldots, d_{j,n}\right).
\end{equation}
Lastly, in view of Corollary \ref{Corollary:GC1}, we conclude the existence of some $\rho_1>0$ such that \eqref{Eq:proof:rho1:00} holds for all $n$ and $j$. 
\end{proof}

\subsubsection{Proof of Lemma \ref{PROPO:GC01:Rho:02}}

Proposition \ref{Prop:Rho:01:01} lies at the heart of the proof of Lemma \ref{PROPO:GC01:Rho:01}. It characterizes the three closed cylinders containing a given $\xi\in \Di_8(\zeta_1)$. We will revisit this argument in the proof of Lemma \ref{PROPO:GC01:Rho:02} and afterward. We will also use the fact that the only closed cylinders intersecting $\Esq$ are $\oclC_0(\veps)\colon= \omfF$, $\oclC_1(1+i)$, $\oclC_1(1-i)$, $\oclC_1(-1+i)$, and $\oclC_1(-1-i)$ (see Proposition \ref{GC:PROPO:02}).

\begin{propo01}\label{Lemma:GC01:rho:02}
Let $n\in\Na$, $\bfa\in\sfR(n)$, and $b,c\in\scD$ such that $b\neq c$, and $\bfa b,\bfa c\in \sfR(n+1)$. Assume that there is some $w\in \Cx$ such that
\[
w\in \Fr\left( \mfF_n(\bfa) \right) \cap \oclC_1(b)\cap \oclC_1(c).
\]
Then, the next assertions hold:
\begin{enumerate}[\rm i.]
\item \label{Lemma:GC01:rho:02:01} $\tau_{-b}\iota (w)\in \Esq\cup \Di_8(\zeta_1)$.
\item \label{Lemma:GC01:rho:02:02} The numbers $b$ and $c$ are the only solutions $d\in\scD$ to
\[
\tau_{-d}\iota(w)\in \omfF 
\;\text{ and }\;
\bfa d \in \sfR(n+1).
\]
\item \label{Lemma:GC01:rho:02:03} If $\tau_{-b}\iota(w)\in \Esq$, then the disc $\iota[\Dx(\iota(w);\frac{1}{10})]$ intersects exactly one side of $\Fr(\mfF_n(\bfa))$.
\end{enumerate}
\end{propo01}
\begin{proof}
The three parts are shown by checking the possible forms of $\omfF_n(\bfa)$. For clarity, we briefly discuss one option in the first point.
Assume that $\omfF_n(\bfa)=\omfF(-2)$. By Figure \ref{Fig:iotaF1-2}, when $|b|\leq \sqrt{5}$ and $|c|\leq \sqrt{5}$, any number belonging to $\iota[\Fr\left( \mfF_n(\bfa) \right) \cap \oclC_1(b)\cap \oclC_1(c)]$ is a Gaussian-integral translate of some element in $\Di_8(\zeta_1)$. And if $|b|\geq  3$ or $|c|\geq 3$, then any number belonging to $\iota[\Fr\left( \mfF_n(\bfa) \right) \cap \oclC_1(b)\cap \oclC_1(c)]$ is a Gaussian-integral translate of some point in $\Esq$.
\end{proof}
\begin{proof}[Proof of Lemma \ref{PROPO:GC01:Rho:02}]
Let $k\in\{1,\ldots,n-1\}$ be the minimal integer for which there is some $\bfc=(c_1,\ldots, c_n)\in \sfR(n)$ such that
\[
z\in \oclC_n(\bfc), \quad
(a_1,\ldots, a_k)=(c_1,\ldots, c_k),
\quad\text{ and }\quad
a_{k+1} \neq c_{k+1}.
\]
Then, we have
\[
w
\colon=
T_k(z;(a_1,\ldots, a_k))
\in 
\oclC_{n-k}(a_{k+1},\ldots, a_n)\cap \oclC_{n-k}(c_{k+1},\ldots, c_n)\cap \omfF_k(a_1,\ldots, a_k).
\]
From $a_{k+1}\neq c_{k+1}$ we get $z\in \Fr(\clC_k(a_1,\ldots, a_k))$ and, thus, $w\in \Fr(\mfF_k(a_1,\ldots, a_k))$. Therefore, 
\[
w
\in 
\oclC_1(a_{k+1})\cap \oclC_1(c_{k+1})\cap \Fr(\mfF_k(a_1,\ldots, a_k))
\]
and, by Proposition \ref{Lemma:GC01:rho:02}.\ref{Lemma:GC01:rho:02:01}, we must have either
\[
\tau_{-a_{k+1} }\iota(w)\in \Esq
\;\text{ or }\;
\tau_{-a_{k+1}}\iota(w)\in \Di_8(\zeta_1).
\]
We discuss both possibilities separately.

\textbf{Case $\tau_{-a_{k+1}}\iota(w)\in \Esq$.} Proposition \ref{Lemma:GC01:rho:02}.\ref{Lemma:GC01:rho:02:02} tells us there is a single possibility for $c_{k+1}$. Moreover, in view of
\[
\iota(w)
\in
\tau_{a_{k+1}}
\left[\oclC_{n-k-1}(a_{k+2},\ldots, a_n) \cap \Esq\right],
\]
we must have $k\in \{n-1,n-2\}$ by Proposition \ref{GC:PROPO:02}.

First, assume that $k=n-1$. Then, $\bfc=(a_1,\ldots, a_{n-1}, c_n)$ and, by Proposition \ref{Lemma:GC01:rho:02}.\ref{Lemma:GC01:rho:02:03}, the disc $\iota[\Dx(\iota(w);\frac{1}{10})]$ intersects only one side of $\Fr(\mfF_{n-1}(a_1,\ldots,a_{n-1}))$. From 
\[
\Dx
\left( \iota(w);\frac{1}{10} \right)
\cap
\iota\left[ \omfF_k(a_1,\ldots, a_{n-1}) \right]
\subseteq
\left(  \tau_{a_{n}}[\omfF]\cup \tau_{c_n} [\omfF] \right)
\cap 
\iota \left[ \omfF_{n-1}(a_1,\ldots, a_{n-1})\right]
\]
we get
\[
\iota\left[
\Dx\left(\iota(w);\frac{1}{10}\right)
\cap
\left[ \omfF_{n-1}(a_1,\ldots, a_{n-1})\right]
\right]
\subseteq
\left( 
\oclC_1 (a_{n})  \cup \oclC_1(c_n)
\right)
\cap 
\omfF_{n-1}(a_1,\ldots, a_{n-1}).
\]
We now apply $Q_k(\,\cdot\,;(a_1,\ldots, a_k))=Q_{n-1}(\,\cdot\,;(a_1,\ldots, a_{n-1}))$ and use Proposition \ref{Lemma:GC01:rho:02} to conclude the existence of some absolute constant $\rho_{2,2}>0$ such that 
\[
\Dx\left( z; \frac{\rho_{2,2} |\clC_{n-1}(a_1,\ldots, a_{n-1})|}{|a_n|^2}\right)
\cap 
\oclC_{n-1}(a_1,\ldots, a_{n-1})
\subseteq 
\oclC_n(\bfa)\cup \oclC_n(\bfc).
\]
Hence, by Proposition \ref{Prop:SizeEstimates}.\ref{Prop:SizeEstimates:i}, there is some absolute constant $\rho_{2,3}>0$ such that
\[
\Dx(z;\rho_{2,3}|\oclC_n(\bfa)|)
\cap 
\oclC_{n-1}(a_1,\ldots, a_{n-1})
\subseteq 
\oclC_n(\bfa)\cup \oclC_n(\bfc).
\]
To finalize the proof of this case, we must verify that
\begin{equation}\label{Eq:Prop:6.6:CorrectInclusion}
\Dx(z;\rho_{2,3}|\oclC_n(\bfa)|)
\cap
\omfF
=
\Dx(z;\rho_{2,3}|\oclC_n(\bfa)|)
\cap 
\oclC_{n-1}(a_1,\ldots, a_{n-1}).
\end{equation}

To this end, note that the disc $\Dx(z;\rho_{2,3}|\clC_n(\bfa)|)$ intersects exactly one side of $\clC_n(\bfa)$, which is the one containing $z$. This side is a line segment contained in $\Fr(\mfF)$. Hence, $\Dx(z;\rho_{2,3}|\clC_n(\bfa)|)$ only intersects one side of $\clC_{n-1}(a_1,\ldots, a_{n-1})$. If there was some $\bfb \in \sfR(n-1)$, $\bfb\neq (a_1,\ldots, a_{n-1})$, such that $\Dx(z;\rho_{2,3}|\oclC_n(\bfa)|)$ had an inner point of $\oclC_{n-1}(\bfb)$, then this disc would have to intersect at least two sides of $\clC_{n-1}(a_1,\ldots, a_{n-1})$. This proves \eqref{Eq:Prop:6.6:CorrectInclusion}.

When $n=k-2$, we proceed in a similar fashion, so we omit the proof. However, we note that in this case we must have $|a_n|=|c_n|=|1+i|=\sqrt{2}$ (see Figure \ref{Fig:mfFPart}).

\textbf{Case $\tau_{-a_{k+1}}\iota(w)\in \Di_8(\zeta)$.} After checking all the possible options for $\mfF_n(\bfa)$, we may say that there are $f_1,f_2\in \Di_8$ such that 
\[
\iota(w)
\in 
\{\tau_{a_{k+1}}f_1(\zeta_1)\}
\cap
\{\tau_{c_{k+1}}f_2(\zeta_1)\}
\cap
\Fr(\iota[\mfF_k(a_1,\ldots, a_k)]).
\]
Figure \ref{Fig:GC02:05} depicts all the possible configurations$\mod \RotaG$ for $\iota(w)$ (in red) and $\iota[\mfF_k(a_1,\ldots, a_k)]$ (gray). For example, case (c) arises from
\[
z\in 
\oclC_3(-2,4i, 1+i)\cap \oclC_3(-2,4i, 2+i)\cap \Fr(\mfF).
\]
\footnotesize
\begin{figure}[h!]
\begin{center}
\includegraphics[scale=0.85,  trim={5cm 21cm 4cm 4.5cm},clip]{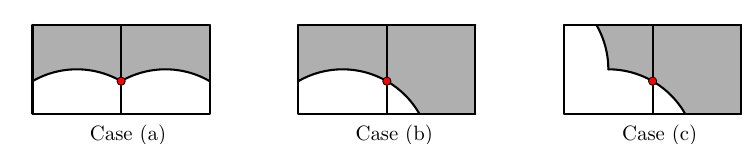}
\caption{Possible location of $\iota(w)$ in $\iota[\omfF_{k}(a_1,\ldots, a_{k})]$. \label{Fig:GC02:05}}
\end{center}
\end{figure}
\normalsize
We are now essentially in the situation of Lemma \ref{PROPO:GC01:Rho:01}. Arguing accordingly, we conclude the existence of some constant $\rho_{2,3}>0$ for which
\begin{align*}
\Dx\left( \iota(w); \rho_{2,3}|\oclC_{n-k-2}(a_{k+2},\ldots, a_n)|\right) 
&\cap
\iota[ \mfF_k(a_1,\ldots, a_k)]
\subseteq \\
&\subseteq 
\oclC_{n-k-2}(a_{k+2},\ldots, a_n)
\cup
\oclC_{n-k-2}(c_{k+2},\ldots, c_n).
\end{align*}
As in the proof of \eqref{Eq:Prop:6.6:CorrectInclusion}, we may conclude the existence of some constant $\rho_{2,4}>0$ such that 
\[
\Dx\left( z; \rho_{2,4}|\oclC_{n}(\bfa)|\right) 
\cap
\omfF
\subseteq
\oclC_n(\bfa)
\cup
\oclC_n(a_1,\ldots, a_k, b_{k+1}, \ldots, b_n).
\]
Lastly, we choose $\rho_2\colon=\min\{\rho_{2,2},\rho_{2,4}\}$ and the lemma follows. Observe that in case (a) of Figure \ref{Fig:GC02:05}, we have $\# G_n(z,\bfa,\rho_2)=2$ and $\# G_n(z,\bfa,\rho_2)=3$ in cases (b) and (c). 
\end{proof}

\subsubsection{Proof of Lemma \ref{PROPO:GC01:Rho:03}}

Let $n$, $\bfa$, and $z$ be as in the statement. For $\rho_2>0$ as in Lemma \ref{PROPO:GC01:Rho:02} and define
\[
\delta
\colon=
\min
\left\{ 
\frac{1}{10}, 
 \kappa_1\left( 1- \frac{\sqrt{2}}{2}\right)^2  \frac{ \rho_2}{2}\right\}.
\]
Since $z\in \Fr(\mfF)\cap \oclC_n(\bfa)$, we have $z\in \Fr(\clC_{n-1}(a_1,\ldots, a_{n-1} ))$, so
\[
\iota T_{n-1}(z;(a_1,\ldots, a_{n-1}))
\in 
\Fr\left(\iota[\mfF_{n-1}(a_1,\ldots, a_{n-1}]\right)
\]
and there are at most two solutions $d\in\scD$ of
\begin{equation}\label{EQ:GC02:03}
\tau_{d}[\mfF]\cap \iota[\mfF_{n-1}(a_1,\ldots, a_{n-1}) ]
\cap 
\Dx(\iota(T_{n-1}(z;(a_1,\ldots, a_{n-1})); \delta))
\neq \vac
\end{equation}
one of them being $d=a_{n}$ (cf. Lemma \ref{Lemma:GC01:rho:02}.\ref{Lemma:GC01:rho:02:02}). 

When $d=a_n$ is the only solution of \eqref{EQ:GC02:03}, we may show as in Lemma \ref{PROPO:GC01:Rho:01} that
\[
\Dx\left(z; \frac{\delta}{2\left(1+\frac{\sqrt{2}}{2}\right) \left( 2 + \frac{\sqrt{2}}{2}\right)} |\clC_n(\bfa)|\right) 
\cap \omfF
\subseteq 
\oclC_n(\bfa).
\]
Assume that \eqref{EQ:GC02:03} has two different solutions, say $d=a_{n}$ and $d=b$. Then, we can single out some $\xi\in\Di_8(\zeta_1)\cup\Esq$ for which
\[
\Dx\left(\iota\left(T_{n-1}(z;(a_1,\ldots, a_{n-1}) \right); \delta) )\right)
\cap
\Fr(\mfF_{n-1}(a_1,\ldots, a_{n-1}))
\cap
\tau_{a_n}\left[\,\omfF\,\right]
\cap 
\tau_{b}\left[\,\omfF\,\right]
=
\{\tau_{a_n}(\xi)\}.
\]
Hence, $|\xi - T_n(z;\bfa)|< \delta$ and
\begin{align*}
|Q_n(\xi;\bfa) - z| &= 
\left| Q_n(\xi;\bfa) - Q_n(T_n(z;\bfa);\bfa) \right| \\
&= 
\frac{|\xi - T_n(z;\bfa)|}{|q_n(\bfa)|^2\left| 1- \frac{q_{n-1}}{q_{n}} Q_n(\xi;\bfa)\right|\left| 1- \frac{q_{n-1}}{q_{n}} Q_n(T_n(z;\bfa);\bfa)\right| } \\
&< 
\frac{\delta}{\kappa_1 \left( 1- \frac{\sqrt{2}}{2}\right)^2} |\clC_n(\bfa)| && \text{(by }Q_n(\xi;\bfa)\in\omfF\text{)} \\
&\leq 
\frac{ \rho_2}{2}|\clC_n(\bfa)|.
\end{align*}
Then, for all $0<r< \frac{\rho_2}{2}$, we have $\Dx(z;r |\clC_n(\bfa)|) \subseteq  \Dx(Q_n(\xi;\bfa); \rho_2|\clC_n(\bfa)|)$. The result follows from Lemma \ref{PROPO:GC01:Rho:02} after picking $\rho_3=\frac{ \rho_2}{3}$.

\begin{rema01}\label{Rem:rho3}
Our choice of $\rho_3$ implies that if $n\in\Na$, $\bfa\in\sfR(n)$, and $\xi\in \Di_8(\zeta_1)\cup\Esq$ is such that $Q_n(\xi;\bfa)\in \Fr(\clC_n(\bfa))$, then
\[
Q_n\left(  \Di_8(\zeta_1) \cup\Esq;\bfa) \cap \Dx(Q_n(\xi;\bfa); \rho_3\right)
=
\{Q_n(\xi;\bfa)\}.
\]
\end{rema01}
 \subsubsection{Proof of Lemma \ref{PROPO:GC01:Rho:05}}
\begin{figure}[h!]
\begin{center}
\includegraphics[scale=0.85,  trim={4.5cm 19cm 10cm 4cm},clip]{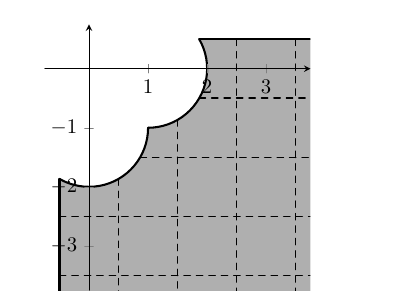}
\caption{ The set $\iota[\omfF_1(1+i)]$. \label{Fig:GC02:06}}
\end{center}
\end{figure}

Let $n$, $\bfa$, and $z$ be as in the statement. Let $k\in \{0,1\ldots, n-1\}$ be such that 
\[
z\in \clCc_k(a_1,\ldots, a_k)\cap \Fr(\clC_{k+1}(a_1,\ldots, a_k, a_{k+1})).
\]
As before, we interpret $\clC_0(\epsilon)=\mfF$ and $T_0(\,\cdot\,:\epsilon)$ as the identity map. Let us first show that
\begin{equation}\label{Eq:Rho5:01}
w\colon=
\iota T_k(z;(a_1,\ldots, a_k))
\in
\Fr(\tau_{a_{k+1}}[\mfF]).
\end{equation}
For clarity, we discuss the conclusion on an example. When $\omfF_k(a_1,\ldots,a_k)=\omfF_1(1+i)$, the set $\iota[\mfFc_k(a_1,\ldots,a_k)]$ is the one depicted in Figure \ref{Fig:GC02:06}. Then, the inclusion \eqref{Eq:Rho5:01} means that $w$ belongs to one of the dashed lines.

To show \eqref{Eq:Rho5:01}, note that 
\[
w
\in
\iota\left[\mfFc_k(a_1,\ldots, a_k) \right]
\cap
\iota 
T_k
\left(
\Fr\left(
\clC_{k+1}(a_1,\ldots, a_k, a_{k+1})\right) ;
(a_1,\ldots, a_k)
\right).
\]
Then, in view of  
\begin{align*}
\iota 
T_k
\left(
\clC_{k+1}(a_1,\ldots, a_k, a_{k+1}));
(a_1,\ldots, a_k)\right)
&=
\iota\left[\mfF_k(a_1,\ldots, a_k)\cap \clC_1(a_{k+1})\right] \nonumber\\
&=
\iota\left[\mfF_k(a_1,\ldots, a_k) \right]
\cap 
\tau_{a_{k+1}}[\mfF] 
\cap 
\iota [\mfF] \nonumber\\
&=
\iota\left[\mfF_k(a_1,\ldots, a_k) \right]
\cap 
\tau_{a_{k+1}}[\mfF],
\end{align*}
we have
\[
w
\in
\iota\left[\mfFc_k(a_1,\ldots, a_k) \right]
\cap
\Fr\left(
\iota\left[\mfF_k(a_1,\ldots, a_k) \right]
\cap 
\tau_{a_{k+1}}[\mfF]
 \right).
 \]
From here, \eqref{Eq:Rho5:01} is easily deduced.
Let us further assume that $k\leq n-2$. Take $0<\tilde{\rho}<\rho_3$, we impose later further restrictions on $\tilde{\rho}$. For $\kappa_1$ and $\kappa_2$ as in Proposition \ref{Prop:SizeEstimates}, put
\[
\delta
\colon= 
\frac{\kappa_1}{\kappa_2}  \left( 1- \frac{\sqrt{2}}{2}\right)^2.
\]
We now consider two cases depending on how far is $w$ from the boundary of $\iota\left[ \mfFc_{k}(a_1,\ldots, a_k)\right]$. More precisely, we first assume that
\begin{equation}\label{EQ:GC02:06}
\Dx
\left(
w; \frac{\delta \tilde{\rho} }{2}|\clC_{n-k-2}(a_{k+2},  \ldots, a_n)|
 \right)
\subseteq
\iota\left[ \mfFc_{k}(a_1,\ldots, a_k)\right].
\end{equation}
By $k\leq n-2$ and Proposition \ref{GC:PROPO:02} we have $w\not\in \tau_{a_{k+1}}[\Esq]$, so Remark \ref{Rem:rho3} tells us that the set
\[
\tau_{a_{k+1}}[\omfF]
\cap
\Dx
\left(
w; \frac{\delta \tilde{\rho} }{2}|\clC_{n-k-2}(a_{k+2},  \ldots, a_n)|
 \right)
\]
is contained in at most three sets of the form $\tau_{a_{k+1}}[\clC_{n-k-2}(\bfb)]$, $\bfb\in \sfR(n-k-2)$. Moreover, there is a single $b\in \scD$ with $|b-a_{k+1}|=1$ and 
\[
\Dx
\left(
w; \frac{\delta \rho_{3}}{2}|\clC_{n-k-2}(a_{k+2},  \ldots, a_n)|
 \right)
\subseteq
\tau_{a_{k+1}}[\omfF]
\cup
\tau_{b}[\omfF].
\]
By Lemma \ref{PROPO:GC01:Rho:03} and Subsection \ref{Subsection:Symmetries oF Regular Cylinders}, the previous disc is contained in the closure of at most six integral translates of closed regular cylinders of level $n-k-2$. 

Now assume that \eqref{EQ:GC02:06} fails. Then, there is some $\xi\in \Di_8(\zeta_1)\cup\Esq$ such that 
\[
\tau_{a_k}(\xi) 
\in 
\Fr\left(\iota(\mfF_k(a_1,\ldots, a_k))\right) 
\cap
\Dx
\left(
w; \frac{\delta \tilde{\rho} }{2}|\clC_{n-k-2}(a_{k+2},  \ldots, a_n)|
\right).
\]
Using $w=\iota T_k(z;(a_1,\ldots, a_{k}))$ and $\tau_{-a_{k+1}}(w) = T_{k+1}(z;(a_1,\ldots, a_{k+1}))$, we get
\[
\left|T_{k+1}\left(z;(a_1,\ldots,a_{k+1} )\right) -  \xi\right|
< 
\frac{\delta \tilde{\rho}}{2}|\clC_{n-k-2}(a_{k+2},  \ldots, a_n)|.
\]
Therefore, writing $\xi'=Q_{k+1}(\xi;(a_1,\ldots, a_{k+1}))$, we obtain
\begin{align*}
|z - \xi'| 
&= 
|Q_{k+1}\left(T_{k+1}(z;(a_1,\ldots, a_{k+1});(a_1,\ldots, a_{k+1}) \right) - Q_{k+1}(\xi;(a_1,\ldots, a_{k+1}))|  \\
&<
\frac{\left|T_{k+1}\left(z;(a_1,\ldots,a_{k+1} )\right) -  \xi\right|}{|q_{k+1}(a_1,\ldots, a_{k+1})|^2 \left( 1- \frac{\sqrt{2}}{2}\right)^2} \\ 
&<
\frac{\left| \clC_{k+1}(a_1,\ldots,a_{k+1})\right|}{\kappa_1 \left( 1- \frac{\sqrt{2}}{2}\right)^2} \, \left|T_{k+1}\left(z;(a_1,\ldots,a_{k+1} )\right) -  \xi\right| \\ 
&\leq
\frac{\delta \tilde{\rho} }{2 \kappa_1 \left( 1- \frac{\sqrt{2}}{2}\right)^2} 
\left| \clC_{k+1}(a_1,\ldots, a_{k+1})\right||\clC_{n-k-2}(a_{k+2},  \ldots, a_n)| \\
&=
\frac{ \tilde{\rho} }{2 \kappa_2} 
\left| \clC_{k+1}(a_1,\ldots, a_{k+1})\right||\clC_{n-k-2}(a_{k+2},  \ldots, a_n)| \\
&< \frac{ \tilde{\rho} }{2} |\clC_n(\bfa)|.
\end{align*}
As a consequence, if $\rho_4:=\frac{\tilde{\rho}}{2}$, we have 
\[
\Dx\left(z; \rho_4  |\clC_n(\bfa)|\right)
\subseteq 
\Dx\left( Q_{k+1}(\xi;(a_1,\ldots, a_{k+1})); \tilde{\rho} |\clC_n(\bfa)|\right).
\]
Next, we explain why we can choose $\tilde{\rho}$ so small that the larger disc in the previous inclusion is contained in the interior of at most six regular cylinders of level $k+1$. For concreteness, assume $\mfFc_k(a_1,\ldots,a_k)=\mfFc_1(1+i)$; the other options are treated similarly. First, if $\tau_{a_k}(\xi)$ is in one of the arches in $\Fr(\iota[\mfF_k(a_1,\ldots,a_k)])$ (see Figure \ref{Fig:GC02:06}), then $\iota\tau_{a_k}(\xi)$ belongs to one of the line segments in $\Fr(\mfF_1(1+i))$ (see Figure \ref{Fig:ClosedPrototypeSets}), we apply Lemma \ref{PROPO:GC01:Rho:03} and the proof is complete. If $\tau_{a_k}(\xi)$ belongs to one of the rays in $\Fr(\iota[\mfF_k(a_1,\ldots,a_k)])$, then $\iota\tau_{a_k}(\xi)$ belongs to one of the arches in $\Fr(\mfF_1(1+i))$ and, thus,
\[
\iota\tau_{a_{k-1}}\iota\tau_{a_{k}}(\xi)
=
Q_2(\xi;(a_{k-1},a_k))
\in\Fr(\mfF).
\]
To see it, apply $\tau_{-(1+i)}\iota$ on $\oclC_1(1+i)$ and refer to Figure \ref{Fig:mfFPart}. From this point, we can conclude after Lemma \ref{PROPO:GC01:Rho:03}.

The proof for $k=n-1$ is similar. In fact, we only used $k\leq n-2$ to avoid the possibility $\xi\in \Esq$. However, this situation is easier to deal with. We leave the details to the reader.

\subsubsection{Proof of Lemma \ref{PROPO:GC01:Rho:06}} 
For $\rho_4>0$ as in Lemma \ref{PROPO:GC01:Rho:05}, define
\[
\rho_5 \colon= 
\frac{\kappa_1\left( 1 - \frac{\sqrt{2}}{2}\right)^2}{2\left( 1 + \frac{\sqrt{2}}{2}\right)\left( 2 + \frac{\sqrt{2}}{2}\right)}\, \frac{ \rho_4}{2}.
\]
Let $n$, $\bfa$, and $z$ be as in the statement. We may show as before that if 
\begin{equation}\label{EQ:GC02:07}
\Dx\left(T_n(z;\bfa);\frac{\rho_5}{2}\kappa_1\left( 1 - \frac{\sqrt{2}}{2}\right)^2\right)
\subseteq
\mfFc_n(a_1,\ldots,a_n) 
=
T_n\left( \clCc_n(\bfa);\bfa \right),
\end{equation}
then
\[
\Dx\left(z; \rho_5 |\clC_n(\bfa)|\right)
\subseteq
\clCc_n(\bfa)
\subseteq \oclC_n(\bfa).
\]
When \eqref{EQ:GC02:07} is false, there is some $\xi\in \omfF$ such that
\[
\xi \in\Fr(\mfF_n(\bfa))
\quad\text{ and }\quad
|\xi - T_n(z; \bfa)|< \frac{ \rho_5}{2}\kappa_1\left( 1 - \frac{\sqrt{2}}{2}\right)^2.
\]
Then, as in the proof of Lemma \ref{PROPO:GC01:Rho:05}, we get
\[
|z - Q_n(\xi;\bfa)|
 < 
\frac{ \rho_5}{2} |\clC_n(\bfa)|,
\]
so 
\[
\Dx\left(z; \frac{ \rho_5}{2} |\clC_n(\bfa)|\right)
\subseteq
\Dx\left(Q_n(\xi;\bfa); \rho_5 |\clC_n(\bfa)|\right).
\]
The result now follows from Lemma \ref{PROPO:GC01:Rho:05}.

We highlight two important consequences of the proofs of lemmas \ref{PROPO:GC01:Rho:01} and \ref{PROPO:GC01:Rho:02}. First, two neighboring cylinders have a similar diameter (Corollary \ref{CORO:NeighCylProp}). Second, the partial quotients determining two neighboring cylinders are of similar size (Corollary \ref{CORO:NeighDigh}). 
The details are left to the reader.

\begin{coro01}\label{CORO:NeighCylProp}
There is a constant $\kappa_3>1$ such that for any $n\in\Na$, all $\bfa\in\sfR(n)$, and all $\bfb\in\Veci_n(\bfa)$
we have
\[
\frac{1}{\kappa_3} 
\leq 
\frac{|\clC_n(\bfa)|}{|\clC_n(\bfb)|} 
\leq 
\kappa_3.
\]
\end{coro01}
\begin{coro01}\label{CORO:NeighDigh}
For all $n\in\Na$ and $\bfa=(a_1,\ldots,a_n)\in\sfR(n)$, if $\bfb=(b_1,\ldots,b_n)\in \Veci_n(\bfa)$ then
\[
\|a_j\|- 1 \leq \|b_j\|\leq \|a_j\| +1
\;\text{ for all }j\in\{1,\ldots, n\}.
\]
\end{coro01}

\subsection{Second geometric construction}\label{Subsection:SecondGeomCons}
We dedicate this section to the proof of Theorem \ref{TEO:GC:02}. It counts how many sets of a certain form are intersected by discs with small radius. More precisely, for any $n\in\Na$, $\bfa\in\sfR(n)$, and $0<\mathcal{l} < \mathcal{u}$, define
\[
S_n(\bfa;\mathcal{l},\mathcal{u}) 
\colon=
\bigcup_{\substack{\mathcal{l}\leq \|d\|\leq \mathcal{u} \\ \bfa d \in \sfR(n+1)}}  \oclC_{n+1}(\bfa d).
\]
\begin{teo01}\label{TEO:GC:02}
There is a constant $\rho>0$ with the following property: for any two natural numbers $\cll$ and $\clu$ such that $33\leq  \cll <\clu$ and $\frac{27}{100}< \frac{\cll}{\clu}< \frac{3}{4}$, any $\bfa\in\sfR(n)$, and any $z\in S_n(\bfa;\mathcal{l}, \mathcal{u})$, we have
\[
\#
\left\{
\bfc\in \sfR(n):
\Dx(z;\rho|S_n(\bfa;\cll,\clu)|)
\cap 
S_n(\bfc;\mathcal{l},\mathcal{u})\neq \vac
\right\}
\leq 2.
\]
\end{teo01}
We introduce two types of sets underpinning our arguments in this section.
\begin{def01}
For any positive real numbers $0<L<U$, define
\[
\mfTreb(L,U) \colon= \Cl\left( \left\{z\in \Cx : U^{-1}  \leq \|\iota(z)\|\leq L^{-1}  \right\} \right).
\]
For each $R>0$, put
\[
\mfLT(R)
\colon=
\left\{z\in \Cx: \|\iota(z)\|= R^{-1} \right\}.
\]
\end{def01}

The next properties follow immediately from the definition of $\mfTreb(L,U)$ and $\mfLT(R)$.
\begin{propo01}\label{Propo:BasicProps:Treb:mfLT}
The next assertions hold:
\begin{enumerate}[\rm i.]
\item\label{Propo:BasicProps:Treb:mfLT:01} For every $R>0$, we have $R\, \mfLT(1)=\mfLT(R)$.
\item\label{Propo:BasicProps:Treb:mfLT:02} For all $0<L<U$, we have
\[
\Fr(\mfTreb(L,U))
=
\mfLT(L)
\cup
\mfLT(U).
\]
\item\label{Propo:BasicProps:Treb:mfLT:03} For all $R>0$, we have
\[
\mfLT(R)
=
\Fr\left( 
\Dx\left( \frac{R}{2}; \frac{R}{2}\right)
\cup
\Dx\left( i\frac{R}{2}; \frac{R}{2}\right)
\cup
\Dx\left( -\frac{R}{2}; \frac{R}{2}\right)
\cup
\Dx\left( -i\frac{R}{2}; \frac{R}{2}\right)
\right).
\]
\item\label{Propo:BasicProps:Treb:mfLT:04} For all $n\in\Na$, $\bfa\in \sfR(n)$, and $\mathcal{l},\mathcal{u} \in \Na$ with $2\leq \mathcal{l} < \mathcal{u}$, we have
\[
S_n(\bfa; \mathcal{l}, \mathcal{u}) 
=  
Q_n\left(\omfF_n(\bfa)\cap \mfTreb\left((\mathcal{u} + \tfrac{1}{2})^{-1}, \left(\mathcal{l}-\tfrac{1}{2})^{-1}\right)\right); \bfa\right).
\]
\end{enumerate}
\end{propo01}
We consider cases depending on the form of $\omfF_n(\bfa)$ in the proof of Theorem \ref{TEO:GC:02}:

\noindent{{\bf Case I}: $\omfF_n(\bfa)\equiv \omfF,\, \omfF_1(2+i)\mod \RotaG$.}

\noindent{{\bf Case II}: $\omfF_n(\bfa) \equiv \omfF_1(2),\, \omfF_1(1+i)\mod \RotaG$.}

\subsubsection{\bf Case I} 
In short, we first note that the choice of parameters yields
\[
\bigcup_{\substack{\mathcal{l}\leq \|d\|\leq \mathcal{u} \\ d \in \scD } }   \oclC_{1}(d)
\subseteq 
\omfF_n(\bfa).
\]
Then, we apply the map $Q_n(\,\cdot\,;\bfa)$ on the sets $\mfTreb(\frac{L}{2},2U)$ and $\mfTreb(L,U)$ and conclude. We divide the full argument into lemmas \ref{LEM:GC:02:Case01:01}, \ref{LEM:GC:02:Case01:02}, and \ref{LEM:GC:02:Case01:03}.
\begin{lem01}\label{LEM:GC:02:Case01:01}
There is an absolute constant $\kappa'>0$ such that for all $R>0$ we have
\[
\dist\left( \mfLT(R), \mfLT(2R)\right)
=
\kappa' R.
\]
\end{lem01}
\begin{proof}
The compact sets $\mfLT(1)$ and $\mfLT(\frac{1}{2})$ are disjoint, then $\kappa' \colon= 2\dist\left(\mfLT(1),\mfLT(\frac{1}{2})\right)>0$. These sets are depicted in Figure \ref{Fig:GC02:01:Clovers}. The lemma follows from Proposition \ref{Propo:BasicProps:Treb:mfLT}.\ref{Propo:BasicProps:Treb:mfLT:01}, because for all $R>0$ the map $\Cx\to\Cx$, $z\mapsto Rz$ is a similarity. It is not hard to show that $\kappa'= \sqrt{5}-1$, but we do not need the exact value.
\begin{figure}[h!]
\begin{center}
\includegraphics[scale=0.5,  trim={4.5cm 12.5cm 5cm 4cm},clip]{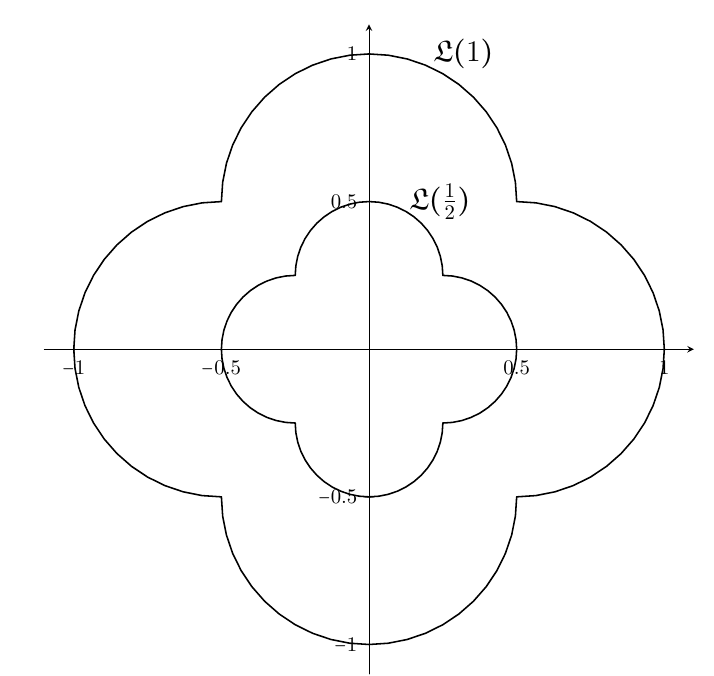}
\caption{ The sets $\mfLT(1)$ and $\mfLT(\frac{1}{2})$. \label{Fig:GC02:01:Clovers}}
\end{center}
\end{figure}
\end{proof}
\begin{lem01}\label{LEM:GC:02:Case01:02}
Let $n\in\Na$, $\bfa\in\sfR(n)$, and $R\in(0,\tfrac{1}{8})$ be arbitrary, and $\kappa'>0$ as in Lemma \ref{LEM:GC:02:Case01:01}, then
\begin{align*}
\frac{\kappa'R}{(1 + 2R)(1+4R)} \frac{1}{|q_n(\bfa)|^2} 
&\leq 
\dist\left(Q_n(\mfLT(R);\bfa),Q_n(\mfLT(2R);\bfa)\right) \\
&\leq 
\frac{\kappa' R}{(1 - 2R)(1-4R)} \frac{1}{|q_n(\bfa)|^2}.
\end{align*}
\end{lem01}
\begin{proof}
Write $q_n=q_n(\bfa)$ and $q_{n-1} = q_{n-1}(a_1,\ldots, a_{n-1})$.  If $\alpha\in \mfLT(R)$ and $\beta\in \mfLT(2R)$, then
\begin{align*}
\left| Q_n(\alpha;\bfa) - Q_n(\beta;\bfa)\right| &= \left| \frac{\alpha}{q_n^2\left( 1 - \frac{q_{n-1}}{q_n} \alpha\right)}
 \,-\, 
\frac{\beta}{q_n^2\left( 1 - \frac{q_{n-1}}{q_n} \beta\right)} \right| \\
&= \frac{1}{|q_n|^2} \frac{|\alpha-\beta|}{\left| \left( 1 - \frac{q_{n-1}}{q_n}\alpha \right) \left( 1 - \frac{q_{n-1}}{q_n}\beta \right)\right|}.
\end{align*}
So, in view of $|\alpha|\leq R$, $|\beta|\leq 2R$ and $|q_{n-1}|<|q_n|$, we have
\[
\frac{1}{|q_n|^2} \frac{|\alpha-\beta|}{\left| \left( 1 + R \right) \left( 1 + 2R \right)\right|}
\leq 
\left| Q_n(\alpha;\bfa) - Q_n(\beta;\bfa)\right|
\leq 
\frac{1}{|q_n|^2} \frac{|\alpha-\beta|}{\left| \left( 1 - R \right) \left( 1 - 2R \right)\right|}.
\]
The result now follows from Lemma \ref{LEM:GC:02:Case01:01}.
\end{proof}
\begin{lem01}\label{LEM:GC:02:Case01:03}
There is a constant $\rho_1'>0$ such that for every $0<L<U< \frac{1}{32}$ satisfying $\frac{1}{4} < \frac{L}{U} < \frac{3}{4}$, all $n\in \Na$, $\bfd\in\sfR(n)$, and $z\in Q_n\left( \mfTreb( L,U); \bfd \right)$, we have
\begin{equation}\label{EQ:GC02:04:04}
\Dx\left( z;\rho_1'\frac{U}{|q_n(\bfd)|^2}\right)
\subseteq 
Q_n\left( \mfTreb\left( \frac{L}{2},2U\right) ; \bfd \right).
\end{equation}
\end{lem01}
\begin{proof}
Let $L,U>0$ be as in the statement. From
\[
\mfTreb(L,U)
\subseteq
\inte\left(\mfTreb\left(\frac{L}{2},2U\right)\right),
\]
we get
\[
Q_n\left(\mfTreb(L,U);\bfa\right)
\subseteq
\inte\left(Q_n\left( \mfTreb\left(\frac{L}{2},2U\right);\bfa\right)\right).
\]
It is, thus, enough to find some $\rho_1'>0$ such that
\[
\dist\left(
\Fr\left( Q_n\left(\mfTreb(L,U);\bfa\right)\right),
\Fr\left( Q_n\left( \mfTreb\left(\frac{L}{2},2U\right);\bfa\right)\right)
\right)
\geq
\rho_1' \frac{U}{|q_n(\bfa)|^2}.
\]
Proposition \ref{Propo:BasicProps:Treb:mfLT}.\ref{Propo:BasicProps:Treb:mfLT:04} implies
\[
\Fr\left( Q_n\left(\mfTreb(L,U);\bfa\right)\right)
=
Q_n( \mfLT(L);\bfa)\cup Q_n( \mfLT(U);\bfa), 
\]
and
\[
\Fr\left( Q_n\left(\mfTreb\left(\frac{L}{2},2U\right);\bfa\right)\right)
=
Q_n\left( \mfLT\left(\frac{L}{2}\right);\bfa\right)\cup Q_n( \mfLT(2U);\bfa ).
\]
We argue as in the proof of Lemma \ref{LEM:GC:02:Case01:02} to conclude
\[
\dist\left(
Q_n( \mfLT(U);\bfa ) , Q_n( \mfLT(2U);\bfa )
\right)
\leq 
\dist\left(
Q_n( \mfLT(U);\bfa ),Q_n\left( \mfLT \left(\frac{L}{2} \right);\bfa \right)
\right)
\]
and
\[
\dist\left(
Q_n ( \mfLT(L);\bfa ), Q_n\left( \mfLT\left(\frac{L}{2} \right)\right)
\right)
\leq 
\dist\left(
Q_n( \mfLT(L);\bfa ),Q_n ( \mfLT ( 2U );\bfa )
\right).
\]
Hence, the result will follow after we prove that
\begin{equation}\label{Eq:GC02:Lem03:01}
\dist\left(
Q_n( \mfLT(U);\bfa ),Q_n( \mfLT(2U);\bfa )
\right)
\geq 
\dist\left(
Q_n( \mfLT(L);\bfa ),Q_n\left( \mfLT\left(\frac{L}{2}\right);\bfa \right)
\right)
\end{equation}
and  show the existence of some $\rho_1'>0$ for which
\begin{equation}\label{Eq:GC02:Lem03:02}
\dist\left(
Q_n( \mfLT(L);\bfa ),Q_n\left( \mfLT\left(\frac{L}{2}\right);\bfa \right)
\right) 
\geq
\rho_1' \frac{U}{|q_n(\bfa)|^2}. 
\end{equation}
We first check \eqref{Eq:GC02:Lem03:01}. The assumptions on $L$ and $U$ lead us to 
\[
\frac{(1-L)(1-2L)}{(1+U)(1+2U)}
>
\frac{(1-U)(1-2U)}{(1+U)(1+2U)}
>
\frac{3}{4}
>
\frac{L}{U}.
\]
So,
\[
\kappa'
\frac{L}{(1-L)(1-2L)}
\frac{1}{|q_n(\bfa)|^2}
<
\kappa'
\frac{U}{(1+U)(1+2U)}
\frac{1}{|q_n(\bfa)|^2}
\]
and we get \eqref{Eq:GC02:Lem03:01} by Lemma \ref{LEM:GC:02:Case01:02}. 

In order to get \eqref{Eq:GC02:Lem03:02}, first observe that $(1+L)(1+2L)\leq 4$, so
\[
\frac{L}{U}
\geq
\frac{1}{4}
\geq
\frac{1}{16}(1+L)(1+2L)
\;\text{ and }\;
\frac{L}{(1+L)(1+2L)}\geq \frac{U}{16}.
\]
The result is now a consequence of Lemma \ref{LEM:GC:02:Case01:02}.
\end{proof}
\begin{rema01}\label{Rem:LEM:GC:02:Case01:03:02}
By Proposition \ref{Prop:SizeEstimates}.\ref{Prop:SizeEstimates:i}, we may shrink $\rho_1'>0$ so that we can replace the last inclusion in Lemma \ref{LEM:GC:02:Case01:02} with
\[
\Dx\left( z;\rho_1 U|\clC_n(\bfd)| \right)
\subseteq 
Q_n\left( \mfTreb\left( \frac{L}{2},2U\right) ; \bfd \right).
\]
\end{rema01}
\begin{rema01}\label{Rem:LEM:GC:02:Case01:03:03}
We will use Lemma \ref{LEM:GC:02:Case01:03} for Case II as well. Assume for a moment that $\omfF_n(\bfd) \equiv \omfF_1(2),  \,\omfF_1(1+i)\mod \RotaG$. Then, Lemma \ref{LEM:GC:02:Case01:03} implies 
\begin{align*}
\Dx\left( z;\rho_1 U|\clC_n(\bfd)| \right)
\cap
\oclC_n(\bfd)
&\subseteq 
Q_n\left( \mfTreb\left( \frac{L}{2},2U\right) ; \bfd \right)
\cap
\oclC_n(\bfd) \\
&=
Q_n\left( \mfTreb\left( \frac{L}{2},2U\right)\cap \omfF_n(\bfd) ; \bfd \right).
\end{align*}
In this case, the disc $\Dx\left( z;\rho_1 U|\clC_n(\bfd)| \right)$ might certainly not be contained in $\oclC_n(\bfd)$. However, it intersects at most one side of $\oclC_n(\bfd)$. 
\end{rema01}
\begin{proof}[Proof of Theorem \ref{TEO:GC:02}: Case I]
Let $\mathcal{l}\in\Na_{\geq 33}$ and $\mathcal{u}\in\Na$ be as in the statement and put 
\[
L= \left(\mathcal{u} + \frac{1}{2} \right)^{-1}
\;\text{ and }\;
U= \left(\mathcal{l} - \frac{1}{2}\right)^{-1}.
\]
Then, we have 
\[
\frac{1}{4}
<
\frac{L}{U}
<
\frac{3}{4}
\;\text{ and }\;
U
< 
\frac{1}{32}.
\]
If $\omfF_n(\bfa)\equiv\omfF, \, \omfF_1(2+i) \mod\RotaG$, then $\mfTreb(L,U)\subseteq \mfFc_n(\bfa)$ and, by Proposition \ref{Propo:BasicProps:Treb:mfLT}.\ref{Propo:BasicProps:Treb:mfLT:04},
\[
S_n(\bfa;\mathcal{l},\mathcal{u})
=
Q_n\left( \mfTreb(L,U);\bfa\right).
\]
We conclude the proof of this case by choosing $\rho'=\rho_1'$ as in Lemma \ref{LEM:GC:02:Case01:03}.
\end{proof}

\subsubsection{\bf Case II}
 We now assume $\omfF_n(\bfa) \equiv \omfF(2), \, \omfF(1+i) \mod \RotaG$. 

\begin{propo01}\label{PROPO:GC02:03}
Let $n\in\Na$ and $\cll,\clu\in \Na$ such that $4<\cll<\clu$. Every $z\in \omfF$ satisfies
\[
\#\left\{ \bfa\in \sfR(n): z\in S_n(\bfa;\cll,\clu)\right\}\leq 2.
\]
\end{propo01}
\begin{proof}
Take any $z\in \omfF$ and assume there are three different words $\bfa,\bfb,\bfc\in \sfR(n)$ such that
\[
z\in 
S_n(\bfa;\cll,\clu)
\cap
S_n(\bfb;\cll,\clu)
\cap
S_n(\bfc;\cll,\clu).
\]
Let $k$ be the largest non-negative integer such that 
\[
\prefi(\bfa;k)
=
\prefi(\bfb;k)
=
\prefi(\bfc;k).
\]
Without loss of generality, we may assume that $a_{k+1}\neq b_{k+1}$. We must have either $c_{k+1}=a_{k+1}$ or $c_{k+1}=b_{k+1}$. Otherwise, we would have
\[
\iota T_k\left(z;(a_1,\ldots, a_k)\right)
\in 
\tau_{a_{k+1}}[\mfF]
\cap 
\tau_{b_{k+1}}[\mfF]
\cap
\tau_{c_{k+1}}[\mfF]
\]
and, hence, $T_{k+1}(z;(a_1,\ldots, a_{k}, a_{k+1})) \in \Esq$. However, since $k+1\leq n$, Proposition \ref{GC:PROPO:02} would forbid $z\in S_n(\bfa;\mathcal{l}, \mathcal{u})$. Assume that $c_{k+1}=b_{k+1}$, so $c_{k+1}\neq a_{k+1}$. Let $j\in\Na$ be the least integer such that $b_{j+k}\neq c_{j+k}$, so $k+j\leq n$. Then, $z'=T_{k+j-1}(z;(b_1,\ldots, b_{k+j-1}))$ satisfies
\[
z'
\in 
\Fr \left(\mfF_{k+j-1}(b_1,\ldots, \ldots, b_{k+j-1}) \right)
\cap 
\oclC_1(b_{k+j})
\cap 
\oclC_1(c_{k+j}).
\]
The inclusion $\Fr \left(\mfF_{k+j-1}(b_1,\ldots, \ldots, b_{k+j-1}) \right)$ follows from $z'\in  \Fr\left( \Fr \left(\mfF_{k}(b_1,\ldots, \ldots, b_{k}) \right)\right)$. Therefore, in view of Lemma \ref{Lemma:GC01:rho:02}.\ref{Lemma:GC01:rho:02:01}, we conclude
\[
T_{k+j}(z;(b_1,\ldots, b_{k+j}))
\in 
\Di_8(\zeta_1)\cup \Esq.
\]
However, this would contradict $z\in S_n(\bfb;\mathcal{l},\mathcal{u})$. $(the partial quotients would be small)$.
\end{proof}
\begin{propo01}\label{PROP:Possibleforms}
Let $n\in\Na$, $\bfa\in\sfR(n)$ be arbitrary, and let $c\in\scD$ be such that $\bfa c\in \sfR(n+1)$. If $c\in\scD$ satisfies $\bfa c\in \sfR(n+1)$, then $\clCc_1(c) \setminus \mfFc_k(\bfa) \neq \vac$ implies $\omfF_n(\bfa)=\omfF_1(2+i) \mod \RotaG$.
\end{propo01}
\begin{proof}
Simply verify four possible forms of $\omfF_n(\bfa)\mod\RotaG$.
\end{proof}

\begin{propo01}\label{PROPO:GC02:06}
Let $\cll,\clu\in\Na$, $\cll< \clu$, be given. Assume that
\[
\left( \ell - \frac{1}{2}\right)^{-1} 
<
\frac{1}{10}
\quad\text{ and }\quad
\frac{1}{4}< \frac{ \cll - \frac{1}{2}}{ \clu + \frac{1}{2}}<\frac{3}{4}. 
\]
There exists a number $\rho_2'>0$ with the following property: for every $k,j\in\Na$, $\bfa\in\sfR(k)$, $\bfd,\bfe\in\sfR(j)$ such that $\bfa\bfd,\bfa\bfe\in\sfR(k + j)$, and
\[
z \in S_{j+k}(\bfa\bfd;\cll, \clu)\cap S_{j+k}(\bfa\bfe;\cll, \clu),
\]
then we have
\[
\Dx\left( T_k(z;\bfa);\frac{\rho_2}{\left(\cll - \frac{1}{2}\right)|q_j(\bfd)|^2}\right)
\subseteq
S_j\left( \bfd; \left\lfloor \frac{\cll-1}{2}\right\rfloor, 2(\clu+1)\right)
\cup
S_j\left( \bfe; \left\lfloor \frac{\cll-1}{2}\right\rfloor, 2(\clu+1)\right).
\]
\end{propo01}
\begin{proof}
The proof is divided into two steps. First, we show that the cylinders $\oclC_{j-1}(d_2,\ldots, d_j)$ and $\oclC_{j-1}(e_2,\ldots, e_j)$ are a reflection of each other. Then, we apply Remark \ref{Rem:LEM:GC:02:Case01:03:03} to obtain the existence of $\rho_2'>0$.

\noindent {\bf Step I.} Let us show that 
\begin{equation}\label{EQ:GC:02:01}
\Mir\left[ \oclC_{j-1}(d_2,\ldots,d_j)\right]
=
\oclC_{j-1}(e_2,\ldots,e_j)
\quad\text{ for some }
\Mir\in\{\Mir_1,\Mir_2\}.
\end{equation}

Define $w\colon=T_k(z;\bfa)$, then 
\begin{equation}\label{EQ:GC:02:02}
w\in 
S_{j}(\bfd;\cll,\clu)
\cap
S_{j}(\bfe;\cll,\clu).
\end{equation}
Pick $b,b'\in\scD$ such that
\begin{equation}\label{EQ:GC:02:03}
\cll \leq \|b\|\leq \clu,
\quad
\cll \leq \|b'\|\leq \clu,
\quad
z\in 
\oclC_{k+j+1}(\bfa \bfd b)
\cap
\oclC_{k+j+1}(\bfa \bfd b').
\end{equation}

From $\iota(w) \in \tau_{e_1}[\omfF]\cap \tau_{d_1}[\omfF]$ we conclude that $1\leq |e_1-d_1|\leq \sqrt{2}$. If we had $|e_1-d_1|=\sqrt{2}$, we would have $w\in\Esq$. However, by Proposition \ref{GC:PROPO:02}, this contradicts \eqref{EQ:GC:02:02}; therefore, we must have $|d_1-e_1|=1$. Without loss of generality, assume that $e_1=d_1+1$.
We claim that
\begin{equation}\label{EQ:GC:02:03:aux}
\tau_{d_1}\left[ \oclC_{1}(d_2) \right]
\cup
\tau_{e_1}\left[ \oclC_{1}(e_2) \right]
\subseteq 
\iota\left[ \omfF_n(\bfa)\right].
\end{equation}
Certainly, by Lemma \ref{Lemma:GC01:rho:02}.\ref{Lemma:GC01:rho:02:01}, Proposition \ref{PROP:Possibleforms} and $e_1=d_1+1$, if we had $\clCc_1(d_2) \setminus \mfFc_n(\bfa)\neq \vac$, we would have $\mfFc_{n+1}(\bfa d_1)=\mfFc_{1}(1+i)$ and $\iota(w)= \tau_{d_1} \, \Mir_1(\zeta_1)$. The proof of Proposition \ref{Prop:Rho:01:01} would then yield
\[
\|d_2\| = \ldots = \|d_j\| = \|b\| = 2, 
\]
which contradicts \eqref{EQ:GC:02:02} and $\cll\geq 10$. The same argument shows $\tau_{e_1}\left[ \oclC_{1}(e_2) \right]\subseteq \iota\left[ \omfF_n(\bfa)\right]$. We, thus, conclude \eqref{EQ:GC:02:03:aux}.

Now, we show that $\oclC_{j-1}(d_2,\ldots, d_j)$ is the single regular cylinder of level $j-1$ whose closure contains $T_{k+1}(z;\bfa d_1)$. Assume there was some $(d_2',\ldots, d_j')\in\sfR(j-1)$ different from $(d_2,\ldots, d_j)$ such that
\[
T_{k+1}(z;\bfa d_1)
\in 
\oclC_{j-1}(d_2',\ldots, d_j')
\cap 
\oclC_{j-1}(d_2,\ldots, d_j).
\]
Then, since $T_{k+1}(z;\bfa d_1)\in \Fr(\mfF)$, an application of Proposition \ref{Lemma:GC01:rho:02}.\ref{Lemma:GC01:rho:02:01} would forbid 
\[
T_{k+1}(z;\bfa d_1)\in S_{j-1}((d_2,\ldots, d_j);\cll,\clu).
\]
Similarly, we may show that $T_{k+1}(z;\bfa e_1)$ belongs to the closure of a single regular cylinder of level $j-1$; namely, $\oclC_{j-1}(e_2,\ldots, e_j)$. 

Observe that, using our assumption $e_1=d_1+1$, we have 
\[
\Mir_1\left(T_{k+1}(z;\bfa d_1)\right)
=
T_{k+1}(z;\bfa e_1).
\]
As a consequence, $T_{k+1}(z;\bfa e_1)\in \Mir_1\left[\oclC_{j-1}(d_2,\ldots,d_j) \right]$ and, because of the uniqueness of $(e_2,\ldots, e_j)$, we conclude \eqref{EQ:GC:02:01}.

\noindent {\bf Step II.}  Next, for  $\cll'=\lfloor \frac{\ell-1}{2}\rfloor$ and $\clu'=2(\clu+1)$, we show that there are absolute constant $\rho_2'>0$ such that the disc
\begin{equation}\label{Eq:Prop:GC:02:02:01}
\Dx\left(\iota(w); \frac{\rho_2'}{\left( \cll'-\frac{1}{2}\right)|q_{j-1}(d_2,\ldots,d_j)|^2} \right)
\end{equation}
is contained in 
\begin{equation}\label{Eq:Prop:GC:02:02:02}
\tau_{d_1}\left[Q_{j-1}\left( \mfTreb(\cll', \clu'); (d_2,\ldots, d_{j-1} \right) \right] \cup
\tau_{e_1}\left[
Q_j\left( \mfTreb(\cll', \clu'); (e_2,\ldots, e_{j-1} \right)\right].
\end{equation}
To this end, note that the side of $\oclC_{j-1}(d_2,\ldots, d_j)$ containing $T_{k+1}(z;\bfa d_1)$ is a segment of $\Fr(\mfF)$. Then, by Remark \ref{Rem:LEM:GC:02:Case01:03:03}, the left half of the disc
\[
\Dx\left(T_{k+1}(z;\bfa d_1); \frac{\rho_1}{\left( \cll'-\frac{1}{2}\right)} \left| \oclC_{j-1}(d_2,\ldots, d_j)\right| \right)
\]
is contained in $\oclC_{j-1}(d_2,\ldots, d_j)$. Similarly, the right half of the disc 
\[
\Dx\left(T_{k+1}(z;\bfa e_1);  \frac{\rho_1}{\left( \cll'-\frac{1}{2}\right)} \left| \oclC_{j-1}(d_2,\ldots, d_j)\right| \right)
\]
is contained in $\oclC_{j-1}(e_2,\ldots, e_j)$. Therefore, the disc in \eqref{Eq:Prop:GC:02:02:01} is contained in the set in \eqref{Eq:Prop:GC:02:02:02} after taking an absolutely small $\rho_2'$ (see Proposition \ref{Prop:SizeEstimates}).
\end{proof}

As in the proof of Lemma \ref{LEM:GC:02:Case01:01}, we may show the existence of some absolute constant $\rho_3'>0$ such that for all $0<L<U<\frac{1}{8}$, $n\in\Na$, and $\bfb\in \sfR(n)$ we have
\[
\Dx\left( z; \rho_3' |S_{k+j}(\bfa\bfd;\cll,\clu)| \right)
\subseteq
S_{j+k}(\bfa\bfe;\cll',\clu')\cap S_{j+k}(\bfa\bfd;\cll',\clu').
\]
This observation combined with Remark \ref{Rem:LEM:GC:02:Case01:03:03} and Proposition \ref{PROPO:GC02:06} yield the next result.
\begin{propo01}\label{PROPO:GC02:07}
There is some $\rho_4'>0$ such that for all $\cll, \clu, \bfa,\bfd,\bfe$ as in Proposition \ref{PROPO:GC02:06}, every $z\in S_{j+k}(\bfa\bfe;\cll,\clu)\cap S_{j+k}(\bfa\bfd;\cll,\clu)$ verifies
\[
\Dx\left( z; \rho_4' |S_{k+j}(\bfa\bfd;\cll,\clu)| \right)
\subseteq
S_{j+k}(\bfa\bfe;\cll',\clu')\cap S_{j+k}(\bfa\bfd;\cll',\clu').
\]
\end{propo01}

\begin{proof}[Proof of Theorem \ref{TEO:GC:02}: Case II]
First, we assume that there are two different $\bfa,\bfb\in\sfR(n)$ for which
\[
z\in S_n(\bfa;\cll,\clu)\cap S_n(\bfb;\cll,\clu).
\]
Take $\rho_4'>0$ as in Proposition \ref{PROPO:GC02:07}, so
\begin{equation}\label{EQ:GC:02:DEM:TEO:01}
\Dx\left( z; \rho_4' |S_{n}(\bfa;\cll,\clu)| \right) 
\subseteq 
S_n(\bfa;\cll',\clu')\cup S_n(\bfb;\cll',\clu').
\end{equation}
Next,  let us show that
\begin{equation}\label{EQ:GC:02:DEM:TEO:02}
\left\{ \bfc\in \sfR(n): S_{n}(\bfc;\cll,\clu)\cap \Dx(z; \rho_3'|S_n(\bfa;\cll,\clu)|)\right\} 
\subseteq 
\{\bfa,\bfb\}.
\end{equation}
If this was not true, then there would be some $w\in \Dx(z;2\rho|S_n(\bfa;\cll,\clu)|)$ and $\bfc\in \sfR(n)\setminus\{\bfa,\bfb\}$ such that $w\in S_n(\bfc;\cll,\clu)$. Take $c\in\scD$ with $\cll\leq \|c\|\leq \clu$ and $w\in \oclC_{n+1}(\bfc c)$. Since $\oclC_{n+1}(\bfc c)=\Cl(\clCc_{n+1}(\bfc c))$, we would have
\[
\Dx(z;\rho_4' |S_n(\bfa;\cll,\clu)|) \cap \clCc_{n+1}(\bfc c)\neq \vac;
\]
however, since $\bfc\neq \bfa$ and $\bfc\neq \bfb$, we know that
\[
\clCc_n(\bfc)\cap \left( S_n(\bfa;\cll',\clu')\cup S_n(\bfb;\cll',\clu')\right)
\subseteq
\clCc_n(\bfc)\cap \left( \oclC_n(\bfa)\cup \oclC_n(\bfb)\right) =\vac,
\]
which is a contradiction. Therefore, \eqref{EQ:GC:02:DEM:TEO:02} holds and the theorem has been proven in this subcase.

Now assume that $\bfa$ is the only element $\bfc\in\sfR(n)$ for which $z\in S_n(\bfc;\ell,\clu)$ holds. If the disc $\Dx\left( z;  \frac{\rho_4'}{2} |S_{n}(\bfa;\cll,\clu)|\right)$ only intersects one set of the form $S_n(\bfc;\cll,\clu)$, we are done. Otherwise, there would be some $\bfb\in\sfR(n)$, $\bfa\neq \bfb$, and some  $w\in S_n(\bfa;\cll,\clu)\cap S_n(\bfb;\cll,\clu)$. We may, thus, use the previous subcase along with 
\[
\Dx\left( z; \frac{\rho_4'}{2} |S_{n}(\bfa;\cll,\clu)| \right) 
\subseteq 
\Dx\left( w; \rho_4' |S_{n}(\bfa;\cll,\clu)| \right) 
\subseteq 
S_n(\bfa;\cll',\clu')\cup S_n(\bfb;\cll',\clu')
\]
to conclude the theorem.
\end{proof}


\section{Approaching the irregular by the regular}\label{SEC:IRRREG}
So far, our discussion has revolved around regular numbers and their associated objects (cylinders, prototype sets, words). However, as noted before, there are uncountably many irregular numbers. Theorem \ref{TEO:REG:IRREG:CLOS} helps us understand why we can overlook irregular numbers in the proof of Theorem \ref{TEO:DIMHB}.

\begin{teo01}\label{TEO:REG:IRREG:CLOS}
For every $n\in\Na_{\geq 2}$, $\bfa=(a_1,\ldots,a_n)\in\sfIr(n)$, and $f\in\Di_{8}$, there is some $\bfb=(b_1,\ldots,b_n)\in\sfR(n)$ such that 
\[
f\left[ \clC_n(\bfa)\right]\subseteq \oclC_n(\bfb)
\quad \text{ and } \quad
\|a_j\|=\|b_j\| 
\;
\text{ for all } j\in\{1,\ldots,n\}.
\]
\end{teo01}
\begin{coro01}\label{CoroShiftSpace}
There is a single continuous extension $\overline{\Lambda}$ of $\Lambda|_{\sfR}$ to the closure $\overline{\sfR}$ of $\sfR$. Moreover, $\overline{\Lambda}[\overline{\sfR}] = \omfF\setminus \QU(i)$.
\end{coro01}
Using a slightly different notation, Theorem \ref{TEO:REG:IRREG:CLOS} and Corollary \ref{CoroShiftSpace} are proven in \cite{GRGRH24}. Moreover, that paper provides a procedure to obtain $\bfb$. We, thus, omit the proofs and refer the reader to \cite[Subsection 7.2]{GRGRH24}. 

\begin{rema01}
We can deal with irregular numbers in a different way.  It is not hard to show that the set of irregular numbers has Hausdorff dimension $1$. To see this, note that an irregular cylinder is either a line segment or an arch and there are countably many of them. Hence, we may ignore them while proving Theorem \ref{TEO:MAIN} and Theorem \ref{TEO:DIMHB}. The downside of this strategy is the loss of Corollary \ref{CoroShiftSpace}.
\end{rema01}

\section{Proof of Theorem \ref{TEO:DIMHB}} \label{Sec:TeodimBprf}

We start the proof of Theorem \ref{TEO:DIMHB} by showing the existence of the limit in the statement. Then, we divide the proof into two parts: the upper bound and the lower bound. For the upper bound, we use a natural covering argument. For the lower bound, we apply the Mass Distribution Principle (Lemma \ref{LE:MDP} below).

\subsection{The functions $g_n$}\label{Sec:AuxFun}
The numbers $N_1\colon=11$ and $\psi\colon=\sqrt{\frac{1+\sqrt{5}}{2}}$ will remain fixed throughout this section. Write
\[
\forall M\in\Na_{\geq 4}
\qquad
\scD(M)
\colon=\
\{a\in\scD:\|a\|\leq M\}.
\]

\begin{def01}
Given $B>1$, $n\in\Na$, and $\scA\subseteq \scD$, define 
\[
\forall \rho\in\RE_{>0}
\qquad
g_n(\rho;\scA,B)
\colon=
 \sum_{\bfa\in \scA^n\cap \sfR(n)}
  \left(\frac{ 5\sqrt{2} }{B^n\|q_n(\bfa)\|^2}\right)^{\rho}.
\]
\end{def01}
When $\scA$ and $B$ are fixed, we write $g_n(\rho)$ instead of $g_n(\rho;\scA,B)$.

\begin{propo01}\label{Prop:PropertiesOfgn}
Take $B>1$, $n\in \Na_{\geq N_1}$, and let $\scA$ be either $\scE$, $\scD(M)$ for some $M\geq 4$ or $\scD$.
\begin{enumerate}[\rm i.]
\item \label{Prop:PropertiesOfgn:i} The function $\rho\mapsto g_n(\rho;\scA,B)$ is non-increasing.
\item \label{Prop:PropertiesOfgn:ii} If $\scA'$ is  $\scE$, $\scD(M)$ for some $M\geq 4$ or $\scD$, and  $\scA\subseteq \scA'$; then, $g_n(\rho;\scA,B)\leq g_n(\rho;\scA',B)$.
\end{enumerate} 
\end{propo01}
\begin{proof}
\begin{enumerate}[\rm i.]
\item From  $|q_n(\bfa)|> \psi^{n-1} \geq \psi^{10}>10$, we have $5\sqrt{2}< \|q_n(\bfa)\|$, hence $\frac{5\sqrt{2}}{B^n\|q_n(\bfa)\|^2}<1$. The definition of $g_n$ then tells us that $0<\rho<\rho'$ implies $g_n(\rho';\scA,B)\leq g_n(\rho;\scA,B)$.
\item The proof is straightforward.
\end{enumerate}
\end{proof}
 
For any $n\in\Na$, $B>1$, and $\scA$ equal to $\scE$, or $\scD(M)$ for some $M\geq 4$, or $\scD$, we define
\[
s_{n,B}(\scA)\colon=\inf\{\rho>0: g_n(\rho;\scA,B)\leq 1\}.
\]
For notational ease, for any $j\in\Na$, $B>1$, $\rho>0$, and $\bfa\in\Omega(j)$, put
\[
\tilde{h}_j(\bfa;\rho)\colon= (B^j\|q_j(\bfa)\|^2)^{-\rho}.
\]
The next estimates follow from Proposition \ref{Prop:SizeEstimates}.\ref{Prop:SizeEstimates:ii}. 
\begin{propo01}\label{Propo:BndsOnFunctions-h_j}
Take $B>1$. For all $m,n\in\Na$, $\bfa\in\Omega(n)$, $\bfb\in\Omega(m)$ such that $\bfa\bfb\in \Omega(m+n)$, and $\rho>0$, we have
\[
\frac{1}{6^{2\rho}} \tilde{h}_n(\bfa;\rho)\tilde{h}_m(\bfb;\rho)
\leq
\tilde{h}_{n+m}(\bfa\bfb;\rho)
\leq
(5\sqrt{2})^{2\rho} \tilde{h}_n(\bfa;\rho)\tilde{h}_m(\bfb;\rho).
\]
\end{propo01}

\begin{propo01}\label{Propo:gn+m_leq_gngm}
Let $\scA$ be either $\scD$, $\scD(M)$ for some $M\geq 4$, or $\scE$, and $B>1$. For every $m,n\in\Na$ and $\rho>0$, we have
\[
g_{n+m}(\rho)\leq g_{n}(\rho)g_{m}(\rho).
\]
\end{propo01}
\begin{proof}
For $m,n\in\Na$ and $\rho>0$, we have 
\begin{align*}
g_{n+m}(\rho) &= \sum_{\bfc\in \scA^n \cap \sfR(n) } (5\sqrt{2})^{2\rho} \tilde{h}_{n+m} (\bfc;\rho) \\
&\leq \sum_{\substack{\bfa\in \scA^n \cap \sfR(n) \\ \bfb\in \scA^m \cap \sfR(m) }} (5\sqrt{2})^{2\rho} \tilde{h}_{m} (\bfa;\rho) (5\sqrt{2})^{2\rho} \tilde{h}_{m} (\bfb;\rho) \\
&\leq  \left(\sum_{ \bfa\in \scA^n \cap \sfR(n) } (5\sqrt{2})^{2\rho} \tilde{h}_{n} (\bfa;\rho)\right)\left(\sum_{ \bfb\in \scA^m \cap \sfR(m)} (5\sqrt{2})^{2\rho} \tilde{h}_{m} (\bfb;\rho)\right)\\
&= g_n(\rho)g_m(\rho).
\end{align*}
\end{proof}
\begin{propo01}\label{Propo:gngm-LB}
Let $\scA$ be either $\scD$, $\scD(M)$ for some $M\geq 4$, or $\scE$, and $B>1$. Write $\kappa_4\colon= \frac{3\sqrt{2}+1}{\sqrt{2}-1}$. For every $m,n\in\Na$ and $\rho>0$, the next estimate holds:
\[
\frac{1}{1800^{\rho}  + 20(3600\kappa_4^2)^\rho} g_n(\rho)g_m(\rho
)\leq g_{n+m}(\rho).
\]
\end{propo01} 
\begin{proof}
Take any $m,n\in\Na$ and $\rho>0$. Then,
\begin{align}
g_m(\rho)g_n(\rho) &=  g_m(\rho)(5\sqrt{2})^{2\rho} \sum_{\bfa\in \scA^n\cap\sfR(n)} \tilde{h}_n(\bfa;\rho) \nonumber \\
&= g_m(\rho) (5\sqrt{2})^{2\rho} \sum_{\substack{\bfa\in\scA^n\cap\sfR(n) \\ \|a_n\|\geq 3}} \tilde{h}_n(\bfa;\rho) + g_m(\rho) (5\sqrt{2})^{2\rho} \sum_{\substack{\bfa\in\scA^n\cap\sfR(n) \\ \|a_n\|\leq 2}} \tilde{h}_n(\bfa;\rho).\label{Eq:SeparaProd}
\end{align}

We now bound separately each of the two addends in \eqref{Eq:SeparaProd}. By Proposition \ref{Propo:BndsOnFunctions-h_j}, the first term can be bounded as follows:
\begin{align*}
(5\sqrt{2})^{2\rho}\left(\sum_{\substack{\bfa\in \scA^n\cap\sfR(n) \\ \|a_n\|\geq 3}} \tilde{h}_n(\bfa;\rho)\right)g_m(\rho)  
&= (5\sqrt{2})^{4\rho} \left(\sum_{\substack{\bfa\in\scA^n\cap\sfR(n) \\ \|a_n\|\geq 3}} \tilde{h}_n(\bfa;\rho)\right) \left(\sum_{\bfb\in\scA^m\cap\sfR(m)} \tilde{h}_m(\bfb;\rho)\right) \\
& \leq  (5\sqrt{2})^{4\rho} 6^{2\rho} \sum_{\substack{\bfa\in \scA^n\cap   \sfAF(n)\\ \bfb\in\scA^m\cap\sfR(m)}} \tilde{h}_{n+m}(\bfa\bfb;\rho) \\
&\leq (5\sqrt{2})^{2\rho} 6^{2\rho} \sum_{\bfc\in \scA^{n+m}\cap\sfR(n+m)(\scA) } (5\sqrt{2})^{2\rho} \tilde{h}_{n+m}(\bfc;\rho);
\end{align*}
that is
\begin{equation}\label{Eq:PropComparegn+m-02:01}
(5\sqrt{2})^{2\rho}\left(\sum_{\substack{\bfa\in\scA^n\cap \sfR(n) \\ \|a_n\|\geq 3}} \tilde{h}_n(\bfa;\rho)\right)g_m(\rho)
\leq 
1800^{\rho} g_{n+m}(\rho).
\end{equation}
We need a more intricate argument for the second term in \eqref{Eq:SeparaProd}. For each $\bfa\in\scA^n\cap\sfR(n)$ such that $\|a_n\|\leq 2$, choose $d(\bfa)\in\scD$ satisfying
\[
\Pm(d(\bfa))=\|d(\bfa)\|=3 
\quad\text{ and }\quad
\bfx(\bfa)\colon= (a_1,\ldots,a_{n-1},d)\in \scA^n \cap \sfF(n),
\]
(see Proposition \ref{PROP:SQRT8-REG:NUVO}). Then, we have
\[
\left|\frac{q_n(\bfx(\bfa))}{q_n(\bfa)}\right| =
\left| \frac{d(\bfa)q_{n-1}(\bfa)+q_{n-2}(\bfa)}{a_{n}q_{n-1}(\bfa)+q_{n-2}(\bfa)} \right| 
= \left| \frac{d(\bfa) + \frac{q_{n-2}(\bfa)}{q_{n-1}(\bfa)}}{a_{n} + \frac{q_{n-2}(\bfa)}{q_{n-1}(\bfa)}} \right|\leq \frac{|d(\bfa)|+1}{|a_n|-1} \leq \frac{3\sqrt{2}+1}{\sqrt{2}-1} = \kappa_4,
\]
which implies
\[
\|q_n(\bfx(\bfa))\| \leq  |q_n(\bfx(\bfa))| \leq \kappa_4 |q_n(\bfa)| \leq \kappa_4 \sqrt{2}\|q_n(\bfa)\|.
\]
Hence
\begin{equation}\label{Eq:BndTrm}
\frac{1}{(B^n\|q_n(\bfa)\|^2)^{\rho}}
\leq 
\frac{ (2\kappa_4^2)^{\rho} }{(B^n\|q_n(\bfx(\bfa))\|^2)^{\rho}}.
\end{equation}

The function $\bfa\mapsto \bfx(\bfa)$ might not be injective, but it is at most twenty-to-one. Indeed, for any $\bfb=(b_1,\ldots,b_{n-1},b_n)\in \scA^n\cap \sfF(n)$ such that $\Pm(b_n)=\|b_n\|=3$, the set of words $\bfa=(a_1,\ldots,a_n)\in \scA^n\cap\sfR(n)$ satisfying $\|a_n\|\leq 2$ and $\bfx(\bfa)=\bfb$ is contained in the set
\[
\{\bfa\in \scA^n\cap\sfR(n): (a_1,\ldots,a_{n-1})=(b_1,\ldots,b_{n-1}), \|a_n\|\leq 2\},
\]
which has at most $20$ elements, because $\#\{a\in\scD: \|a\|\leq 2\}=20$. This observation and \eqref{Eq:BndTrm} give
\begin{align*}
\sum_{\substack{\bfa\in\scA^n\cap\sfR(n) \\ \|a_n\|\leq 2}} \tilde{h}_n( \bfa;\rho) &\leq
(2\kappa_4^2)^{\rho} \sum_{\substack{\bfa\in\scA^n\cap\sfR(n) \\ \|a_n\|\leq 2}} \tilde{h}_n(\bfx(\bfa);\rho) \\
&\leq 20 (2\kappa_4^2)^{\rho} \sum_{\substack{\bfc\in\scA^n\cap\sfR(n) \\ \Pm(c_n)=\|c_n\|= 3}} \tilde{h}_n(\bfc;\rho) \\
&\leq 20 (2\kappa_4^2)^{\rho} \sum_{\bfc\in \scA^n \cap\sfF(n)} \tilde{h}_n(\bfc;\rho).
\end{align*}
Then, by Proposition \ref{Propo:BndsOnFunctions-h_j}, we have
\begin{align*}
(5\sqrt{2})^{2\rho} \left( \sum_{\substack{\bfa\in \scA^n\cap\sfR(n) \\ \|a_n\|\leq 2}} \tilde{h}_n( \bfa;\rho) \right) g_m(\rho) &\leq (5\sqrt{2})^{4\rho} 20 (2\kappa_4^2)^{\rho} \left(  \sum_{\bfa\in \scA^n\cap\sfF(n)} \tilde{h}_n(\bfa;\rho) \right)\left(\sum_{\bfb\in\scA^m\cap\sfR(m)} \tilde{h}_m(\bfb;\rho) \right) \nonumber\\
&\leq 
(5\sqrt{2})^{2\rho}  20 (2\kappa_4^2)^{\rho} 6^{2\rho}
\sum_{\bfc\in\scA^{n+m}\cap\sfR(n+m)  }(5\sqrt{2})^{2\rho} \tilde{h}_{n+m}(\bfc;\rho) \\
&=
 20 (3600\kappa_4^2)^{\rho} g_{n+m}(\rho),
\end{align*}
and the proof is done.
\end{proof}
 
\subsection{The dimensional number $s_B$}
In this subsection, we give an alternative definition of the number $s_B$. After establishing some relevant properties, we show that this definition is equivalent to the one in \eqref{Eq:Defsb}.

\begin{propo01}\label{Le:LemPrel-ENT}
If $p,q\in\Na$ are coprime and $p<q$, then 
\[
\Na_{\geq p(q-1)} \subseteq \{j\in \Na: j=pk+tq \text{ for some } k,t\in\Na_0\}.
\]
\end{propo01}
\begin{proof}
See \cite[Section 3. Problem 14]{AndAndFen2007}.
\end{proof}
\begin{lem01}\label{Le:LemPrel-01}
Let $B>1$ be arbitrary. If $\scA=\scD$, $\scA=\scE$, or $\scA=\scD(M)$ for some $M\geq 4$, then the following limit exists
\[
s_B(\scA) \colon = \lim_{n\to\infty} s_{n,B}(\scA).
\]
\end{lem01}
\begin{proof}
For clarity, we divide the proof into six steps.
\begin{enumerate}[\rm i.]
\item We show that for $n,t\in\Na_{\geq N_1}$ and $r, k\in\Na$ we have
\begin{equation}\label{Eq:Bnd_s_m+n}
s_{rn+tk,B}(\scA)
\leq
\max\{s_{rn,B}(\scA),s_{tk,B}(\scA)\} 
\leq 
\max \{ s_{n,B}(\scA),s_{t,B}(\scA)\}.
\end{equation}
To this end, take any $m,n\in\Na_{\geq N_1}$. We obtain \eqref{Eq:Bnd_s_m+n} from the stronger inequality 
\begin{equation}\label{Eq:LeTech01Eq01}
s_{m+n,B}(\scA) \leq \max\{ s_{m,B}(\scA),s_{n,B}(\scA)\}.
\end{equation}
Take $\veps>0$ and $\rho=\max\{s_{m,B}(\scA),s_{n,B}(\scA)\}+\veps$. Then, $g_m(\rho)<1$, $g_n(\rho)<1$ and, by Proposition \ref{Propo:gn+m_leq_gngm},  $g_{n+m}(\rho)\leq g_n(\rho)g_m(\rho)<1$. Therefore, 
\[
s_{m+n,B}(\scA) <  \max\{s_{m,B}(\scA), s_{n,B}(\scA)\} + \veps,
\]
which proves \eqref{Eq:LeTech01Eq01}.
\item \label{Parteii} Define 
\[
\overline{s}\colon=\limsup_{n\to\infty} s_{n,B}(\scA).
\]
We claim that 
\begin{equation}\label{Eq:Le01-Claim}
\#\{ n\in\Na:s_{n,B}(\scA)\geq \overline{s} \} = \infty.
\end{equation}
If \eqref{Eq:Le01-Claim} were false, there would be some $N_0\in\Na_{\geq N_1}$ satisfying
\[
s_{n,B}(\scA)<\overline{s}
\;\text{ for all } n\in\Na_{\geq N_0}.
\]
Define
\[
\overline{r}\colon= \max\{s_{N_0+r,B}\colon r\in\{0,1,\ldots,N_0 - 1\}\}< \overline{s}.
\]
Given $n\in\Na_{\geq 2N_0}$, let $q\in\Na_{\geq 2}$ and $r\in\{0,1,\ldots,N_0-1\}$ be such that $n=qN_0+r$, then
\begin{align*}
s_{n,B}(\scA) &= s_{(q-1)N_0 + (N_0+r),B}(\scA)\\
&\leq \max\{ s_{(q-1)N_0,B}(\scA),s_{N_0+r,B}(\scA)\} &&\text{(by \eqref{Eq:Bnd_s_m+n})} \\
&\leq \max \{s_{(q-1)N_0,B}(\scA),\overline{r}\} \\
&< \overline{s}.
\end{align*}
However, this leads us to the contradiction
\[
\overline{s} =\limsup_{n\to\infty} s_{n,B}(\scA)\leq \overline{r} <\overline{s},
\]
so \eqref{Eq:Le01-Claim} must hold.
\item \label{Parteiii} Observe that 
\begin{equation}\label{Eq:Le01Claim02}
\forall n\in\Na_{\geq N_1(N_1+1)}
\qquad
s_{n,B}(\scA)\leq \max\{s_{N_1,B}(\scD), s_{N_1+1,B}(\scD)\}.
\end{equation}
Certainly, for any $n\in\Na_{\geq N_1(N_1+1)}$, Proposition \ref{Prop:PropertiesOfgn}.\ref{Prop:PropertiesOfgn:ii}, \eqref{Eq:Bnd_s_m+n}, and Proposition \ref{Le:LemPrel-ENT} give
\[
s_{n,B}(\scA)
\leq
s_{n,B}(\scD)
\leq 
\max\{s_{N_1,B}(\scD), s_{N_1+1,B}(\scD)\}.
\]

\item \label{Parteiv} Write 
\[
\hat{s}\colon=\max\{s_{N_1,B}(\scD),s_{N_1+1,B}(\scD)\}, \quad
\kappa(\hat{s}) \colon= (1800^{\hat{s}}  + 20(3600\kappa_4^2)^{\hat{s}})^{-1}.
\]
We now show that
\begin{equation}\label{Eq:Le01Claim02.5}
\forall\delta>0
\qquad
g_n(s_{jn,B}(\scA)+\delta)\leq \frac{1}{\kappa(\hat{s}) ^{\frac{1}{2}}}.
\end{equation}
By Proposition \ref{Propo:gngm-LB}, we have
\begin{equation}\label{Eq:Le01Claim02:s:hat}
\forall \rho\in (0,\hat{s}) 
\quad 
\forall m, n\in\Na 
\qquad
\kappa(\hat{s})  g_{n}(\rho)g_{m}(\rho)\leq g_{n+m}(\rho).
\end{equation}
Hence, if $0<\rho<\hat{s}$, $j\in\Na_{\geq 2}$, and $n\in\Na$, we obtain $g_{jn}(\rho) \geq \kappa(\hat{s})^{j-1} g_n(\rho)^j$ and
\[
g_n(\rho)
\leq 
\frac{1}{\kappa(\hat{s}) ^{1-\frac{1}{j}}} g_{jn}(\rho)^{\frac{1}{j}} 
\leq 
\frac{1}{ \kappa(\hat{s})^{\frac{1}{2}}} g_{jn}(\rho)^{\frac{1}{j}}.
\]
Using the definition of $s_{nj,B}(\scA)$, we conclude \eqref{Eq:Le01Claim02.5}.

\item Let us show that
\begin{equation}\label{Eq:Le01Claim03}
\forall r\in\Na_{\geq N_1(N_1+1)}
\qquad
s_{r,B}(\scA)\geq \overline{s}.
\end{equation}
Assume that $r\in \Na_{\geq N_1(N_1+1)}$ satisfies $s_{r,B}(\scA)< \overline{s}$. Consider $\veps\colon=\overline{s}- s_{r,B}(\scA)>0$. By \eqref{Eq:Le01-Claim}, we may choose some $p\in\Na$ such that 
\[
s_{p,B}(\scA)\geq \overline{s} 
\quad\text{ and }\quad
\left(\frac{ B^p\psi^{p-1} }{5\sqrt{2}}\right)^{\frac{\veps}{2}} > \kappa(\hat{s}) ^{\frac{1}{2}}.
\]
Then, taking $\delta=\tfrac{\veps}{2}$ in \eqref{Eq:Le01Claim02.5}, we have
\begin{align*}
g_p\left( s_{rp,B}(\scA) + \veps \right)  &=  \sum_{\bfa\in \scA^p\cap\sfR(p) }  \frac{(5\sqrt{2})^{2s_{rp,B}(\scA) + \veps}}{(B^p\|q_p(\bfa)\|^2)^{s_{rp,B}(\scA)+\veps }} \\
&\leq \left(\frac{5\sqrt{2}}{ B^p\psi^{p-1} }\right)^{\frac{\veps}{2}} \sum_{\bfa\in  \scA^p\cap\sfR(p) }  \frac{(5\sqrt{2})^{2s_{rp,B}(\scA) + \frac{\veps}{2} }}{(B^p\|q_p(\bfa)\|^2)^{ s_{rp,B}(\scA)+\frac{\veps}{2} }}  &&\text{(by Proposition \ref{Propo3.2}.\ref{Propo3.2_ii})}\\
&=\left(\frac{5\sqrt{2}}{ B^p\psi^{p-1} }\right)^{\frac{\veps}{2}} g_p\left(s_{rp,B}(\scA) + \frac{\veps}{2} \right) \\
&< 1.
\end{align*}
Thus, 
\begin{equation*}\label{Eq:Le01-Eqa}
s_{p,B}(\scA)\leq s_{rp,B}(\scA)+\frac{\veps}{2}
\end{equation*}
and, by \eqref{Eq:Bnd_s_m+n},
\begin{equation*}\label{Eq:Le01-Eqb}
s_{p,B}(\scA)\geq \overline{s} = s_{r,B}(\scA)+\veps \geq s_{rp,B}(\scA)+\veps.
\end{equation*}
These inequalities imply
\[
\frac{\veps}{2} + s_{rp,B}(\scA) 
\geq 
s_{p,B}(\scA) 
\geq 
s_{rp,B}(\scA) + \veps,
\;\text{ so }\; 0>\frac{\veps}{2}.
\]
In view of this contradiction, \eqref{Eq:Le01Claim03} must hold. 
\item We conclude the result, because \eqref{Eq:Le01Claim03} yields 
\[
\liminf_{n\to\infty} s_{n,B}(\scA) 
\geq 
\overline{s}
= 
\limsup_{n\to\infty} s_{n,B}(\scA).
\]
\end{enumerate}
\end{proof}
For any $p\in\Na_{\geq N_1}$, write $s_{p,B}\colon= s_{p,B}(\scD)$ and $s_{p,B}(M)\colon=s_{p,B}(\scD(M))$ for all $M\in\Na_{\geq 4}$.
\begin{lem01}\label{Le:LemPrel-02}
If $p\in\Na_{\geq N_1}$, then $\displaystyle\lim_{M\to\infty} s_{p,B}( M )=s_{p,B}$.
\end{lem01}
\begin{proof} 
By Proposition \ref{Prop:PropertiesOfgn}.\ref{Prop:PropertiesOfgn:ii}, the map $M\mapsto s_{p,B}( M )$ is non-decreasing for $M\geq 4$ and bounded above by $s_{p,B}$. Then, the next limit exists:
\[
t_{p,B}\colon= \lim_{M\to\infty} s_{p,B}(M)\leq s_{p,B}.
\]
The definition of $s_{p,B}$ and $p\geq N_1$ give
\[
\forall \rho\in \RE_{>0}
\qquad
g_p(\rho;\scD,B) = \lim_{M\to\infty} g_p(\rho;\scD(M),B).
\]
Take $\veps>0$. Then, by $g_p(s_{p,B}-\veps;\scD, B)>1$ and the above limit, every sufficiently large $M$ satisfies
$
g_p(s_{n,B} - \veps, \scD(M), B)>1,
$
so $s_{p,B}-\veps\leq s_{p,B}( M )\leq t_{p,B}$ and $s_{p,B}\leq t_{p,B}$.
\end{proof}
Write $s_B\colon= s_B(\scD)$ and $s_B(M)\colon= s_B(\scD(M))$ for all $M\in\Na_{\geq 2}$.
\begin{lem01}\label{Le:LemPrel-03}
$\displaystyle \lim_{M\to\infty} s_{B}( M )=s_{B}$.
\end{lem01}
\begin{proof}
Take $\veps>0$. Let $\widetilde{N} = \widetilde{N}(B,\veps) \in\Na_{\geq N_1}$ be such that 
\[
\kappa(\hat{s})^{\frac{1}{2}}
< 
\left( \frac{B^{ \widetilde{N} } \psi^{ \widetilde{N}  }}{5\sqrt{2}}\right)^{\frac{\veps}{2}}.
\]
Consider any $n\in\Na_{\geq \widetilde{N}}$. As in part \ref{Eq:Le01-Eqa} in the proof of Lemma \ref{Le:LemPrel-01}, we conclude that if $\scA = \scD$, $\scA=\scE$, or $\scA=\scD(M)$, $M \geq 4$, then every $m\in\Na$ satisfies
\[
s_{n,B}(\scA)\leq s_{mn,B}(\scA) + \frac{\veps}{2}.
\]
Moreover, \eqref{Eq:Bnd_s_m+n} gives $s_{mn,B}(\scA)\leq s_{n,B}(\scA)$, so
\[
|s_{n,B}(\scA) - s_{mn,B}(\scA)|<\veps.
\]
We let $m\to\infty$ and use Lemma \ref{Le:LemPrel-01} to get
\begin{equation}\label{Eq:Lema27-02}
|s_{n,B}(\scA) - s_{B}(\scA)|<\veps.
\end{equation}
Also, by Lemma \ref{Le:LemPrel-02}, every sufficiently large $M$ satisfies
\begin{equation}\label{Eq:Lema27-03}
|s_{n,B}( M ) - s_{n,B}|< \veps.
\end{equation}
Inequalities \eqref{Eq:Lema27-02} and \eqref{Eq:Lema27-03} lead us to the desired conclusion:
\[
|s_B( M ) - s_B| 
\leq |s_B( M ) - s_{n,B}( M )| + 
|s_{n,B}( M ) - s_{n,B}| + 
|s_{n,B}-s_B|< 3\veps.
\]
\end{proof}
For $B>1$ and $\kappa_1$ as in Proposition \ref{Prop:SizeEstimates}.\ref{Prop:SizeEstimates:i}, define
\[
N_2
\colon= 
\max\left\{ N_1, \frac{\log(50) - \log\pi - \log \kappa_1}{2\log B}+1 \right\}.
\]

\begin{propo01}\label{Propo:Bnd_snBM}
For any $B>1$ and $M\geq 4$, we have
\[
\sup\left\{ s_{n,B}( M ): n\in\Na_{\geq N_2}, M\in\Na \right\}\leq 2.
\]
\end{propo01}
\begin{proof}
Take $n\in\Na_{\geq N_2}$. For any $M\geq 4$, we have $s_{n,B}( M )\leq s_{n,B}$ by Proposition \ref{Prop:PropertiesOfgn}. Moreover, $s_{n,B}<2$ holds, because of
\[
g_n(2;\scD,B) 
= \sum_{\bfa\in\sfR(n)} \left( \frac{5\sqrt{2}}{B^n\|q_n(\bfa)\|^2}\right)^2 
= \frac{50}{B^{2n}} \sum_{\bfa\in\sfR(n)} \frac{ 1 }{\|q_n(\bfa)\|^4} 
\leq \frac{50}{B^{2n}\kappa_1\pi} \sum_{\bfa\in\sfR(n)}  \leb(\clC_n(\bfa)) \\
<1.
\]
\end{proof}
The next estimate can be shown inductively on $n$.
\begin{propo01}\label{Propo:Cota_qn}
If $n\in\Na$ and $\bfa=(a_1,\ldots, a_n)\in\sfR(n)$, then
\[
\frac{1}{\sqrt{2}}\prod_{j=1}^n (|a_j| - 1)
\leq
\|q_n(\bfa)\|
\leq 
\prod_{j=1}^n (|a_j| + 1).
\]
\end{propo01}

\begin{lem01}\label{Le:LemPrel-04}
The function $B\mapsto s_B$, mapping $\RE_{>0}$ to $\RE$, has the following properties:
\begin{enumerate}[\rm i.]
\item \label{Le:LemPrel-04-03} it is continuous.
\item \label{Le:LemPrel-04-01} $\displaystyle\lim_{B\to 1} s_B=2$,
\item \label{Le:LemPrel-04-02} $\displaystyle\lim_{B\to \infty} s_B=1$,
\end{enumerate}
\end{lem01}
\begin{proof}
\begin{enumerate}[\rm i.]
\item 
Let $B_0>1$ be given. Let $\veps\in(0,1)$ be arbitrary. Put
\[
\delta_1\colon= \frac{1}{2}\left( B_0 - B_0^{\frac{1}{1+\frac{\veps}{4}}} \right), 
\quad 
\delta_2 \colon= \frac{1}{2}\left( B_0^{ 1 + \frac{\veps}{12}} - B_0 \right), 
\quad 
\delta\colon= \min\{\delta_1,\delta_2\}
\]
and
\[
N
\colon=
\max\left\{
N_2, \;
\frac{\log(5\sqrt{2})}{2\left( \left(1+\frac{\veps}{4}\right) \log(B_0 - \delta) - \log B_0\right)}, \;
\frac{\log(5\sqrt{2})}{2\left( \left(1+\frac{\veps}{4}\right) \log B_0 - \log (B_0 + \delta) \right)}
\right\}.
\]
We chose the parameters for the next estimates to hold:
\[
\forall n\in\Na_{\geq N}
\qquad
\left(\frac{B_0}{(B_0 - \delta)^{1+\frac{\veps}{4} }}\right)^{2n} 
< \frac{1}{5\sqrt{2}}, 
\quad
\left(\frac{B_0 + \delta}{B_0^{1+\frac{\veps}{4} }} \right)^{2n} 
<\frac{1}{5\sqrt{2}}.
\]
Given $n\in\Na_{\geq N(B_0,\delta)}$, we have
\begin{align*}
g_{n}\left(s_{n,B_0} + \frac{\veps}{2};\scD, B_0-\delta\right) 
&= 
\sum_{\bfa\in\sfR(n)} \left( \frac{5\sqrt{2}}{(B_0 - \delta)^n \|q_n(\bfa)\|^2}\right)^{s_{n,B_0} + \frac{\veps}{2}} \\
&\leq 
\frac{(5\sqrt{2})^{\frac{\veps}{2}}}{(B_0 - \delta)^{\frac{n\veps}{2} }} 
\left(\frac{B_0}{B_0 - \delta}\right)^{ns_{n,B_0}} 
\sum_{\bfa\in\sfR(n)} \left( \frac{5\sqrt{2}}{B_0^n \|q_n(\bfa)\|^2}\right)^{s_{n,B_0}} \\
&=\frac{(5\sqrt{2})^{\frac{\veps}{2}}}{(B_0 - \delta)^{\frac{n\veps}{2} }} 
\left(\frac{B_0}{B_0 - \delta}\right)^{ns_{n,B_0}}
g_n(s_{n,B}; \scD , B_0) \\
&\leq 5\sqrt{2}  \left(\frac{B_0}{(B_0 - \delta)^{1+\frac{\veps}{4}}}\right)^{2n} \\
&< 1,
\end{align*}
and, thus, $s_{n,B_0} + \frac{\veps}{2} \geq s_{n, B_0 - \delta} \geq s_{n,B_0}$. Lastly, we obtain $s_{n,B_0} - \frac{\veps}{2} \leq s_{n, B_0 + \delta} \leq s_{n,B_0}$ from
\begin{align*}
g_n\left(s_{n,B_0 + \delta} + \frac{\veps}{2}; \scD, B_0\right) 
&= 
\sum_{\bfa\in\sfR(n)} \left( \frac{5\sqrt{2}}{B_0^n \|q_n(\bfa)\|^2}\right)^{s_{n,B_0 + \delta} + \frac{\veps}{2} } \\
&\leq 
\frac{(5\sqrt{2})^{ \frac{\veps}{2} }}{B_0^{ \frac{n\veps}{2} }} \left(\frac{B_0 + \delta}{B_0} \right)^{ns_{n,B_0+\delta}}  
\sum_{\bfa\in\sfR(n)} \left( \frac{5\sqrt{2}}{(B_0 + \delta)^n \|q_n(\bfa)\|^2}\right)^{s_{n,B_0 + \delta}} \\
&< 5\sqrt{2}\left(\frac{B_0 + \delta}{B_0^{1+ \frac{\veps}{4}}} \right)^{2n}  \qquad\qquad (\text{by }0<\veps<1 \text{ and } s_{n,B_0}<2)\\
&< 1.  
\end{align*}
The result is now a consequence of Proposition \ref{Le:LemPrel-01}.
\item Consider $0<\veps<2$ and let $B$ be such that $1<B^2<\psi^{\veps}$. If $n\in\Na$ satisfies
\[
n \geq
\max\left\{ N_2, \frac{2\veps\log(\psi) + \log\pi }{\veps\log\psi - 2 \log B} \right\},
\]
then
\begin{align*}
g_n(2-\veps;\scD,B) &= \sum_{\bfa\in\sfR(n)} \left( \frac{5\sqrt{2}}{B^n\|q_n(\bfa)\|^2} \right)^{2-\veps}\\
&=(5\sqrt{2})^{ 2-\veps} \sum_{\bfa\in\sfR(n)} \left( \frac{1}{B^n\|q_n(\bfa)\|^2} \right)^{2-\veps} \\
&\geq \sum_{\bfa\in\sfR(n)} \frac{\|q_n(\bfa)\|^{2\veps}}{B^{2n}\|q_n(\bfa)\|^{4}} \\
&\geq \frac{1}{(2\psi)^{\veps}}\left( \frac{\psi^{\veps}}{B^2}\right)^n  \sum_{\bfa\in\sfR(n)} \frac{1}{\|q_n(\bfa)\|^{4}} \\
&\geq \frac{1}{(2\psi)^{\veps}}\left( \frac{\psi^{\veps}}{B^2}\right)^n  \sum_{\bfa\in\sfR(n)} \frac{1}{|q_n(\bfa)|^{4}} \\
&\geq \frac{1}{(2\psi)^{\veps}}\left( \frac{\psi^{\veps}}{B^2}\right)^n  \sum_{\bfa\in\sfR(n)} \frac{1}{\pi} \leb(\clC_n(\bfa)) \\
&\geq \frac{1}{(2\psi)^{\veps}}\left( \frac{\psi^{\veps}}{B^2}\right)^n  \frac{1}{\pi}  \\
&>1.
\end{align*}
Hence, $2-\veps <s_{n,B} <2$ (by Proposition \ref{Le:LemPrel-02}) and we conclude that $s_{n,B}\to 2$ as $B\to 1$.
\item Take any $n\in\Na_{\geq N_1}$. First, note that $s_{n,B}\geq 1$, since  
\begin{align*}
g_n(1,\scD,B) &= \sum_{\bfa\in\sfR(n)} \frac{5\sqrt{2}}{B^n\|q_n(\bfa)\|^2} \\ 
&\geq \frac{5\sqrt{2}}{B^n} \sum_{\bfa\in\sfR(n)} \prod_{j=1}^n \frac{1}{(|a_j|+1)^2} \\ 
&\geq \frac{5\sqrt{2}}{B^n} \sum_{\bfa\in\scE^n} \prod_{j=1}^n \frac{1}{(|a_j|+1)^2} \\ 
&\geq \frac{5\sqrt{2}}{B^n} \left( \sum_{ a \in\scE} \frac{1}{(|a|+1)^2}\right)^n \\
&\geq \frac{(5\sqrt{2})}{2^nB^n} \left( \sum_{ a \in\scE} \frac{1}{|a|^2}\right)^n \\
&= \infty.
\end{align*}
Take $\veps\in (0,1)$ and choose $N\in\Na_{\geq N_1}$ so large that $n\in\Na_{\geq N}$ implies
\[
200 \left( \sum_{|a|\geq \sqrt{2}} \frac{1}{(|a_j|-1)^{2+2\veps}} \right)^n 
<
\left( \sum_{|a|\geq \sqrt{2}} \frac{1}{(|a_j|-1)^{2+\veps}} \right)^n.
\]
For such $n$, we have
\begin{align*}
\sum_{\bfa\in \sfR(n)} \frac{(5\sqrt{2})^{1+\veps}}{\|q_n(\bfa)\|^{2+2\veps}} 
&\leq 
(5\sqrt{2})^{1+\veps} (\sqrt{2})^{2+2\veps} \sum_{\bfa\in \sfR(n)} \prod_{j=1}^n \frac{1}{(|a_j|-1)^{2+2\veps}} \\
&<200 \left( \sum_{ |a|\geq \sqrt{2}} \frac{1}{(|a|-1)^{2+2\veps}}\right)^n \\
&<\left( \sum_{ |a|\geq \sqrt{2}} \frac{1}{(|a|-1)^{2+\veps}}\right)^n.
\end{align*}
Therefore, if
\[
B_{\veps} 
\colon= 
\left( \sum_{ |a|\geq \sqrt{2}} \frac{1}{(|a|-1)^{2+\veps}}\right)^{\frac{1}{1+\veps}}, 
\]
we have $g_n(1+\veps;\scD,B_{\veps} ) < 1$ and $s_{n,B_{\veps}}\leq 1 + \veps$. Taking $n\to\infty$, we conclude $s_{B_{\veps}}\leq 1+\veps$. Finally, since $B\mapsto s_{n,B}$ is decreasing for every $n\in\Na$, the function $B\to s_B$ is also decreasing, so $s_{B}\to 1$ when $B\to \infty$. 
\end{enumerate}
\end{proof}
Next, we show that $s_B$, $B>1$, is precisely the number defined in \eqref{Eq:Defsb}. In other words, in the definition of $g_n$, we can replace $5\sqrt{2}$ with $1$ and $\|\,\cdot\,\|$ with $|\,\cdot\,|$. For each $n\in\Na$, let $G_n:\RE_{>0}\to\RE$ be given by
\[
\forall \rho\in \RE_{>0}
\qquad
G_n(\rho;B)
\colon=
\sum_{\bfa\in\sfR(n)}
\frac{1}{\left(B^n |q_n(\bfa) |^2\right)^{\rho}}.
\]
Define
\[
S_{n,B}
\colon=
\inf\left\{ \rho>0 : G_n(\rho;B)\leq  1 \right\}.
\]
\begin{coro01}
We have
\[
\lim_{n\to\infty} S_{n,B} = s_B.
\]
\end{coro01}
\begin{proof}
Note that, for all $n\in\Na$ and $\bfa\in\sfR(n)$, we have 
\[
\frac{1}{B^n |q_n(\bfa)|^2}
\leq
\frac{5\sqrt{2}}{B^n \|q_n(\bfa)\|^2},
\]
so $G_n(\rho;B)\leq g_n(\rho;B)$ for all $\rho>0$, thus, $S_{n,B}\leq s_{n,B}$. Reversely, for any $\delta>0$, every large $n\in \Na$ and all $\bfa\in\sfR(n)$ satisfy
\[
\frac{5\sqrt{2}}{(B+\delta)^n\|q_n(\bfa)\|^2}
\leq 
\frac{1}{B^n|q_n(\bfa)|^2}.
\]
Therefore, $g_n(\rho;B+\delta)\leq G_n(\rho;B)$, so $s_{n,B+\delta}\leq S_{n,B}$. By Lemma \ref{Le:LemPrel-01}, we obtain
\[
s_{B+\delta}
\leq 
\liminf_{n\to\infty} S_{n,B}
\leq
\limsup_{n\to\infty} S_{n,B}
\leq 
s_B.
\]
To finalize the proof, let $\delta\to 0$ and apply Lemma \ref{Le:LemPrel-04}.\ref{Le:LemPrel-04-01}.
\end{proof}
\subsection{Upper bound of Hausdorff dimension}
In this subsection, $B>1$ will remain fixed.
For $n\in\Na$ and $\bfa=(a_1,\ldots,a_n)\in\sfR(n)$, define
\[
J_n(\bfa) \colon= \Cl\left( \{ z\in \mfF\colon a_k(z)=a_k \text{ for all } k\in\{1,\ldots,n\} \text{ and } \|a_{n+1}(z)\|\geq B^{n+1} \}\right).
\]
\begin{lem01}\label{Lem:TEO:UPPERBD:01}
For each $n\in\Na$ and $\bfa\in \sfR(n)$,
\[
|J_n(\bfa)|\leq \frac{2\lambda_1(\lambda_1+1)}{B^{n+1}\|q_n(\bfa)\|^2},
\]
where
\[
\lambda_1 \colon= \left( 1- \frac{1}{\sqrt{2}}\right)^{-1}.
\]
\end{lem01}
\begin{proof}
The proof is similar to that of Lemma \ref{LEM:GC:02:Case01:02}. 
\end{proof}
Recall that 
\[
F(B)\colon= \{z=[a_1,a_2,\ldots]\in\mathfrak{F}: \|a_n\|\geq B^n \text{ for infinitely many } n\in\mathbb{N}\}.
\]
\begin{propo01}\label{Propo-29}
We have
\[
F(B)\subseteq \bigcap_{m\in\Na}\bigcup_{n\geq m} \bigcup_{\bfb\in\sfR(n)} J_n(\bfb).
\]
\end{propo01}
\begin{proof}
Take $z=[a_1,a_2,\ldots]\in F(B)$. The result is obvious when $\bfa=\sanu\in \sfR$, so we assume that $\bfa\in \sfIr$. Theorem \ref{TEO:REG:IRREG:CLOS} ensures that for each $n\in\Na$ there is some $\bfb=(b_{n,1},\ldots, b_{n,n})\in\sfR(n)$ such that $z\in \oclC_n(\bfb)$ and $\|b_{n,j}\|=\|a_j\|$ for each $j\in \{1,\ldots,n\}$. Hence, if $m\in\Na$ is such that $\|a_m\|\geq B^m$, then $\|b_{n,m}\|\geq B^m$ for all $n\in\Na_{\geq m}$. 
\end{proof}

\begin{proof}[Proof of Theorem \ref{TEO:DIMHB}: Upper Bound]
Let $\veps>0$ be arbitrary. By Lemma \ref{Le:LemPrel-01}, we may pick $N\in\Na_{\geq N_1}$ large enough so that any $n\in\Na_{\geq N}$ verifies $s_{n,B}< s_B + \frac{\veps}{2}$, which gives $s_{n,B} + \tfrac{\veps}{2}< s_B +  \veps$ and thus
\[
g_n(s_B+\veps;\scD,B)<1.
\]
Put $\lambda_2\colon=2\lambda_1(\lambda_1+1)$. The $(s_b+2\veps)$-Hausdorff measure of $F(B)$, denoted by $\clH^{s_B+2\veps}(F(B) )$, is bounded above as follows:
\begin{align*}
\clH^{s_B+2\veps}&\left(F(B)\right) \leq \liminf_{N\to\infty} \sum_{n\geq N} \sum_{\bfa\in\sfR(n)} |J_n(\bfa)|^{s_B + 2\veps} \\
&\leq \lambda_2^{s_B+2\veps} \liminf_{N\to\infty} \sum_{n\geq N} \sum_{\bfa\in\sfR(n)} \left(\frac{1}{B^{n+1} \|q_n(\bfa)\|^2} \right)^{s_B + 2\veps} \\
&= \frac{\lambda_2^{s_B+2\veps}}{(5\sqrt{2}B)^{s_B+\veps}}  \liminf_{N\to\infty} \sum_{n\geq N} \sum_{\bfa\in\sfR(n)} \left(\frac{5\sqrt{2} }{B^{n}\|q_n(\bfa)\|^2} \right)^{s_B+\veps} \frac{1}{(B^{n+1} \|q_n(\bfa)\|^2)^{\veps}} &&\text{(by Lemma \ref{Lem:TEO:UPPERBD:01})} \\
&\leq \frac{\lambda_2^{s_B+2\veps}}{(5\sqrt{2}B)^{s_B+\veps}}  \liminf_{N\to\infty} \sum_{n\geq N} \left( \frac{2\psi}{B^{n+1}\psi^n}\right)^{\veps} \sum_{\bfa\in\sfR(n)} \left(\frac{5\sqrt{2} }{B^{n}\|q_n(\bfa)\|^2} \right)^{s_B+\veps}  \\
&= \frac{2^{\veps}\psi^{\veps}\lambda_2^{s_B+2\veps}}{B^{\veps}(5\sqrt{2}B)^{s_B+\veps}}  \liminf_{N\to\infty} \sum_{n\geq N}  \frac{1}{(B^{\veps}\psi^{\veps})^{n}} g_n(s_B+\veps;\scD,B)\\
&\leq  \frac{2^{\veps}\psi^{\veps} \lambda_2^{s_B+2\veps}}{B^{\veps}(5\sqrt{2}B)^{s_B+\veps}}  \liminf_{N\to\infty} \sum_{n\geq N}  \frac{1}{(B^{\veps}\psi^{\veps})^{n}} \\
&=0. 
\end{align*}
Therefore, $\dimh(F(B))\leq s_B$.
\end{proof}
\subsection{Lower bound of Hausdorff dimension}
In this subsection, $B>1$ will remain fixed.
To prove the lower bound,   we construct a family of Cantor subsets $\{F_M(B): M\in\Na_{\geq 4}\}$ such that $\dimh(F_M(B))\leq \dimh(F(B))$ for all $M\in\Na$.
Later, we estimate their Hausdorff dimension using the Mass Distribution Principle (see, for example, \cite[Theorem 4.2]{Falconer2014} for a proof).
\begin{lem01}[Mass Distribution Principle]\label{LE:MDP}
Let $F\subseteq \Cx$ be non-empty and let $\mu$ be a finite measure such that $\mu(F)>0$. If there are constants $c>0$, $r_0>0$ and $t\geq 0$ such that for all $z\in F$ and all $r\in (0,r_0)$ we have $\mu\left( \Dx(z;r)\right) \leq c r^t$, then $\dimh(F)\geq t$.
\end{lem01}

\subsubsection{Construction of Cantor subsets}
 Take $M\in \Na_{\geq 4}$. Let $(n_k)_{k\geq 1}$ be a sequence of natural numbers such that 
\begin{align}
B^{n_1-1} &> 11, \label{Eq:CondOnn1} \\
n_1 + \ldots + n_k &\leq \frac{n_{k+1}}{k+1}\,\text{ for all }k\in\Na. \label{Eq:CondOnnj}
\end{align}

For each $n\in\Na$, let $D_n$ be the set of words $\bfa=(a_1,\ldots,a_n)\in\sfR(n)$ such that every $j\in\{1,\ldots, n\}$ satisfies:
\begin{enumerate}[\rm i.]
\item If $j=n_k$ for some $k\in\Na$, then $\lfloor B^{n_k} \rfloor +2 \leq \|a_{n_k}\| \leq 2 \lfloor B^{n_k} \rfloor + 1$;
\item If $j\in\{1,\ldots,n\}\setminus \{n_k:k\in\Na\}$, then $\|a_j\|\leq M$.
\end{enumerate}
We put $D_0\colon=\vac$ and $D\colon= \bigcup_{n\in\Na_0} D_n$. For $n\in\Na$ and $\bfa\in \sfR(n)$, write
\[
A_n(\bfa)\colon= \bigcup_{\substack{b\in\scD \\ \bfa b\in D_{n+1}}}   \oclC_{n+1}(\bfa b).
\]
We refer to the sets $A_n(\bfa)$ as \textbf{fundamental sets} of level $n$. Define
\[
F_M(B) \colon= \bigcap_{n\in\Na} \bigcup_{\bfa\in D_n} A_n(\bfa).
\]
Observe that each union in the definition of $F_M(B)$ is finite. 

\begin{lem01}\label{LE:8.1}
$ \dimh(F_M(B)) \leq \dimh(F(B))$.
\end{lem01}
\begin{proof}
The countable stability of the Hausdorff dimension and its invariance under bi-Lipschtiz maps give $\dimh(F_M(B))= \dimh(F_M(B)\cap \clC_1(a))$ for all $a\in \scD$. Fix any $a$ such that $\oclC_1(a)\subseteq \mfFc$; for instance, $a=3+3i$. Let us show that 
\[
F_M(B)\cap \oclC_1(a)
\subseteq 
F(B).
\]
 Take any $\xi=[a_1,a_2,\ldots]\in F_M(B)$. Let $(\bfc_n)_{n\geq 1}$ be a sequence in $\scD$ such that $\bfc_n\in D_n$ and $\xi\in \oclC_n(\bfc_n)$ for all $n\in\Na$. Pick any $k\in\Na$. If $\xi\in\clC_{n_k}(\bfc_{n_k})$, then $\|a_{n_k}\|\geq \lfloor B^{n_k}\rfloor +2 > B^{n_k}$. Assume that $\xi\in\oclC_{n_k}(\bfc_{n_k})\setminus \clC_{n_k}(\bfc_{n_k})$ and write $\bfc_{n_k}=(c_{1,n_k},\ldots,c_{n_k,n_k})$. If $\xi$ is regular, then $\xi\in\oclC_{n_k}(\bfc_{n_k})\cap \clC_n(a_1,\ldots, a_{n_k})$ and, by Corollary \ref{CORO:NeighDigh},
\[
\|a_{n_k}\|
\geq \|c^{n_k}_{n_k}\|-1 
\geq \lfloor B^{n_k}\rfloor +1 
>B^{n_k}. 
\]
If $\xi$ is irregular, there is some $\bfb=(b_1,\ldots,b_{n_k})\in\sfR(n)$ such that $\xi\in \oclC_{n_k}(\bfb)$ and $\|b_j\|=\|a_j\|$ for all $j\in\{1,\ldots, n_k\}$. Therefore, $\xi\in \oclC_{n_k}(\bfc_{n_k})\cap \oclC_{n_k}(\bfb)$ and 
\[
\|a_{n_k}\|
= \|b_{n_k}\|
\geq \|c^{n_k}_{n_k}\|-1 
\geq \lfloor B^{n_k}\rfloor +1 
>B^{n_k}. 
\]
Since $k$ was arbitrary, we conclude that $\xi\in F(B)$.
\end{proof}

\subsubsection{Probability measure on $F_M(B)$}
Let $(m_j)_{j\geq 1}$ be the sequence given by 
\[
m_1\colon=n_1-1
\quad\text{ and }\quad
m_j\colon= n_j - n_{j-1}-1 \text{ for all } j\in \Na_{\geq 2}.
\]

First, take $n\in\{1,\ldots, n_1\}$.
\begin{enumerate}[\rm i.]
\item If $n=n_1-1$, for each $\bfa\in D_{n_1-1}$, put
\[
\mu\left( A_{n_1-1}(\bfa)\right) \colon= 
\left( \frac{5\sqrt{2}}{B^{m_1}\|q_{m_1}(\bfa)\|^2}\right)^{s_{m_1,B}(M)}.
\]
\item If $n=n_1$, for each $\bfa=(a_1,\ldots,a_{n_1})\in D_{n_1}$, put
\[
\mu\left( A_{n_1}(\bfa)\right) \colon= 
\frac{\mu\left( A_{n_1-1}(a_1,\ldots, a_{n_1-1})\right)}{\#\{d\in \scD: (a_1,\ldots, a_{n_1-1},d)\in D_{n_1}\}}.
\]\item If $n\in \{1,\ldots, n_1-2\}$, for each $\bfa\in D_{n}$, put
\[
\mu\left( A_{n}(\bfa)\right) \colon= 
\sum_{\substack{\bfb\in\sfR(n_1-n-1)\\ \bfa \bfb\in D_{n_1 - 1}}} \mu\left( A_{n_1-1}(\bfa\bfb)\right).
\]
\end{enumerate}
We now consider $n\in \{n_1+1, \ldots, n_2-1\}$.
\begin{enumerate}[\rm i.]
    \item If $n=n_2-1$, for each $\bfa=(a_1,\ldots, a_{n_2-1})\in D_{n_2-1}$, put
\[
\mu\left( A_{n_{2}-1}(\bfa)\right)
\colon=
\mu(A_{n_1}(a_1,\ldots, a_{n_1})) \left( \frac{5\sqrt{2}}{B^{m_{2}}\|q_{m_{2}}(a_{n_1+1}, \ldots, a_{n_{2}-1})\|^2} \right)^{s_{m_{2},B}(M)}.
\]
    \item If $n\in\{n_1+1, \ldots, n_2-2\}$, for each $\bfa=(a_1,\ldots, a_{n})\in D_{n}$, put
\[
\mu\left( A_n(\bfa)\right)
\colon=
\sum_{\substack{\bfb\in \sfR(n_{2} - 1 - n)\\ \bfa\bfb\in D_{n_{2}-1}}} \mu\left( A_{n_{2}-1}(\bfa\bfb)\right).
\]
\end{enumerate}
Assume that, for some $k\in\Na_{\geq 2}$, we have already defined $\mu(A_j(\bfa))$ for all $j\in \{1, \ldots, n_{k}-1\}$ and $\bfa \in D_j$. Take $n\in \{n_{k}, \ldots, n_{k+1}-1\}$ and $\bfa=(a_1,\ldots, a_n)\in D_n$.
\begin{enumerate}[\rm i.]
\item If $n=n_{k}$, put
\[
\mu\left(A_{n_k}(\bfa)\right)
\colon=
\frac{\mu\left(A_{n_k-1}(a_1,\ldots,a_{n_k-1}) \right)}{\#\{d\in\scD:\bfa d\in D_{n_k}\}}.
\]
\item If $n=n_{k+1} - 1$, put
\[
\mu\left( A_{n_{k+1}-1}(\bfa)\right)
\colon=
\mu(A_{n_k}(a_1,\ldots, a_{n_k})) \left( \frac{5\sqrt{2}}{B^{m_{k+1}}\|q_{m_{k+1}}(a_{n_k+1}, \ldots, a_{n_{k+1}-1})\|^2} \right)^{s_{m_{k+1},B}(M)}.
\]
\item If $n\in \{n_{k}+1, \ldots, n_{k+1} -2 \}$, put
\[
\mu\left( A_n(\bfa)\right)
\colon=
\sum_{\substack{\bfb\in \sfR(n_{k+1} - 1 - n)\\ \bfa\bfb\in D_{n_{k+1}-1}}} \mu\left( A_{n_{k+1}-1}(\bfa\bfb)\right).
\]

\end{enumerate}
Finally, for each $n\in\Na$, for $\bfa,\bfb\in D_n$ with $\bfa\neq\bfb$ we define $\mu(A_n(\bfa)\cap A_n(\bfb))=0$. 

Let us explain why $\mu$ is well defined. First, for all $m\in\Na$, the infimum defining $s_{m,B}(M)$ is actually a minimum, because the involved sum is finite. As a consequence, we have
\[
\mu\left( \bigcup_{\bfa\in D_{n_{1}-1}} A_{n_{1}-1}(\bfa)\right)
=
\sum_{\bfa\in D_{n_{1}-1}} \mu\left(A_{n_1-1}(\bfa)\right)
=1.
\]
Second, given $\bfa\in D_{n_1-1}$, we defined $\mu$ on the sets $\{A_{n_1}(\bfa d):\bfa d\in D_{n_1}\}$ by uniformly distributing the mass $\mu(A_{n_1-1}(\bfa))$ among them, so
\[
\mu(A_{n_1-1}(\bfa))
=
\sum_{\substack{d\in\scD\\ \bfa d \in D_{n_1}}} \mu(A_{n_1}(\bfa d)).
\]
The definition of $\mu$ ensures that the previous equation holds if we replace $n_1$ with any $n\in\{1,\ldots, n_1-1\}$. Moreover, it still holds for any $n\in\Na$. Therefore, by the Daniell-Kolmogorov consistency theorem (see, for example,  \cite[Corollary 8.22]{Kal2021}), the function $\mu$ extends to a unique Borel probability measure supported on $F_M(B)$.

\subsubsection{Bounding $\mu(A_n(\bfa))$}
In this section, we give an upper estimate of the $\mu$ measure of the sets $A_n(\bfa)$ in terms of $|A_n(\bfa)|$ (Lemma \ref{LEM:BND_MSR_AN}). The proof relies on estimates of a geometric (Propositions \ref{Propo:BndCardDn}, \ref{Propo:CotaAn}) and arithmetic nature (Proposition \ref{Propo:UpperBndsOnqn}).

\begin{propo01}\label{Propo:BndCardDn}
If $k\in\Na$ and $\bfa\in D_{n_k-1}$, then
\[
2\lfloor B^{n_k}\rfloor^2
\leq
\#\{d\in\scD: \bfa d\in D_{n_k}\} 
\leq
14\lfloor B^{n_k}\rfloor^2.
\]
\end{propo01}
\begin{proof}
Pick any $k\in\Na$ and write $B_1\colon=\lfloor B^{n_k}\rfloor + 2$ and $B_2\colon= 2\lfloor B^{n_k}\rfloor +1$. For the upper bound, we have
\[
\#\{d\in\scD: \bfa d\in D_{n_k}\}
\leq
\#\{d\in\scD: B_1\leq \|d\|\leq B_2\}
\]
and, since $B_1\geq 3$,
\begin{align*}
(2B_2+1)^2 - (2B_1-1)^2 &= 4(B_2-B_1)(B_2+B_1+1) \nonumber\\
&= 4(\lfloor B^{n_k}\rfloor-1) (3\lfloor B^{n_k} \rfloor +3) \nonumber\\
& \leq 14\lfloor B^{n_k}\rfloor^2. \nonumber
\end{align*}
For the lower bound, note that $\#\{d\in\scD: \bfa d\in D_{n_k}\}$ attains its minimum when $\omfF_n(\bfa)=\omfF_1(1+i)\pmod\RotaG$. Then, 
\begin{align*}
(B_2+1)^2 - (B_1+1)^2 &= (B_2-B_1)(B_2+B_1+2) \nonumber\\
&= (\lfloor B^{n_k}\rfloor-1)(3\lfloor B^{n_k}\rfloor + 4) \nonumber\\
& \geq 2\lfloor B^{n_k}\rfloor^2. \nonumber
\end{align*} 
\end{proof}
Define $\kappa_5\colon=\frac{\kappa_1}{\sqrt{2}}$, with $\kappa_1$ as in Proposition \ref{Prop:SizeEstimates}.\ref{Prop:SizeEstimates:i}, and $\kappa_6\colon= \frac{20}{3}$. 
\begin{propo01}\label{Propo:CotaAn}
The following estimates hold for all $n\in\Na$:
\begin{enumerate}[\rm i.]
\item \label{Propo:CotaAn:i} If $n=n_k-1$ for some $k\in\Na$, then every $\bfa\in D_{n_k-1}$ satisfies
\[
\frac{\kappa_5 }{\|q_{n_k-1}(\bfa)\|^2 \lfloor B^{n_k}\rfloor} 
\leq 
|A_{n_k-1}(\bfa)| 
\leq
\frac{\kappa_6}{\|q_{n_k-1}(\bfa)\|^2 \lfloor B^{n_k}\rfloor}.
\]
\item \label{Propo:CotaAn:ii} If $n=n_k$ for some $k\in\Na$, then every $\bfa\in D_{n_k}$ satisfies
\[
\frac{\kappa_5 }{\|q_{n_k}(\bfa)\|^2} 
\leq 
|A_{n_k}(\bfa)| 
\leq
\frac{\kappa_6}{\|q_{n_k}(\bfa)\|^2}.
\]
\item \label{Propo:CotaAn:iii} If $n_{k-1}<n<n_k-1$ for some $k\in\Na$, then every $\bfa\in D_{n}$ satisfies
\[
\frac{\kappa_5 }{\|q_{n}(\bfa)\|^2}
\leq 
|A_{n}(\bfa)| 
\leq
\frac{\kappa_6}{\|q_{n}(\bfa)\|^2}.
\]
\end{enumerate}
\end{propo01}
\begin{proof}
The estimates are a direct consequence of Proposition \ref{Prop:SizeEstimates}. We point out that when $n\in\Na$ falls in cases \ref{Propo:CotaAn:ii} and \ref{Propo:CotaAn:iii}, we have $|\clC_n(\bfa)|=|A_n(\bfa)|$ for any $\bfa\in D_n$.
\end{proof}
The next proposition contains two important bounds. The first one can be obtained inductively. The second one is a consequence of Proposition \ref{Prop:SizeEstimates}.\ref{Prop:SizeEstimates:ii}. We omit the details.
\begin{propo01}\label{Propo:UpperBndsOnqn}
\begin{enumerate}[\rm i.]
\item \label{Propo:UpperBndsOnqn-01} For any $n\in\Na$, $M\in\Na$, and $\bfa\in \Omega(n)\cap \scD(M)$, we have $\|q_n(\bfa)\|\leq (2M+1)^n$.
\item \label{Propo:UpperBndsOnqn-02} Let $(t_i)_{j\geq 1}$ be a sequence of natural numbers and define $(T_j)_{j\geq 1}$ by $T_j\colon=t_1+ \ldots +t_j$, $j\in\Na$. For any $n\in\Na$ and $\bfa=(a_1,\ldots,a_{T_n}) \in \Omega(T_n)$, 
\[
6^{n-1} \|q_{T_n}(\bfa)\| \leq  \|q_{t_1}(a_1,\ldots,a_{t_1})\| \cdots \|q_{t_n} (a_{T_{n-1}+1}, \ldots, a_{T_n})\|.
\]
\end{enumerate}
\end{propo01}

\begin{lem01}\label{LEM:BND_MSR_AN}
For each $0<t<s_B(M)$, there are $k_0\in\Na$ and $c_0>0$ such that every $k\in \Na_{\geq k_0}$, $n\in\Na_{\geq n_k}$, and $\bfa\in D_n$ satisfy
\[
\mu(A_n(\bfa)) \leq c_0|A_n(\bfa)|^t.
\]
\end{lem01}
\begin{proof}
Pick any $0<t<s_{B}(M)$. Choose a number $\gamma$ such that $1< \gamma <\psi$. Define $\veps$ and $\tilde{\veps}$ by
\[
\veps \colon=  \frac{s_B(M)-t}{4}, 
\quad
\tilde{\veps}\colon= \frac{\log \gamma}{\log(B) + 2\log(2M+1)}\,\veps. 
\]
Let $N_3\in\Na_{\geq N_2}$ be such that
\begin{align*}
|s_{n,B}(M)-s_{l,B}(M)| < \tilde{\veps} &\;\text{ for all }n,l\in\Na_{\geq N_3},\nonumber\\
t+2\veps < s_{n,B}(M) &\;\text{ for all }n \in\Na_{\geq N_3}. \nonumber
\end{align*}
Without any reference, we will use that
\[
t<t+\veps<t + 2\veps<s_{n,B}(M)\leq 2 
\quad\text{ for all } 
n\in\Na_{\geq N_3}. 
\]
Define $\kappa_7\colon= 81000$ and let $k_0\in\Na$ be such that every $k\in\Na_{\geq k_0+1}$ has the following properties:
\begin{equation}\label{EQ:COND-ON-nk-01}
m_k\geq N_3,
\end{equation}
\begin{equation}\label{EQ:COND-ON-nk-02}
\frac{n_k}{m_k} \leq 1 + \frac{2\veps}{t+2\veps},
\end{equation}
\begin{equation}\label{EQ:COND-ON-nk-03}
\frac{\log \kappa_7}{2k\veps \log B} + \frac{\log(B^{2N_3}(2M+1)^{4N_3}) + \log(6^8\cdot 50) - \log\kappa_5^2 - \log \kappa(\hat{s})}{2k\veps \log B} 
\leq 
\frac{1}{k-1} \sum_{j=k_0+1}^{k} n_j.
\end{equation}
The inequality \eqref{EQ:COND-ON-nk-03} is actually attained by every large $k$ because of \eqref{Eq:CondOnnj}. For notational convenience, for any $k\in\Na$ and any $j\in \{1, \ldots, k\}$, we write $q_{m_j}$ rather than $q_{m_j}(a_{n_{j-1}+1}, \ldots, a_{n_j-1})$. Lastly, for
\[
c_0\colon=  (2M+1)^{4n_{k_0}}B^{4(n_1 + \cdots + n_{k_0})},
\]
Proposition \ref{Propo:UpperBndsOnqn}.\ref{Propo:UpperBndsOnqn-01} guarantees 
\begin{equation}\label{Eq:CondOn_c0}
\prod_{j=1}^{k_0} \frac{1}{\lfloor B^{n_j}\rfloor} \left( \frac{5\sqrt{2}}{B^{m_j}\|q_{m_j}\|^2}\right)^{s_{m_j,B}(M)}
\leq
1
\leq
c_0\prod_{j=1}^{k_0} \left( \frac{5\sqrt{2}}{B^{2n_j}\|q_{m_j}\|^2}\right)^{t+\veps}.
\end{equation}
Although all of these parameters depend on $t$, the order we used to define them is relevant. We have the following dependencies: $\veps=\veps(t)$, $\tilde{\veps}=\tilde{\veps}(\veps)$, $N_3=N_3(\veps, \tilde{\veps})$, $k_0=k_0(\veps, N_3)$, and $c_0=c_0(k_0)$.

Take $n\in\Na$ such that $n\geq n_{k_0+1}$ and $\bfa\in D_n$. We now estimate $\mu(A_n(\bfa))$.
\subsubsection{Case $n=n_k$}
For any $\bfa=(a_1,\ldots, a_{n_k})\in D_{n_k}$, the definition of $\mu$ and Proposition \ref{Propo:BndCardDn} yield
\begin{align}
\mu(A_{n_k}(\bfa)) 
&\leq 
\frac{1}{2^k}
\prod_{j=1}^k \frac{1}{\lfloor B^{n_j}\rfloor^2}
 \left( \frac{5\sqrt{2}}{B^{m_j} \|q_{m_j}\|^2} \right)^{s_{m_j,B}(M)} \nonumber\\
&= 
\frac{1}{2^k}
\prod_{j=1}^{k_0} \frac{1}{\lfloor B^{n_j}\rfloor^2} 
\left( \frac{5\sqrt{2}}{B^{m_j} \|q_{m_j}\|^2} \right)^{s_{m_j,B}(M)}\;
\prod_{j=k_0+1}^k \frac{1}{\lfloor B^{n_j}\rfloor^2} 
\left( \frac{5\sqrt{2}}{B^{m_j} \|q_{m_j}\|^2} \right)^{s_{m_j,B}(M)} \nonumber\\
&\leq \frac{1}{2^k} c_0 
\prod_{j=1}^{k_0} \left( \frac{5\sqrt{2}}{B^{2n_j}\|q_{m_j}\|^2} \right)^{t + \veps} \; 
\prod_{j=k_0+1}^k \frac{1}{\lfloor B^{n_j}\rfloor^2} 
\left( \frac{5\sqrt{2}}{B^{m_j} \|q_{m_j}\|^2} \right)^{t + 2\veps} \qquad \text{(by \eqref{EQ:COND-ON-nk-01} and \eqref{Eq:CondOn_c0})} \nonumber\\
&\leq 
\frac{50^{k_0}}{2^k} c_0 
\prod_{j=1}^{k_0}  \frac{ 1 }{(B^{2n_j}\|q_{m_j}\|^2)^{t + \veps}}  \; 
\prod_{j=k_0+1}^k \frac{1}{\lfloor B^{n_j}\rfloor^2} 
\left( \frac{5\sqrt{2}}{B^{m_j} \|q_{m_j}\|^2} \right)^{t +2\veps}. \label{Eq:n=n_k01}
\end{align}
Let us focus on the last product:
\begin{align}
\prod_{j=k_0+1}^k \frac{1}{\lfloor B^{n_j}\rfloor^2} \left( \frac{5\sqrt{2}}{B^{m_j} \|q_{m_j}\|^2} \right)^{t+2\veps}
&\leq  \prod_{j=k_0+1}^k \frac{2^2}{ B^{2n_j}} \left( \frac{5\sqrt{2}}{B^{m_j} \|q_{m_j}\|^2} \right)^{t+2\veps} \nonumber\\
&= \prod_{j=k_0+1}^k 4  \left( \frac{5\sqrt{2}}{\|q_{m_j}\|^2} \right)^{t+2\veps}\; \prod_{j=k_0+1}^k \frac{1}{B^{2n_j + m_j (t+2\veps)}} \nonumber\\
&\leq \prod_{j=k_0+1}^k 4  \left( \frac{5\sqrt{2}}{\|q_{m_j}\|^2} \right)^{t+2\veps} \;\prod_{j=k_0+1}^k \frac{1}{B^{(n_j + m_j)(t+2\veps)}}  \nonumber\\
&\leq 
\prod_{j=k_0+1}^k \frac{200}{\|q_{m_j}\|^{2(t+2\veps)}}  \;
\prod_{j=k_0+1}^k \frac{1}{B^{(n_j + m_j)(t+2\veps)}}  \nonumber\\
&\leq  
\prod_{j=k_0+1}^k \frac{200}{\|q_{m_j}\|^{2(t+2\veps)}} \; 
\prod_{j=k_0+1}^k \frac{1}{B^{2n_j( t + \veps)}} \qquad \text{(by \eqref{EQ:COND-ON-nk-02})} \nonumber\\
&\leq 
200^{k-k_0}
\prod_{j=k_0+1}^k \frac{1}{\|q_{m_j}\|^{2(t+\veps)}} \; 
\prod_{j=k_0+1}^k \frac{1}{B^{2n_j( t + \veps)}} \nonumber\\
&= 
200^{k-k_0}
\prod_{j=k_0+1}^k \frac{ 1 }{(B^{2n_j}\|q_{m_j}\|^{2})^{t+\veps}} \; 
\prod_{j=k_0+1}^k \frac{1}{B^{2n_j \veps}} \label{Eq:n=n_k02}.
\end{align}
We substitute \eqref{Eq:n=n_k02} in \eqref{Eq:n=n_k01} to obtain
\begin{align}
\mu(A_{n_k}(\bfa)) 
&\leq 
100^k \frac{c_0}{4^{k_0}}
\prod_{j=1}^{k_0}  \frac{ 1 }{(B^{2n_j}\|q_{m_j}\|^2)^{t + \veps}}
\prod_{j=k_0+1}^k \frac{ 1 }{(B^{2n_j}\|q_{m_j}\|^{2})^{t+\veps}}
\prod_{j=k_0+1}^k \frac{1}{B^{2n_j \veps}} \nonumber\\
&\leq 
100^k \frac{c_0}{4^{k_0}}
\prod_{j=1}^{k}  \frac{ 1 }{(B^{2n_j}\|q_{m_j}\|^2)^{t + \veps}}
\prod_{j=k_0+1}^k \frac{1}{B^{2n_j \veps}} \nonumber\\
&\leq 
100^k c_0
\prod_{j=1}^{k}  \frac{ 1 }{(\lfloor B^{n_j} \rfloor^2\|q_{m_j}\|^2)^{t + \veps}}
\left(\prod_{j=k_0+1}^k \frac{1}{B^{2n_j}}\right)^{\veps}.\label{Eqn:Aux01}
\end{align}
Since for every $j\in\Na$ we have $\|a_{n_j}\|\leq 2\lfloor B^{n_j}\rfloor + 1$, we get $\|a_{n_j}\|\leq 3\lfloor B^{n_j}\rfloor$ and
\begin{align*}
\prod_{j=1}^k  \frac{1}{(\lfloor B^{n_j}\rfloor^2\|q_{m_j}\|^2)^{t+\veps}} 
&\leq 
\prod_{j=1}^k \left( \frac{9}{\|a_{n_j}\|^2 \|q_{m_j}\|^2} \right)^{t+\veps} \nonumber\\
&\leq 
9^{k(t+ \veps)} \left( \frac{1}{6^{2k-1} \|q_{n_k}(\bfa)\|}\right)^{2(t+\veps)} && \text{(by Proposition \ref{Propo:UpperBndsOnqn}.\ref{Propo:UpperBndsOnqn-02})} \\
&\leq 
\left( \frac{9^{t+\veps}}{6^{4t}}\right)^k 6^{2t} \frac{1}{\|q_{n_k}(\bfa)\|^{2t+2\veps}} \\
&\leq 
81^k6^4
\frac{1}{\|q_{n_k}(\bfa)\|^{2t+2\veps}}.
\end{align*}
Hence, plugging the previous inequality into \eqref{Eqn:Aux01}, we have
\begin{align*}
\mu(A_{n_k}(\bfa))
&\leq 
c_0 8100^k 6^4
\left(\prod_{j=k_0+1}^k \frac{1}{B^{2n_j}}\right)^{\veps}
\frac{1}{\|q_{n_k}(\bfa)\|^{2t+2\veps}} \nonumber\\
&\leq 
c_0 8100^k 6^4 
\left(\prod_{j=k_0+1}^k \frac{1}{B^{2n_j}}\right)^{\veps}
\frac{1}{\kappa_5^t} |A_{n_k}(\bfa)|^{t} &&\text{(by Proposition \ref{Propo:CotaAn})} \\ \nonumber
&\leq 
c_0 8100^k 6^4 \frac{1}{\kappa_5^t}
\left(\prod_{j=k_0+1}^k \frac{1}{B^{2n_j}}\right)^{\veps}
|A_{n_k}(\bfa)|^{t} \nonumber\\
&=
c_0 6^4 \frac{\kappa_7^k}{\kappa_5^2}
\left(\prod_{j=k_0+1}^k \frac{1}{B^{2n_j}}\right)^{\veps}
|A_{n_k}(\bfa)|^{t} \nonumber\\
&\leq c_0 |A_{n_k}(\bfa)|^{t}.
\end{align*}
The last inequality follows from \eqref{EQ:COND-ON-nk-03}. This finishes the proof for the case $n=n_k$. We state the next corollary for reference.
\begin{coro01}\label{Cor:AuxBnd01}
If $r\in\Na_{\geq k_0}$ and $\bfa\in D_{n_r}$, then
\begin{equation*}
\mu(A_{n_r}(\bfa)) 
\leq 
c_0 6^4 \kappa_7^r \left(\prod_{j=k_0+1}^r \frac{1}{B^{2n_j}}\right)^{\veps} \frac{1}{\|q_{n_r}(\bfa)\|^{2t+2\veps }}.
\end{equation*}
\end{coro01}
\subsubsection{Case $n=n_k-1$}
Take $k\in\Na_{k_0+1}$, $\bfa\in D_{n_k-1}$ and let $b\in\scD$ be such that $\bfa b\in D_{n_k}$. Then, using the definition of $\mu$ and Corollary \ref{Cor:AuxBnd01},
\begin{align*}
\mu(A_{n_k-1}(\bfa)) 
&=
\mu(A_{n_{k-1}}(a_1,\ldots,a_{n_{k-1}}))
\left( \frac{5\sqrt{2}}{\|q_{m_k}\|^2 B^{m_k}}\right)^{s_{m_k,B}(M)}\nonumber\\
&\leq c_0 \kappa_7^{k-1} 6^4
\left(\prod_{j=k_0+1}^{k-1} \frac{1}{B^{2n_j}}\right)^{\veps} 
\frac{1}{\|q_{n_{k-1}}\|^{2t+2\veps}} 
\left( \frac{5\sqrt{2}}{\|q_{m_k}\|^2 B^{m_k}}\right)^{ s_{m_k,B}(M) }\nonumber\\
&\leq c_0 \kappa_7^{k-1} 6^4 \left(\prod_{j=k_0+1}^{k-1} \frac{1}{B^{2n_j}}\right)^{\veps} \frac{1}{\|q_{n_{k-1}}\|^{2t+2\veps}} \left( \frac{5\sqrt{2}}{\|q_{m_k}\|^2 B^{m_k}}\right)^{ t+ 2\veps}\nonumber\\
&\leq c_0 \kappa_7^{k-1} 6^4
\left(\prod_{j=k_0+1}^{k-1} \frac{1}{B^{2n_j}}\right)^{\veps} 
\frac{50}{(\|q_{n_{k-1}}\|\|q_{m_k}\|)^{2(t+\veps)}} 
\frac{1}{B^{m_k(t+2\veps)}}\nonumber\\
&\leq c_0 \kappa_7^{k-1} 50\cdot 6^8
\left(\prod_{j=k_0+1}^{k-1} \frac{1}{B^{2n_j}}\right)^{\veps} 
\frac{1}{ \|q_{n_{k}-1}\| ^{2(t+\veps)}} 
\frac{1}{B^{tn_k }}  && \text{(by Proposition \ref{Propo:UpperBndsOnqn}.\ref{Propo:UpperBndsOnqn-02})}\nonumber\\
&\leq c_0 \kappa_7^{k-1} 50\cdot 6^8
\left(\prod_{j=k_0+1}^{k-1} \frac{1}{B^{2n_j}}\right)^{\veps} 
\frac{1}{ (B^{n_k }\|q_{n_{k}-1}\|^{2})^t} \nonumber\\
&\leq c_0 \kappa_7^{k-1} 50\cdot 6^8
\left(\prod_{j=k_0+1}^{k-1} \frac{1}{B^{2n_j}}\right)^{\veps} 
\frac{1}{ (\lfloor B^{n_k} \rfloor \|q_{n_{k}-1}\|^{2})^t} \nonumber\\
&\leq c_0 \kappa_7^{k-1} 50\cdot 6^8
\left(\prod_{j=k_0+1}^{k-1} \frac{1}{B^{2n_j}}\right)^{\veps} 
\frac{1}{\kappa_5^t} |A_{n_k-1}(\bfa)|^{t} &&\text{(by Proposition \ref{Propo:CotaAn})} \nonumber\\
&\leq c_0 
|A_{n_k-1}(\bfa)|^{t}. &&\text{(by \eqref{EQ:COND-ON-nk-03})} \nonumber
\end{align*}
\subsubsection{Case $n_{k-1}<n<n_k-1$}
Put
\[
l'\colon=n_k-n-1 \text{ and } l\colon=n-n_{k-1}. 
\]
Letting $\bfb$ run along the set $\{\bfc\in \sfR(l'): \bfa\bfc\in D_{n_k-1}\}$, applying Corollary \ref{Cor:AuxBnd01}, and writing $\bfa[\alpha,\beta]=(a_{\alpha}, a_{\alpha+1}, \ldots, a_{\beta} )$ for natural numbers $\alpha<\beta$, we obtain
\begin{align}
&\mu(A_n(\bfa)) 
=
\sum_{\bfb} \mu(A_{n_k-1}(\bfa\bfb))\nonumber\\
&=\sum_{\bfb} \mu(A_{n_{k-1}}(a_1,\ldots, a_{n_{k-1}})) 
\left( \frac{5\sqrt{2}}{B^{m_k}\|q_{m_k}(\bfa[n_{k-1}+1,n]\bfb)\|^2}\right)^{s_{m_k,B}(M)} \nonumber\\
&\leq 
c_0 \kappa_7^{k-1} 6^4
\left(\prod_{j=k_0 + 1}^{k-1} \frac{1}{B^{2n_j}}\right)^{\veps} 
\frac{1}{\|q_{n_{k-1}}(\bfa)\|^{2t+2\veps}} \sum_{\bfb} 
\left( \frac{5\sqrt{2}}{B^{m_k}\|q_{m_k}(\bfa[n_{k-1}+1,n]\bfb)\|^2}\right)^{s_{m_k,B}(M)} \nonumber\\
&\leq c_0 \kappa_7^{k-1} 6^4
\left(\prod_{j=k_0 + 1}^{k-1} \frac{1}{B^{2n_j}}\right)^{\veps} 
\frac{1}{\|q_{n_{k-1}}(\bfa)\|^{2t+2\veps}} \sum_{\bfb} \left( \frac{(5\sqrt{2})^2}{B^{m_k}\|q_{l}(\bfa[n_{k-1}+1,n])\|^2\|q_{l'}(\bfb)\|^2}\right)^{s_{m_k,B}(M)} \nonumber\\
&\leq c_0 \kappa_7^{k-1} 6^4
\left(\prod_{j=k_0 + 1}^{k-1} \frac{1}{B^{2n_j}}\right)^{\veps} \frac{(5\sqrt{2})^{s_{m_k,B(M)}}}{\|q_{n_{k-1}}(\bfa)\|^{2t+2\veps}\|q_{l}(\bfa[n_{k-1}+1,n])\|^{2s_{m_k,B}(M)}} \sum_{\bfb} \left( \frac{ 5\sqrt{2} }{B^{m_k}\|q_{l'}(\bfb)\|^2}\right)^{s_{m_k,B}(M)} \nonumber\\
&\leq 
c_0 50 \kappa_7^{k-1} 6^4
\left(\prod_{j=k_0 + 1}^{k-1} \frac{1}{B^{2n_j}}\right)^{\veps}
\left(\frac{ 1}{\|q_{n_{k-1}}(\bfa)\|\,\|q_{l}(\bfa[n_{k-1}+1,n])\|}\right)^{2t+ 2\veps} \sum_{\bfb} \left( \frac{ 5\sqrt{2} }{B^{m_k}\|q_{l'}(\bfb)\|^2}\right)^{s_{m_k,B}(M)} \nonumber\\
&\leq c_0 50 \kappa_7^{k-1}  6^4
\left(\prod_{j=k_0 + 1}^{k-1} \frac{1}{B^{2n_j}}\right)^{\veps}
\left(\frac{6}{\|q_{n}(\bfa)\|}\right)^{2t+ 2\veps} \sum_{\bfb} \left( \frac{ 5\sqrt{2} }{B^{m_k}\|q_{l'}(\bfb)\|^2}\right)^{s_{m_k,B}(M)} \nonumber\\
&\leq 
c_0 50 \kappa_7^{k-1} 6^8
\left(\prod_{j=k_0 + 1}^{k-1} \frac{1}{B^{2n_j}}\right)^{\veps}
\frac{1}{\|q_{n}(\bfa)\|^{2t+ 2\veps}}
\sum_{\bfb} \left( \frac{5\sqrt{2}}{B^{m_k}\|q_{l'}(\bfb)\|^2}\right)^{s_{m_k,B}(M)} \nonumber\\
&\leq 
c_0 50 \kappa_7^{k-1} 6^8
\left(\prod_{j=k_0 + 1}^{k-1} \frac{1}{B^{2n_j}}\right)^{\veps}
\frac{1}{\|q_{n}(\bfa)\|^{2t+ 2\veps}}
g_{l'}(s_{m_k,B}(M);\scD(M),B). \label{Eq:CotaMedCyl01}
\end{align}

Since $l+l'=m_k$, in view of \eqref{Eq:Le01Claim02:s:hat}, we have
\[ 
\kappa(\hat{s})g_{l'}(s_{m_k,B}(M);\scD(M),B)g_l(s_{m_k,B}(M);\scD(M),B)
\leq 
g_{m_k}(s_{m_j,B}(M);\scD(M),B)=1,
\]
so
\[
g_{l'}(s_{m_k,B}(M);\scD(M),B)
\leq \frac{1}{\kappa(\hat{s})g_{l}(s_{m_k,B}(M);\scD(M),B)}.
\]
We, thus, need a lower bound of $g_{l}(s_{m_j,B}(M);\scD(M),B)$. We obtain it by considering two cases: $l\leq N_3$ and $l> N_3$. If $l\leq N_3$, then
\begin{align*}
g_{l}(s_{m_k,B}(M);\scD(M),B) &= \sum_{\bfc\in \scD(M)^l\cap\sfR(l) } \left( \frac{5\sqrt{2}}{B^l\|q_l(\bfc)\|^2}\right)^{s_{m_k,B}(M)} \nonumber\\
&\geq \frac{1}{B^{N_3s_{m_k,B}(M)}} \sum_{\bfc\in \scD(M)^l\cap\sfR(l) } \frac{1}{\|q_l(\bfc)\|^{2s_{m_k,B}(M)}} \nonumber\\
&\geq \frac{1}{B^{2N_3} (2M+1)^{4N_3}},\nonumber
\end{align*}
so, by \eqref{Eq:CotaMedCyl01} and Proposition \ref{Propo:UpperBndsOnqn},
\begin{align*}
\mu(A_n(\bfa))
&\leq 
\frac{c_0 50 \kappa_7^{k-1} 6^8 }{\kappa(\hat{s})}
\left(\prod_{j=k_0 + 1}^{k-1} \frac{1}{B^{2n_j}}\right)^{\veps}
\frac{1}{\|q_{n}(\bfa)\|^{2t+ 2\veps}}
B^{2N_3} (2M+1)^{4N_3} \\
&\leq 
\frac{ c_0 50 \kappa_7^{k-1} 6^8}{\kappa(\hat{s}) }
\left(\prod_{j=k_0 + 1}^{k-1} \frac{1}{B^{2n_j}}\right)^{\veps}
B^{2N_3} (2M+1)^{4N_3}
\frac{|A_n(\bfa)|^{t+\veps}}{\kappa_5^{t+\veps}} \\
&\leq 
\frac{ c_0 50 \kappa_7^{k-1} 6^8 \kappa_5^2}{\kappa(\hat{s}) }
B^{2N_3} (2M+1)^{4N_3}
\left(\prod_{j=k_0 + 1}^{k-1} \frac{1}{B^{2n_j}}\right)^{\veps}
|A_n(\bfa)|^{t+\veps} \\
&\leq c_0 |A_n(\bfa)|^{t}.
\end{align*}

Now assume that $l>N_3$. The definitions of $N_3$ and $\tilde{\veps}>0$ yield
\begin{align*}
g_{l}(s_{m_j,B}(M);\scD(M),B) &= \sum_{\bfc\in \scD(M)^l\cap\sfR(l) } \left( \frac{5\sqrt{2}}{B^l\|q_l(\bfc)\|^2}\right)^{s_{m_k,B}(M)} \nonumber\\
&\geq \sum_{\bfc\in \scD(M)^l\cap\sfR(l) } \left( \frac{5\sqrt{2}}{B^l\|q_l(\bfc)\|^2}\right)^{s_{l,B}(M) + \tilde{\veps}} \nonumber\\
&\geq \sum_{\bfc\in \scD(M)^l\cap\sfR(l) } \frac{1}{(B^l(2M+1)^{2l})^{\tilde{\veps}}} \left( \frac{5\sqrt{2}}{B^l\|q_l(\bfc)\|^2}\right)^{s_{l,B}(M)} \nonumber\\
&= \frac{1}{(B^l(2M+1)^{2l})^{\tilde{\veps}}} g_{l}(s_{l,B}(M);\scD(M),B) \nonumber\\
&= \frac{1}{(B^l(2M+1)^{2l})^{\tilde{\veps}}} \geq \frac{1}{\gamma^{l\veps}} \geq \frac{1}{\gamma^{n\veps}}.
\end{align*}
The previous inequality and \eqref{Eq:CotaMedCyl01} imply that
\begin{align*}
\mu(A_n(\bfa)) 
&\leq 
\frac{c_0 50 \kappa_7^{k-1} 6^8}{\hat{\kappa}(s)} 
\left(\prod_{j=k_0 + 1}^{k-1} \frac{1}{B^{2n_j}}\right)^{\veps}
\frac{\gamma^{n\veps}}{\|q_{n}(\bfa)\|^{2t+ 2\veps}} \\
&\leq 
\frac{c_0 50 \kappa_7^{k-1} 6^8}{\hat{\kappa}(s)} 
\left(\prod_{j=k_0 + 1}^{k-1} \frac{1}{B^{2n_j}}\right)^{\veps}
\frac{\gamma^{n\veps}}{\|q_{n}(\bfa)\|^{2\veps}}
\frac{1}{\|q_{n}(\bfa)\|^{2t}} \\
&\leq 
\frac{c_0 50 \kappa_7^{k-1} 6^8\psi^{\veps}}{\kappa(\hat{s})\kappa_5^t} 
\left(\prod_{j=k_0 + 1}^{k-1} \frac{1}{B^{2n_j}}\right)^{\veps}
\frac{\gamma^{n\veps} }{\psi^{n\veps}}
|A_n(\bfa)|^t\\
&\leq 
c_0 |A_n(\bfa)|^t.
\end{align*}
\end{proof}
\subsubsection{Estimate of $\mu(\Dx(z,r))$}
\begin{figure}[h!]
\begin{center}
\includegraphics[scale=0.75,  trim={5cm 17cm 4cm 4cm},clip]{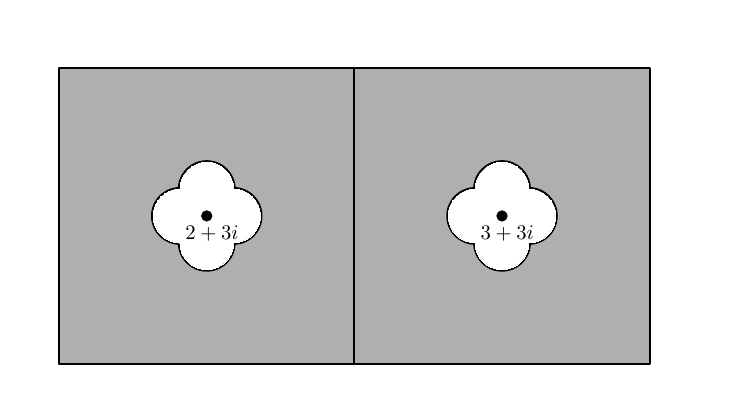}
\caption{ The sets $\iota[A_1(2+3i)]$ and $\iota[A_1(3+3i)]$. \label{Fig-BallsMeas}}
\end{center}
\end{figure}

In the context of regular continued fractions, the measure of the balls is bounded by estimating the length of the gaps between two different fundamental sets (see \cite[Section 3.2]{WanWu2008}). The strategy must be modified for Hcfs. Let us illustrate why, with a simple example. When $M=5$ and $n_1\geq 3$, both $a=3+3i$ and $2+3i$ belong to $D_1$ but $A_1(a)\cap A_1(b)\neq\vac$ (see Figure \ref{Fig-BallsMeas}). We may thus lose any hope of obtaining a non-trivial lower bound of the distance between any two fundamental sets of the same level. This is when our discussion on the geometry of Hcfs pays off.

Let $\tilde{\rho}$ and $\rho$ be as in Theorems \ref{TEO:GC:01} and \ref{TEO:GC:02}, respectively. Set $\rho' \colon= \min\{\tilde{\rho},\rho\}$ and let $\kappa_3$ be as in Corollary \ref{CORO:NeighCylProp}. For each $\bfa\in \scD^{\Na}$ such that $\prefi(\bfa;n)\in D_n$ for all $n\in\Na$, define
\[
\tilde{G}_{n}(\bfa) \colon= \frac{\rho'}{\kappa_3} |A_n(\bfa)|.
\]
Let $r_0$ be
\[
r_0  
\colon=
\min \left\{ \tilde{G}_n(\bfa): n\in \{1,\ldots, n_{k_0}\}, \, \bfa\in D_n\right\}.
\]
\begin{lem01}\label{LEM:dimFB:MeasBalls}
There is a constant $\nu=\nu(M,B,t)>0$ such that for all $z\in F_M(B)$ and every $r\in (0,r_0)$ we have $\mu(\Dx(z;r)) \leq \nu r^t$.
\end{lem01}
\begin{proof}
Take a real number $r$ with $0<r<r_0$, $z\in F_M(B)$, and $\bfa=\sanu \in \scD^{\Na}$ such that $z\in A_n(a_1, \ldots, j)\in D_j$ for all $j\in \Na$. Let $n\in\Na$ be such that 
\[
\tilde{G}_{n+1}(a_1,\ldots, a_n,a_{n+1})
\leq
r
<
\tilde{G}_{n}(a_1,\ldots, a_n).
\]
We examine three different cases for $n$.
\subsubsection{Case $n=n_k-1$} 
First, we further assume $\omfF_n(\bfa)=\omfF$, $\omfF_1(1+2i) \mod \RotaG$. Then, the disc $\Dx(z;r)$ intersects exactly one fundamental set of order $n$. We now consider two possibilities:
\[
r\leq  \kappa_3 |A_{n_k}(a_1,\ldots, a_{n_k})|
\;\text{ or }\;
\kappa_3 |A_{n_k}(a_1,\ldots, a_{n_k})|<r.
\]
In the first option, $\Dx(z;r)$ intersects at most six cylinders of level $n_k$, which are $\clC_n(\prefi(\bfa;n))$ and five of its neighbors, so it intersects at most six fundamental sets of level $n_k$, hence
\[
\mu(\Dx(z;r)) 
\leq 
6\kappa_5 c_0|A_{n_k}(\bfa)|^t
<
6\kappa_5 c_0\frac{\kappa_3}{\tilde{\rho}} r^t.
\]

Now assume that $\kappa_3 |A_{n_k}(a_1,\ldots, a_{n_k})|<r$. Define $\mathsf{N}\colon= \#\left\{\bfd\in D_{n_k}: \Dx(z;r)\cap A_{n_k}(\bfd)\neq\vac\right\}$. Let us show that there is an absolute constant $\nu_2$ such that
\begin{equation}\label{EQ:Caso_nk-1_b1:NBnd}
\mathsf{N}
\leq
6^4 2^{11} \pi \nu_2^2  \|q_{n_k}(a_1,\ldots, a_{n_k-1})\|^4 \lfloor B^{n_k}  \rfloor^4 r^2.
\end{equation}
First, note that there is an absolute constant $\nu_1$ such that for all $c,d\in \scD$ satisfying $(a_1,\ldots, a_{n_k-1}, c)\in D_{n_k}$ and $(a_1,\ldots, a_{n_k-1}, d)\in D_{n_k}$ we have 
\begin{equation}\label{EQ:Caso_nk-1_b2}
\frac{1}{\nu_1} 
\leq 
\frac{|A_{n_k}(a_1,\ldots, a_{n_k-1}, c)|}{|A_{n_k}(a_1,\ldots, a_{n_k-1}, d)|}
\leq 
\nu_1.
\end{equation}
Certainly, since $\frac{1}{3} \leq \frac{\| d\|}{\| c\|}\leq 3$, 
\[
\frac{|A_{n_k}(a_1,\ldots, a_{n_k-1}, c)|}{|A_{n_k}(a_1,\ldots, a_{n_k-1}, d)|} 
= 
\frac{|\clC_{n_k}(a_1,\ldots, a_{n_k-1}, c)|}{|\clC_{n_k}(a_1,\ldots, a_{n_k-1}, d)|} 
\asymp
\frac{|q_{n_k}(a_1,\ldots, a_{n_k-1}, d)|^2}{|q_{n_k}(a_1,\ldots, a_{n_k-1}, c)|^2}
\asymp 
\frac{\|d\|^2}{\|c\|^2} 
\asymp 
1.
\]
Then, for $\nu_2\colon= 1 + \frac{\nu_1}{\kappa_3}$, we have that $A_{n_k}(a_1,\ldots, a_{n_k-1},d)\cap \Dx(z;r)\neq \vac$ implies that $A_{n_k}(a_1,\ldots, a_{n_k-1},d)\subseteq \Dx(z;\nu_2 r)$. Therefore, in view of the estimate 
\begin{align*}
\leb\left( \clC_{n_k}(a_1,\ldots, a_{n_k-1},d)\right) &\geq 
\frac{\kappa_1 }{2^2\|q_{n_k}(a_1,\ldots, a_{n_k-1},d)\|^4} \\
&\geq 
\frac{\kappa_1 6^{-4} }{2^2\|q_{n_k-1}(a_1,\ldots, a_{n_k-1})\|^4 \|d\|^4} \\
&\geq 
\frac{\kappa_1 6^{-4} }{2^2\|q_{n_k-1}(a_1,\ldots, a_{n_k-1})\|^4 2^4\lfloor B^{n_k}  \rfloor^4 \left( 1+  \frac{1}{2\lfloor B^{n_k}\rfloor} \right)^4 } \\
&\geq  \frac{\kappa_1}{6^4 2^{11} }\frac{ 1 }{\|q_{n_k-1}(a_1,\ldots, a_{n_k-1})\|^4 \lfloor B^{n_k}  \rfloor^4 },
\end{align*}
we get \eqref{EQ:Caso_nk-1_b1:NBnd}. We, thus, conclude that
\begin{align}
\mu\left( \Dx(z;r)\right) &\leq \frac{\mu\left( A_{n_k-1}(a_1,\ldots, a_{n_k-1})\right)} {\#\{d\in\scD:\bfa d\in D_{n_k}\}} \mathsf{N} \nonumber\\
&\leq 
\frac{6^42^{11} \pi \nu_2^2}{2\lfloor B^{n_k}\rfloor^2}   \mu\left( A_{n_k-1}(a_1,\ldots, a_{n_k-1})\right) \|q_{n_k-1}(a_1,\ldots, a_{n_k-1})\|^4 \lfloor B^{n_k}  \rfloor^4 r^2 \nonumber\\
&= 
6^4 2^{10} \pi \nu_2^2   \mu\left( A_{n_k-1}(a_1,\ldots, a_{n_k-1})\right) \|q_{n_k-1}(a_1,\ldots, a_{n_k-1})\|^4 \lfloor B^{n_k}  \rfloor^2 r^2. \label{EQ:Cota-muDx-01}
\end{align}
Moreover, since $A_{n_k-1}(\bfa)$ is the single fundamental set of level $n_k-1$ intersecting $\Dx(z;r)$, we have that
\begin{equation}\label{EQ:Cota-muDx-02}
\mu\left( \Dx(z;r) \right)
=
\mu\left( \Dx(z;r)\cap A_{n_k-1}(a_1,\ldots, a_{n_k-1})\right)
\leq 
\mu\left(A_{n_k-1}(a_1,\ldots, a_{n_k-1})\right).
\end{equation}
Combine inequalities \eqref{EQ:Cota-muDx-01} and \eqref{EQ:Cota-muDx-02} and recall that $t<2$, to see that for some constants $\nu_3,\nu_3'>0$ we have
\begin{align*}
\mu\left( \Dx(z;r)\right) 
&\leq \mu\left( A_{n_k-1}(a_1,\ldots, a_{n_k-1})\right) \min\left\{ 1, 6^4 2^{10}\pi \nu_2^2 \|q_{n_k}(a_1,\ldots, a_{n_k-1})\|^4 \lfloor B^{n_k}  \rfloor^2 r^2\right\} \\
&\leq 
\mu\left( A_{n_k-1}(a_1,\ldots, a_{n_k-1})\right) 1^{1-\tfrac{t}{2}} \left(6^4 2^{10} \pi \nu_2^2 \|q_{n_k-1}(a_1,\ldots, a_{n_k-1})\|^4 \lfloor B^{n_k}  \rfloor^2 r^2\right)^{\frac{t}{2}} \\
& \leq 
\nu_3'\mu\left( A_{n_k-1}(a_1,\ldots, a_{n_k-1})\right) \|q_{n_k-1}(a_1,\ldots, a_{n_k-1})\|^{2t} \lfloor B^{n_k}  \rfloor^{t} r^{t} \\
&\leq 
c_0\nu_3 r^{t}.
\end{align*}
Now, assume that either $\omfF_n(\bfa)\equiv \omfF_1(2) \mod \RotaG$ or $ \omfF_n(\bfa) \equiv \omfF_1(1+i) \mod \RotaG$. By the Second Geometric Construction (Theorem \ref{TEO:GC:02}), the disc $\Dx(z;r)$ intersects at most two fundamental sets of level $n_k-1$ and they are subsets of neighboring $n_{k}-1$ cylinders, say $A_{n_k-1}(a_1,\ldots, a_{n_k-1})$ and $A_{n_k-1}(\bfb)$. Hence, 
\[
\mu\left( \Dx(z;r) \right)
= 
\mu\left( \Dx(z;r)\cap A_{n_k-1}(a_1,\ldots, a_{n_k-1}) \right) + \mu\left( \Dx(z;r)\cap A_{n_k-1}(\bfb) \right).
\]
Propping on the above argument we may then assert that, for some constant $\nu_4=\nu_4(M,B,t)>0$, we have
\begin{align*}
\mu\left( \Dx(z;r) \right)
&\leq 
\nu_3 c_0|A_{n_k-1}(a_1,\ldots, a_{k_n-1})|^t + \nu_3c_0|A_{n_k-1}(\bfb)|^t \\
&\leq c_0\nu_3 (1+\kappa_3) r^t.
\end{align*}
\subsubsection{Case $n=n_k-2$}

First, note that for any $\bfb\in D_{n_k-2}$ and $c,d\in \scD$ such that $\bfb c, \bfb d\in D_{n_k-1}$, we have
\begin{align*}
\frac{\mu\left( A_{n_k-1}(\bfb c)\right)}{\mu\left( A_{n_k-1}(\bfb d)\right)} 
&= \left( \frac{\|q_{m_k}(\bfb d)\|^2}{ \| q_{m_k}(\bfb c)\|^2} \right)^{s_{m_k,B}(M)} \\
&\leq \left(1800 \frac{\|d\|^2}{\|c\|^2}\right)^{s_{m_k,B}(M)}  \\
&\leq \left(1800 M^2 \right)^{ s_{m_k,B}(M)}   \\
&\leq 1800^2 M^4. && \text{ (Proposition \ref{Le:LemPrel-03})}
\end{align*}
As a consequence, using $\#\{c\in\scD:\bfb c\in D_{n_k-1}\}\leq \#\{d\in\scD:\|d\|\leq M\} \leq (2M+1)^2$, we have
\[
\mu\left( A_{n_k-2}(\bfb) \right)
\leq
8\kappa_5 1800^2 M^4 (2M+1)^2 c_0 \mu\left( A_{n_k-1}(\bfb c) \right).
\]
Finally, since $\Dx(z;r)$ intersects at most six fundamental sets of level $n_k-2$, one of them being $A_{n_k-2}(a_1, \ldots, a_{n_k-2})$ and the others are contained in the neighboring cylinders of $\clC_{n_k-2}(a_1, \ldots, a_{n_k-2})$, we have
\[
\mu\left( \Dx(z;r)\right)
\leq
48\kappa_5 1800^2 M^4 (2M+1)^2 c_0 r^t.
\]
\subsubsection{Case $n_{k-1} \leq n \leq n_{k}-3$}
Consider any $\bfb=(b_1,\ldots, b_n)\in D_{n}$ and any $c,d\in\scD$ such that $\bfb c$, $\bfc d\in D_{n+1}$. We show the existence of some constant $\nu_5>0$ for which 
\begin{equation}\label{Eq:TeoB:UpperBound:Case3}
\frac{\mu\left(A_{n+1}(\bfb c)\right)}{\mu\left(A_{n+1}(\bfb d)\right)}
\leq 
 4\cdot 3600^{2\nu_5} M^{2\nu_5}.    
\end{equation}
Once this is achieved, the argument in the case $n=n_k-2$ can be applied. 

First, observe that we may express $\mu(A_{n+1}(\bfb c))$ as follows (we sum over all the words $\bfc\in \sfR(n_{k}- n - 2)$ such that $\bfb \,c\, \bfc\in D_{n_{k}-1})$:
\begin{align*}
&\mu\left( A_{n+1}(\bfb c)\right) 
=  
\sum_{ \bfc} \mu\left( A_{n_{k}-1}(\bfb c\bfc )\right)\\
&=  
\sum_{\bfc} \mu\left( A_{n_{k-1}}( \prefi(\bfb;n_{k-1}) )\right) \left( \frac{5\sqrt{2}}{B^{m_{k}} \|q_{m_{k}}(b_{n_{k-1}},\ldots,b_n,c,\bfc)\|^2}\right)^{s_{n_{k}-n-2,B}(M)} \\
&= \frac{\mu\left( A_{n_{k-1}}(\prefi(\bfb;n_{k-1}) )\right) (5\sqrt{2})^{s_{n_k-n-2,B}(M)}}{B^{m_{k} s_{n_k-n-2,B}(M)}} \sum_{\bfc}   \frac{ 1 }{ \|q_{m_{k}}(b_{n_{k-1}},\ldots,b_n,c,\bfc)\|^{2s_{n_k-n-2,B}(M)}}.
\end{align*}
Next, we estimate the sum in the last expression. As an upper estimate, we have
\begin{align*}
\sum_{\bfc}   &\frac{ 1 }{ \|q_{m_{k}}(b_{n_{k-1}},\ldots,b_n,c,\bfc)\|^{2s_{n_k-n-2,B}(M)}} \\
&\leq \left(\frac{50}{\|q_{n-n_{k-1}+1}(b_{n_k},\ldots,b_{n})\| \|c\|} \right)^{2s_{n_k-n-2,B}(M)} \sum_{\bfc}\frac{1}{\|q_{n_{k}- n - 2}(\bfc)\|^{2s_{n_k-n-2,B}(M)}} \\
&\leq \left(\frac{50}{\kappa_5 \| q_{n-n_k+1}(b_{n_k},\ldots, b_{n})\| \|c\|} \right)^{2s_{n_k-n-2,B}(M)} \sum_{\bfc} |\clC_{n_{k}-n-2}(\bfc)|^{s_{n_k-n-2,B}(M)}.
\end{align*}
In the next lower bound, the sums run along those words $\bfd\in \sfR(n_{k}- n - 2)$ such that $\bfb \,d\, \bfd\in D_{n_{k}-1}$:
\begin{align*}
\sum_{\bfd}   &\frac{ 1 }{ \|q_{m_{k}}(b_{n_{k-1}},\ldots, b_n,d,\bfd)\|^{2s_{n_k-n-2,B}(M)}}  \\
&\geq 
\left(\frac{1}{36\|q_{n-n_{k-1}+1}(b_{n_k},\ldots, b_{n})\| \|d\|} \right)^{2s_{n_k-n-2,B}(M)} \sum_{\bfd}\frac{1}{\|q_{n_{k}- n - 2}(\bfd)\|^{2s_{n_k-n-2,B}(M)}} \\
&\geq \left(\frac{1}{72 \| q_{n-n_k+1}(b_{n_k}, \ldots b_{n})\| \|d\|} \right)^{2s_{n_k-n-2,B}(M)} \sum_{\bfd} |\clC_{n_{k}-n-2}(\bfd)|^{s_{n_k-n-2,B}(M)}.
\end{align*}
Hence, if $\nu_5\colon= \nu_5(M,B)\colon=\sup\{s_{n,B}(M):n\in\Na\}$ (see Proposition \ref{Le:LemPrel-03}), then
\begin{align*}
\frac{\mu\left(A_{n+1}(\bfb c)\right)}{\mu\left(A_{n+1}(\bfb d)\right)}
&\leq \left( 3600 \frac{\|d\|}{\|c\|}\right)^{2s_{n_k-n-2,B}(M)} \frac{\sum_{ \bfc } \left|\clC_{n_k-n-2}(\bfc)\right|^{s_{n_k-n-2,B}(M)} }{\sum_{ \bfd } \left|\clC_{n_k-n-2}(\bfd)\right|^{s_{n_k-n-2,B}(M)} } \\
&\leq  3600^{2\nu} M^{2\nu} \frac{\sum_{ \bfc } \left|\clC_{n_k-n-2}(\bfc)\right|^{s_{n_k-n-2,B}(M)} }{\sum_{ \bfd } \left|\clC_{n_k-n-2}(\bfd)\right|^{s_{n_k-n-2,B}(M)} } \\
&\leq  4\cdot 3600^{2\nu_5} M^{2\nu_5}.
\end{align*}
For the last inequality, we considered $\omfF_{n+1}(\bfb c)=\omfF$ in the numerator and $\omfF_{n+1}(\bfb c)=\omfF_1(1+i)$ in the denominator. This shows \eqref{Eq:TeoB:UpperBound:Case3}. 
\end{proof}
\begin{proof}[Proof of Theorem \ref{TEO:DIMHB}: Lower bound]
For any $M\in\Na_{\geq 4}$, the Mass Distribution Principle applied to $F_M(B)$ implies that $s_B(M)\leq \dimh(F_M(B))$. Therefore, by lemmas \ref{LE:8.1} and \ref{Le:LemPrel-03},
\[
s_B
=
\lim_{M\to\infty} s_B(M)
\leq
\dimh(F(B)).
\]
\end{proof}
The proof of Theorem \ref{TEO:DIMHB} is flexible enough to assume that the terms of the sequence $(n_j)_{j\geq 1}$ lie in a prescribed infinite subset of $\Na$. Let $\scL \subseteq \Na$ be infinite and, for each $B>1$, define
\[
E(B;\scL)
\colon=
\left\{ z\in \mfF: \|a_n\|\geq B^n \text{ for infinitely many }n\in \scL\right\}.
\]
\begin{coro01}\label{Coro729}
If $\scL\subseteq \Na$ is infinite and $B>1$, then $\dimh\left(E(B;\scL)\right) = s_B$.
\end{coro01}

\section{A complex {\L}uczak Theorem}\label{Sec:CxLuczak}

In this section, we prove Theorem \ref{TEO:LUC:USUALABS}. First, we observe that we can replace $|\,\cdot\,|$ with $\|\,\cdot\,\|$. For $b,c>1$, define the sets
\begin{align*}
\widetilde{\Xi}(b,c) &\colon=  \{\xi\in\mfF: c^{b^n}\leq \|a_n(\xi)\| \text{ for all } n\in\Na\}, \\
\Xi(b,c) &\colon=  \{\xi\in\mfF: c^{b^n}\leq \|a_n(\xi)\| \text{ for infinitely many } n\in\Na\}.
\end{align*}
On the one hand, since $\|z\|\leq |z|$ for all $z\in \Cx$, we have
\[
\widetilde{\Xi}(b,c)
\subseteq 
\{\xi\in\mfF: c^{b^n}\leq |a_n(\xi)| \text{ for all } n\in\Na\}.
\]
On the other hand, if $\veps>0$ is such that $c-\veps>1$, then every large $n\in\Na$ satisfies
\[
(c-\veps)^{b^n} \leq \frac{1}{\sqrt{2}} c^{b^n}
\]
and, since $|z|\leq \sqrt{2}\|z\|$ for all $z\in\Cx$, we conclude
\[
\{\xi\in\mfF: c^{b^n}\leq  |a_n(\xi) | \text{ for infinitely many } n\in\Na\}
\subseteq
\Xi(b,c- \veps).
\]
Hence, Theorem \ref{TEO:LUC:USUALABS} is equivalent to Theorem \ref{TEO:LUCZAK} below.

\begin{teo01}\label{TEO:LUCZAK}
If $b,c>1$, then
\[
\dimh\left( \Xi(b,c)\right)
=
\dimh\left( \widetilde{\Xi}(b,c)\right)
=
\frac{2}{b+1}.
\]
\end{teo01}
For the rest of this section, the numbers $b,c>1$ will remain fixed. In view of $\widetilde{\Xi}(b,c) \subseteq \Xi(b,c)$, we will obtain Theorem \ref{TEO:LUCZAK} from
\begin{equation}\label{Eq:BndsLuk}
\frac{2}{b+1}
\leq 
\dimh\left( \widetilde{\Xi}(b,c)\right)
\quad\text{ and }\quad
\dimh\left( \Xi (b,c)\right)
\leq 
\frac{2}{b+1}.
\end{equation}

\subsection{Upper bound of Hausdorff dimension}
This subsection is devoted to proving the upper bound in \eqref{Eq:BndsLuk}. Our argument is inspired in \cite{ShangFang24}. 
A previous version of this paper (available on arXiv) contains a proof extending {\L}uczak's original argument.

\begin{propo01}
Let $s,d\in\RE$ be such that 
\[
\frac{1}{b}< s < 1
\;\text{ and }\;
\frac{1}{s} < d < b.
\]
Then, for every $z\in \Xi(b,c)$ there are infinitely many $n\in\Na$ such that
\begin{equation}\label{Eq:Luk:UpBnd:01}
\|a_{n+1}(z)\|
\geq
\left( \|a_{1}(z)\|\cdots \|a_{n}(z)\|\right)^{sd-1}
c^{(1-s)d^{n+1}}.
\end{equation}
\end{propo01}
\begin{proof}
Take any $z\in \Xi(b,c)$. Let $m\in\Na$ be arbitrary and $k\in\Na_{>m}$ such that 
\[
\|a_{1}(z)\|\cdots \|a_{m}(z)\| < c^{b^k d^{m-k}}
\;\text{ and }\;
c^{b^k} < \|a_k(z)\|.
\]
Define the function $f_k\colon \Na\to \RE$ by
\[
f_k(N)= c^{b^kd^{N-k}} 
\quad\text{ for all } N\in\Na.
\]
Define
\[
n=\max\{r \in\Na: m\leq r\leq k-1, \; \|a_1(z)\|\cdots \|a_r(z)\|< f_k(r)\}.
\]
Then,
\begin{align*}
\|a_{1}(z)\| \cdots \|a_{n}(z)\|\|a_{n+1}(z)\|
&\geq
f_{k}(n+1) \nonumber\\
&=
f_{k}^s(n+1)f_{k}^{ 1 - s}(n+1) \nonumber\\
&>
f_{k}^{sd}(n)c^{(1-s)d^{n+1}} \nonumber\\
&>
\left( \|a_{1}(z)\|\cdots \|a_{n}(z)\| \right)^{sd} 
c^{(1-s)d^{n+1}} \nonumber.
\end{align*}
Divide both sides by $\|a_{1}(z)\|\cdots \|a_{n}(z)\|$ and the desired result follows.
\end{proof}
\begin{proof}[Proof of Theorem \ref{TEO:LUCZAK}: Upper bound]
For each $n\in\Na$ and $\bfa\in \Omega$, define the subset $\mathcal{K}_n(\bfa)$ of $\scD$ by
\[
\mathcal{K}_n(\bfa)
:=
\left\{ 
b\in \scD:
\bfa b\in\Omega(\bfa b), 
\; 
\| b\| \geq \left( \|a_{1}(z)\|\cdots \|a_{n}(z)\|\right)^{sd-1}
c^{(1-s)d^{n+1}} 
\right\}
\]
and the subset $K_n(\bfa)$ of $\omfF$ by 
\[
K_n(\bfa)
\colon=
\bigcup_{b\in \mathcal{K}_n(\bfa)}
\oclC_{n+1}(\bfa b)
=
Q_n
\left(
\bigcup_{b\in \mathcal{K}_n(\bfa)}
\oclC_{1}(b)\cap \mfF_n(\bfa);
\bfa
\right).
\]
We may show inductively that $|q_n(\bfa)|\geq \tau^{-n} \|a_1\|\cdots \|a_n\|$ for $\tau= \frac{\sqrt{2}}{\sqrt{2}-1}$. Therefore, Proposition \ref{Propo3.2}.\ref{Propo3.2_i} and the definition of $K_n(\bfa)$ yield
\begin{align*}
|K_n(\bfa)|
&\asymp 
\left| \clC_n(\bfa)\right|
\left|
\bigcup_{b\in \mathcal{K}_n(\bfa)}
\oclC_{1}(b)\cap \mfF_n(\bfa)
\right| \nonumber\\
&\asymp
\left| \clC_n(\bfa)\right|
\frac{1}{\left( \|a_{1}(z)\|\cdots \|a_{n}(z)\|\right)^{sd-1}
c^{(1-s)d^{n+1}}} \nonumber\\
&\asymp
\frac{1}{|q_n(\bfa)|^2}
\frac{1}{\left( \|a_{1}(z)\|\cdots \|a_{n}(z)\|\right)^{sd-1}
c^{(1-s)d^{n+1}}} \nonumber\\
&\ll
\frac{\tau^n}{\left( \|a_{1}(z)\|\cdots \|a_{n}(z)\|\right)^{sd+1}
c^{(1-s)d^{n+1}}}. \nonumber
\end{align*}
Let $F(b,c,s)$ be the set for which \eqref{Eq:Luk:UpBnd:01} holds for infinitely many $n\in\Na$, so $F(b,c,s)\supseteq \Xi(b,c)$. Then, if $t=\frac{2}{s(sd+1)}$, we have
\begin{align*}
\clH^t\left(F(b,c,s)\right)
&\leq 
\liminf_{N\to\infty} \sum_{n\geq N} \sum_{\bfa\in \Omega(n)}
\left| K_n(\bfa)\right|^t \nonumber\\
&\ll
\liminf_{N\to\infty} 
\sum_{n\geq N} \sum_{\bfa\in \Omega(n)}
\frac{\tau^{tn}}{\left( \|a_{1}\|\cdots \|a_{n}\|\right)^{\frac{2}{s}}
c^{(1-s)td^{n+1}}} \nonumber \\
&\ll
\liminf_{N\to\infty} 
\sum_{n\geq N} 
\frac{\tau^{tn}}{c^{(1-s)d^{n+1}}}
\sum_{\bfa\in \scD^n}
\frac{1}{\left( \|a_{1} \|\cdots \|a_{n} \|\right)^{\frac{2}{s}}} \nonumber\\
&=
\liminf_{N\to\infty} 
\sum_{n\geq N} 
\frac{1}{c^{(1-s)d^{n+1}}}
\left(\tau^{t} \sum_{a\in \scD} \frac{1}{ \|a \|^{\frac{2}{s}}}\right)^n \nonumber\\
&=
0.
\end{align*}
Therefore, $\dimh \Xi(b,c)\leq \dimh F(b,c,s)\leq t$ and the upper bound in \eqref{Eq:BndsLuk} follows by letting $s\to 1$ and $d\to b$.
\end{proof}
 
\subsection{Lower bound of Hausdorff dimension}
In this subsection, we obtain the lower bound in \eqref{Eq:BndsLuk} by applying the Mass Distribution Principle on a Cantor subset $E$ contained in $\widetilde{\Xi}(b,c)$.

\subsubsection{Construction of $E$}
Let $\kappa_1$ be as in Proposition \ref{Prop:SizeEstimates}.\ref{Prop:SizeEstimates:i} and take a natural number $M_0\in\Na$ such that 
\[
M_0
\geq 
\log_b\log_c \left( \frac{8\left(1 + \frac{\sqrt{2}}{2} \right)\left(2 + \frac{\sqrt{2}}{2} \right)}{\kappa_1\left(1 - \frac{\sqrt{2}}{2} \right)} + \frac{\sqrt{2}}{2}\right).
\]
For every $n\in\Na$, define the sets
\begin{align*}
\widetilde{D}_n
&\colon= 
\left\{d\in\scD: c^{b^{M_0+n+1}} \leq \|d \|\leq 2c^{b^{M_0+n+1}} \right\} \ \text{and} \ \Pre(n)\colon=
\prod_{j=1}^n \widetilde{D}_j.
\end{align*}

Note that, by $M_0\geq \log_b \log_c 3$ and Proposition \ref{PROP:SQRT8-REG:NUVO}, any sequence $(a_n)_{n\geq 1}$ in $\scD$ such that $(a_1,\ldots, a_n)\in \Pre(n)$ for all $n\in\Na$ is full.

Next, we define the layers of compact sets determining $E$. Put
\[
J_0\colon=\bigcup_{a_1\in\widetilde{D}_1} \oclC_1(a_1)
\]
and for every $n\in\Na$ and $\bfa\in \Pre(n)$ define
\[
J_n(\bfa)
:=
\bigcup_{b\in \widetilde{D}_{n+1}} \oclC_{n+1}(\bfa b).
\]
We define $E$ to be
\[
E
\colon=
\bigcap_{n\in\Na}\bigcup_{\bfa\in \Pre(n)} J_n(\bfa).
\]
We claim that $E\subseteq \widetilde{\Xi}(b,c)$. Certainly, let $z\in E$ be arbitrary. 
Consider any $n\in\Na$ and take $\bfa\in\Pre(n)$, $b\in \tilde{D}_{n+1}$ such that $z\in \oclC_{n+1}(\bfa b)$. By the proof of Proposition \ref{PROP:SQRT8-REG:NUVO}, $\oclC_{n+1}(\bfa b)\subseteq \clCc_n(\bfa)$, so the first $n$ partial quotients of $z$ are precisely $a_1,a_2,\ldots, a_n$. Hence, $z\in \widetilde{\Xi}(b,c)$. The lower bound in \eqref{Eq:BndsLuk} is, thus, a consequence of
\begin{equation}\label{Ec:LoBndLuc}
\dimh(E)\geq \frac{2}{1+b}.
\end{equation}
\subsubsection{A measure on $E$}
We define a probability measure $\tilde{\mu}$ on $E$ using a sequence of consistent measures. First, for every $a_1\in\widetilde{D}_1$, define
\[
\tilde{\mu}_1(J_1(a_1)) 
\colon= 
\frac{1}{\#\widetilde{D}_1}.
\]
Take $n\in\Na$ and assume that we have defined $\tilde{\mu}_n(J_n(\bfa))$ for all $\bfa\in\Pre(n)$. For any $a_{n+1}\in \widetilde{D}_{n+1}$ put
\[
\tilde{\mu}_{n+1} \left( J_{n+1}(\bfa a_{n+1})\right)
\colon=
\frac{1}{\#\widetilde{D}_{n+1}}
\tilde{\mu}_n\left(J_{n}(\bfa)\right).
\]
Since the measures $(\tilde{\mu}_n)_{n\geq 1}$ are consistent, we apply the Daniell-Kolmogorov Consistency Theorem to obtain a probability measure $\tilde{\mu}$ supported on $E$.

\subsubsection{Estimate of $\tilde{\mu}(J_n(\bfa))$}
Now, for $n\in\Na$ and $\bfa\in\widetilde{D}_n$, we estimate $\tilde{\mu}(J_n(\bfa))$ in terms of $|J_n(\bfa)|$. Write $t_0\colon= \frac{2}{1+b}$ and take any $0<t<t_0$.

\begin{propo01}\label{Propo:LUC:LB:EstMeas}
There exists some constant $\gamma_2 =\gamma_2(b,c) >0$ such that for all large $n\in\Na$ and all $\bfa\in\Pre(n)$ we have
\[
\tilde{\mu}\left( J_n( \bfa )\right) \leq \gamma_2  |J_n( \bfa )|^t.
\]
\end{propo01}
First, we estimate the number of elements in each $\#D_n$. Recall that $\#\{a\in \Za[i]: \|a\|\leq m\}=(2m+1)^2$ for every $m\in\Na$.
\begin{lem01}\label{Lem:Bnd:Dk}
For every $k\in\Na$, we have $\#\widetilde{D}_k \geq 8c^{2b^{M_0 + k +1}}$.
\end{lem01}
\begin{proof}
Take $k\in\Na$ and $d\in \scD$. Then, $d\in \widetilde{D}_k$ if and only if 
\[
\left\lceil c^{b^{M_0+k+1}} \right\rceil
\leq 
\|d\|
\leq 
\left\lfloor 2c^{b^{M_0+k+1}}\right\rfloor.
\]
As a consequence, we have
\begin{align*}
\#\widetilde{D}_k
&=
\left( 2\left\lfloor 2c^{b^{M_0 + k + 1}}\right\rfloor +1 \right)^2 
- 
\left( 2 \left\lceil c^{b^{M_0 + k + 1}} -1 \right\rceil + 1 \right)^2 \\
&\geq 
\left( 4c^{b^{M_0 + k + 1} } -1 \right)^2
-
\left( 2c^{b^{M_0+k+1}} + 1 \right)^2 \\
&=
12c^{2b^{M_0+k+1}}
\left( 1 - \frac{1}{c^{b^{M_0+k+1}}} \right) \\ 
&>
8 c^{2b^{M_0+k+1}}. && (\text{by } M_0\geq \log_b \log_c 3)
\end{align*}
\end{proof}

\begin{propo01}\label{Prop:Luczak:LoBo:01}
There exists some constant $\gamma_1=\gamma_1(b,c)>0$ such that for all $n\in\Na$ and all $\bfa\in\Pre(n)$ we have
\begin{equation}\label{EQ:LUC:LB:01}
\tilde{\mu}\left(J_n(\bfa )\right) 
\leq 
\frac{ \gamma_1 }{c^{tb^{M_0+n+1}}|q_n(\bfa)|^{2t}}.
\end{equation}
\end{propo01}
\begin{proof}
Take $n\in\Na$ and $\bfa\in \Pre(n)$. We obtain \eqref{EQ:LUC:LB:01} by showing that 
\begin{equation}\label{EQ:LUC:LB:02}
\log\left( \tilde{\mu}\left(J_n( \bfa )\right)\right) + \log\left( c^{tb^{M_0+n+2}}|q_n(\bfa)|^{2t}\right)
\end{equation}
is bounded from above. We deal with each term of the sum separately. First, we show that
\[
\log\left( \tilde{\mu}\left(J_n( \bfa )\right)\right) 
\leq  
-n\log 8 - b^{M_0 + n +1} \frac{2b}{b-1}\left( 1 - \frac{1}{b^n}\right) \log c.
\]
This follows from
\begin{align*}
\log\left( \tilde{\mu}\left(J_n( \bfa )\right)\right) 
&= 
\sum_{j=1}^n \log\left( \frac{1}{\#\widetilde{D}_j}\right) \\
&\leq 
\sum_{j=1}^n \log\left( \frac{1}{8} \right) + \log\left( \frac{1}{c^{2b^{M_0+j+1}} } \right) &&\text{(by Lemma \ref{Lem:Bnd:Dk})} \\
&=
-n\log 8 - 2b^{M_0 + 2} \log c \sum_{j=1}^n b^{j - 1}\\
&=
-n\log 8 - 2b^{M_0 + 2} \left( \frac{b^n-1}{b-1} \right) \log c \\
&=
-n\log 8 - b^{M_0 + n +1} \frac{2b}{b-1}\left( 1 - \frac{1}{b^n}\right) \log c.
\end{align*}
Now, we prove that the second term in \eqref{EQ:LUC:LB:02} verifies
\[
\log\left( c^{tb^{M_0+n+2}}|q_n(\bfa)|^{2t}\right)
\leq 
tb^{M_0+n+1}\log c  
\left( b + \frac{2b}{b-1} \left( 1 - \frac{1}{b^n}\right)  \right) 
+ 5nt\log 2.
\]
Indeed, from
\[
|q_n(\bfa)| 
\leq \prod_{j=1}^n \left(|a_j| + 1 \right)
\leq \prod_{j=1}^n \left(\sqrt{2}\| a_j \| + 1 \right)
\leq 8^{\frac{n}{2}} \prod_{j=1}^n \|a_j\|,
\]
we obtain
\[
\log |q_n(\bfa)|
\leq 
 \frac{n}{2}\log 8 + \sum_{j=1}^n \log \|a_j\|.
\]
Thus,
\begin{align*}
\log\left( c^{tb^{M_0+n+2}}|q_n(\bfa)|^{2t}\right)
&\leq 
tb^{M_0 + n + 2} \log c + 2t\left( \frac{n}{2}\log 8 + \sum_{j=1}^n \log \|a_j\|  \right)\\
&\leq 
tb^{M_0+n+2} \log c + 3tn\log 2 + 2t \left( \sum_{j=1}^n \log 2 + b^{M_0+j+1}\log c\right) \\
&=
tbb^{M_0+n+1} \log c + 5nt\log 2  +  2t b^{M_0+2} \left(  \frac{b^n-1}{b-1} \right)\log c \\
&=
tb^{M_0+n+1} \log c 
\left( b + \frac{2b}{b-1} \left( 1 - \frac{1}{b^n}\right)  \right) 
+ 5nt\log 2.
\end{align*}
Combining the estimates, we can bound the sum \eqref{EQ:LUC:LB:02} as follows:
\begin{align}
\log &\left( \tilde{\mu}\left(J_n( \bfa )\right)\right) + \log\left( c^{tb^{M_0+n+2}}|q_n(\bfa)|^{2t}\right) \nonumber\\
&<
b^{M_0+n+1}\log c  \left(
 tb + \frac{2tb}{b-1} \left( 1 - \frac{1}{b^n}\right)    -  \frac{2b}{b-1}\left( 1 - \frac{1}{b^n}\right) \right) +  n(5t \log 2 - \log 8).\label{EQ:LUC:LB:02:05}
\end{align}
Using the identity
\[
t_0b + \frac{2t_0b}{b-1}   -  \frac{2b}{b-1} = 0,
\]
we may see that the coefficient of $b^{M_0+n+1}$ in \eqref{EQ:LUC:LB:02:05} is negative for large $n$ and the sum in \eqref{EQ:LUC:LB:02} is bounded above.
\end{proof}
\begin{proof}[Proof of Proposition \ref{Propo:LUC:LB:EstMeas}]
Using the notation of Subsection \ref{Subsection:SecondGeomCons}, for $n\in \Na$ and $\bfa\in \widetilde{D}_n$ we have
\[
J_n(\bfa)
=
Q_n
\left( 
\mfTreb
\left( 
\left(\left\lfloor 2c^{b^{M_0+n+1}}\right\rfloor+\frac{1}{2}\right)^{-1}, 
\left(\left\lceil c^{b^{M_0+n+1}}\right\rceil - \frac{1}{2}\right)^{-1}
\right); \bfa
\right).
\]
Therefore, there are absolute constants for which
\[
\left| J_n(\bfa)\right|
\asymp
\frac{ 1}{c^{b^{M_0+n+1}}|q_n(\bfa)|^{2}},
\]
and the result follows from Proposition \ref{Prop:Luczak:LoBo:01}.
\end{proof}

\subsubsection{Estimates of $\tilde{\mu}(\Dx(z;r))$}
Next, we estimate the measure of small disc. To this end, define
\[
r_0 \colon= \frac{1}{2} \min\{|J_1(a)|:a\in \widetilde{D}_1\}.
\]
\begin{propo01}\label{PROP:BOUND:MU:DISC}
There is a constant $\gamma_3=\gamma_3(b,c)>0$ such that $\mu(\Dx(z;r))<\gamma_3r^t$ whenever $z\in E$ and $0<r<r_0$.
\end{propo01}
\begin{proof}
Take any $z=[a_1,a_2,\ldots] \in E$ and $r\in\RE$ such that $0<r< r_0$. Let $n$ be the single natural number verifying
\[
|J_{n+1}(a_1,\ldots,a_n, a_{n+1})|
<
r
\leq
|J_{n}(a_1,\ldots , a_{n})|.
\]
Define
\[
\widetilde{\gamma} \colon= \frac{1}{16\left( 1+ \frac{\sqrt{2}}{2}\right)\left( 2 +  \frac{\sqrt{2}}{2}\right)}.
\]
We consider two possible cases: either
\begin{equation}\label{EQ:LUC:LB:03}
|J_{n+1}(a_1,\ldots,a_n, a_{n+1})|
< 
r
\leq 
\widetilde{\gamma}\,|\clC_{n+1}(a_1,\ldots,a_n, a_{n+1})|
\end{equation}
or
\begin{equation}\label{EQ:LUC:LB:04}
\widetilde{\gamma}\,|\clC_{n+1}(a_1,\ldots , a_{n+1})|
<
r
\leq 
|J_{n}(a_1,\ldots , a_{n})|.
\end{equation}
Suppose that \eqref{EQ:LUC:LB:03} is true. We claim that $\Dx(z;r)$ intersects exactly one cylinder of level $n+1$, hence at most one set of the form $J_{n+1}(b_1,\ldots, b_{n+1})$. We prove this in two steps. First, we show that
\begin{equation}\label{EQ:LUC:LB:05}
\Dx\left( \frac{p_{n+1}(z)}{q_{n+1}(z)}; \frac{|\clC_{n+1}(a_1,\ldots, a_{n+1})|}{4\left( 1+ \frac{\sqrt{2}}{2}\right)\left( 2+ \frac{\sqrt{2}}{2}\right) } \right)
\subseteq 
\clCc_{n+1}(a_1,\ldots,a_{n+1}).
\end{equation}
We afterward verify that
\begin{equation}\label{EQ:LUC:LB:06}
\Dx(z;  \widetilde{\gamma}|\clC_{n+1}(a_1,\ldots, a_{n+1})|)
\subseteq
\Dx\left( \frac{p_{n+1}(z)}{q_{n+1}(z)}; \frac{|\clC_{n+1}(a_1,\ldots, a_{n+1})|}{4\left( 1+ \frac{\sqrt{2}}{2}\right)\left( 2+ \frac{\sqrt{2}}{2}\right) } \right).
\end{equation}
The inclusion \eqref{EQ:LUC:LB:06} follows directly from Proposition \ref{PROPO:GC01:QnDISCOS} and $\Dx(0;2^{-1})\subseteq \mfFc=\mfFc_{n+1}(a_1,\ldots,a_{n+1})$. To show \eqref{EQ:LUC:LB:05}, we note that the choice of $M_0$ implies 
\[
\left| T^{n+1}(z)\right|
< 
\frac{\kappa_1\left( 1 - \frac{\sqrt{2}}{2}\right)}{8\left( 1 + \frac{\sqrt{2}}{2}\right) \left( 2 + \frac{\sqrt{2}}{2}\right) },
\]
because
\[
\left| \frac{1}{T^{n+1}(z)}\right| 
\geq 
\| a_{n+2}(z) \| - \left|T^{n+2}(z)\right|
\geq c^{b^{M_0}} - \frac{\sqrt{2}}{2}
>
\frac{8\left( 1 + \frac{\sqrt{2}}{2}\right)\left( 2 + \frac{\sqrt{2}}{2}\right)}{\kappa_1 \left( 1 - \frac{\sqrt{2}}{2}\right)}.
\]
Therefore, using Proposition \ref{Prop:QnCorrespRule}, 
\begin{align*}
\left| \frac{p_{n+1}(z)}{q_{n+1}(z)} -  z \right|
&=
\left| Q_{n+1}(0;(a_1,\ldots, a_{n+1})) -  Q_{n+1}(T^{n+1}(z);(a_1,\ldots, a_n, a_{n+1})) \right| \\
&\leq \frac{|T^{n+1}(z)|}{|q_{n+1}(a_1,\ldots, a_{n+1})|^2\left| 1 + \frac{q_{n}}{q_{n+1}}T^{n+1}(z) \right| } \\
&\leq \frac{1}{8\left( 1 + \frac{\sqrt{2}}{2}\right)\left( 2 + \frac{\sqrt{2}}{2}\right)}|\clC_{n+1}(a_1,\ldots,a_n)|.
\end{align*}
As a consequence, for all $w\in \Dx(z;\widetilde{\gamma}\;|\clC_{n+1}(a_1,\ldots,a_{n+1})|)$ we have
\begin{align*}
\left| \frac{p_{n+1}(z)}{q_{n+1}(z)} -  w \right|
&\leq 
\left| \frac{p_{n+1}(z)}{q_{n+1}(z)} -  z \right| + |z - w| \\
&< 
\frac{1}{4\left( 1 + \frac{\sqrt{2}}{2}\right)\left( 2 + \frac{\sqrt{2}}{2}\right)}|\clC_{n+1}(a_1,\ldots,a_n)|.
\end{align*}
This proves \eqref{EQ:LUC:LB:05}, hence, $D(z;r)\subseteq \clCc_{n+1}(a_1,\ldots,a_{n+1})$. We may, thus, conclude that
\begin{align*}
\tilde{\mu}\left( \Dx(z;r)\right)
&=
\tilde{\mu}\left( \Dx(z;r) \cap J_{n+1}(a_1,\ldots, a_{n+1}) \right) \\
&\leq 
\tilde{\mu}\left(J_{n+1}(a_1,\ldots, a_{n+1}) \right) 
\leq
\gamma_2 |J_{n+1}(a_1,\ldots, a_{n+1})|^t 
\leq 
\gamma_2 r^t .
\end{align*}
Let us now assume \eqref{EQ:LUC:LB:04}. Define
\[
\tilde{\sfN}
\colon=
\#\left\{h\in\widetilde{D}_{n+1}: J_{n+1}(a_1,\ldots,a_n,h)\cap \Dx(z;r)\neq \vac\right\}.
\]
In view of
\[
\bigcup_{a\in \widetilde{D}_{n+1}} \oclC_1(a)
\subseteq
\Dx\left( 0;2^{-1}\right)
\subseteq
\mfFc,
\]
for each $a\in \widetilde{D}_{n+1}$ we apply $Q_{n+1}(\,\cdot\,;(a_1,\ldots, a_n,a))$ and obtain 
\[
J_{n+1}(a_1,\ldots, a_n,a)
\subseteq 
Q_{n+1}\left(\Dx\left( 0;2^{-1}\right);(a_1,\ldots, a_n,a)\right)
\subseteq
\clCc_{n+1}(a_1,\ldots, a_n,a).
\]
Therefore, since $\clCc_{n+1}(a_1,\ldots, a_n,a)$ contains only one set of the form $J_{n+1}(b_1,\ldots,b_{n+1})$,
\[
\tilde{\sfN}
\leq 
\#\left\{a\in\widetilde{D}_{n+1}: Q_{n+1}\left(\Dx\left( 0;2^{-1}\right);(a_1,\ldots, a_n,a)\right)\cap \Dx(z;r)\neq \vac\right\}.
\]
Recall that $\leb$ denotes the Lebesgue measure on $\Cx$. Then, there is an absolute constant $\gamma_3>0$ such that
\begin{align*}
\leb\left( Q_{n+1}\left(\Dx\left( 0;2^{-1}\right);(a_1,\ldots, a_n,a)\right)\right) &\geq 
\gamma_3 \frac{\pi}{4} \frac{1}{|q_{n+1}(a_1,\ldots, a_n,a)|^4} \\
&\geq 
\gamma_3 \frac{\pi}{36} \frac{1}{|q_{n}(a_1,\ldots, a_n)|^4|a|^4} \\
&\geq 
\gamma_3 \frac{\pi}{144} \frac{1}{|q_{n}(a_1,\ldots, a_n)|^4c^{4b^{M_0+n+2}}},
\end{align*}
which gives
\[
\tilde{\sfN} 
\leq 
\frac{144}{ \gamma_3}
|q_{n}(a_1,\ldots, a_n)|^4c^{4b^{M_0+n+2}} r^2.
\]
From $0<t<1$ we get $0<\frac{t}{2}<\frac{1}{2}$, hence
\begin{align*}
\tilde{\mu} \left(\Dx(z;r)\right) 
&\leq 
\tilde{\mu} \left(J_n(a_1,\ldots, a_n)\right)
\min\left\{ 1, \tilde{\sfN} \frac{1}{\#\widetilde{D}_{n+1}} \right\} \\
&\leq 
\tilde{\mu} (J_n(a_1,\ldots, a_n)
\min\left\{ 1, \frac{144}{ \gamma_4  }c^{4b^{M_0 + n + 2}}|q_n(\bfa)|^4 r^2 \frac{1}{8c^{2b^{M_0+ n + 2}}}  \right\} \\
&= 
\tilde{\mu} \left(J_n(a_1,\ldots, a_n)\right)
\min\left\{ 1, 18 c^{2b^{M_0 + n + 2}}|q_n(\bfa)|^4 r^2   \right\} \\
&\leq 
\tilde{\mu}\left(J_n(a_1,\ldots, a_n)\right)
18^{\frac{t}{2}} c^{tb^{M_0 + n + 2}}|q_n(\bfa)|^{2t} r^t \\
&<
18\gamma_2 r^t. 
\end{align*}
\end{proof}
\begin{proof}[Proof of Theorem \ref{TEO:LUCZAK}. Lower bound]
The Mass Distribution Principle and Proposition \ref{PROP:BOUND:MU:DISC} imply \eqref{Ec:LoBndLuc}. Since $E\subseteq \widetilde{\Xi}(b,c)$, the lower bound in \eqref{Eq:BndsLuk} follows.
\end{proof}

\begin{rema01}
A slight modification of the proof of the lower bound in Theorem \ref{TEO:LUCZAK} tells us that for every $b,c>1$ and $\delta>0$ we have
\[
\dimh\left( \widetilde{\Xi}(b,c)\setminus \widetilde{\Xi}(b,c+\delta) \right )
=
\dimh\left( \Xi(b,c)\setminus \Xi(b,c+\delta) \right )
=
\frac{2}{b+1}.
\]
\end{rema01}

\section{Proof of Theorem \ref{TEO:MAIN}}\label{Sec:TeoMain}
Recall that 
\[
E_{\infty}(\Phi)
\colon=
\left\{ z\in \mathfrak{F}: \|a_n(z)\|\geq \Phi(n) \text{ \rm for infinitely many }n\in\mathbb{N} \right\}.
\]
\begin{enumerate}[\rm i.]
\item Assume that $B=1$. Let $\veps>0$ be arbitrary and consider the infinite set
\[
\scL_{\veps}
\colon=
\left\{ n\in\Na : \Phi(n) \leq (1+\veps)^n\right\}.
\]
Then, we have
\[
\left\{ z\in \mfF: \|a_n\|\geq (1+\veps)^n \text{ for infinitely many }n\in \scL\right\}
\subseteq  E_{\infty}(\Phi), 
\]
so $s_{1+\veps}\leq \dimh(E_{\infty}(\Phi))$, by Corollary \ref{Coro729}. In view of Proposition \ref{Le:LemPrel-04}, we conclude that
\[
2 = \lim_{\veps \searrow 0} s_{1+\veps}\leq \dimh(E_{\infty}(\Phi)) \leq 2.
\]
\item Assume that $1<B<\infty$. For each $\veps>0$, define the infinite set
\[
\scL_{\veps}'
\colon=
\left\{ n\in\Na: \Phi(n) \leq (B+ \veps)^n\right\}.
\]
Hence,
\[
\left\{z\in \mfF: \|a_n\|\geq (B+\veps)^n \text{ for infinitely many } n\in \scL_{\veps}'\right\} 
\subseteq
E_{\infty}(\Phi)
\]
and
\[
E_{\infty}(\Phi)
\subseteq
\left\{z\in \mfF: \|a_n\|\geq (B - \veps)^n \text{ for infinitely many } n\in \Na \right\}.
\]
Therefore, $s_{B-\veps} \leq \dimh(E_{\infty}(\Phi))\leq s_{B+\veps}$ and, since $B\mapsto s_B$ is continuous (Proposition \ref{Le:LemPrel-04}), $\dimh(E_{\infty}(\Phi))=s_B$.

\item Assume that $B=\infty$. Further, assume that $b=1$. Then, given $\veps>0$, the set
\[
\scL 
\colon=
\left\{n\in\Na: \Phi(n) \leq e^{(1+\veps)^n} \right\}
\]
is infinite and we have
\begin{align*}
\left\{ z\in\mfF: \|a_n(z)\|\geq e^{(1+\veps)^n} \text{ for all }n\in\Na\right\}
&\subseteq
\left\{ z\in\mfF: \|a_n(z)\|\geq e^{(1+\veps)^n} \text{ for all }n\in\scL \right\} \\
&\subseteq E_{\infty}(\Phi),
\end{align*}
so, Theorem \ref{TEO:LUCZAK} gives
\[
\dimh\left( E_{\infty}(\Phi)\right) 
\geq \frac{2}{2+\veps}.
\]
Letting $\veps\to 0$, we conclude $\dimh\left( E_{\infty}(\Phi)\right) \geq 1$. Moreover, $B=\infty$ implies that, for any $C>1$, the inequality $\frac{\log \Phi(n)}{n}\geq \log C$ holds for infinitely many $n\in\Na$, hence
\[
E_{\infty}(\Phi)
\subseteq
\left\{ z\in\mfF: C^n\leq \|a_n(z)\| \text{ for infinitely many } n\in\Na \right\} 
\]
and $\dimh\left( E_{\infty}(\Phi) \right) \leq s_{C}$. We let $C\to\infty$ and recourse to Lemma \ref{Le:LemPrel-04} to conclude $\dimh\left( E_{\infty}(\Phi) \right) \leq 1$, so
\[
\dimh\left( E_{\infty}(\Phi) \right) = 1.
\]
Suppose that $1<b<\infty$. Then, for any $\veps>0$, the inequality $\Phi(n)\leq e^{(b+\veps)^n}$ holds infinitely often while $\Phi(n)\geq e^{(b-\veps)^n}$ is true for every large $n$. As a consequence, since the set
\[
\left\{ z\in \mfF: \|a_n(z)\|\geq e^{(b+\veps)^n} \text{ for all } n\in\Na \right\}
\]
is contained in
\[
\left\{ z\in \mfF: \|a_n(z)\|\geq e^{(b+\veps)^n}\geq \Phi(n) \text{ for infinitely many } n\in\Na \right\}
\subseteq 
E_{\infty}(\Phi)
\]
and
\[
E_{\infty}(\Phi)
\subseteq 
\left\{ z\in \mfF: \|a_n(z)\|\geq \Phi(n) \geq e^{(b-\veps)^n} \text{ for infinitely many } n\in\Na \right\},
\]
Theorem \ref{TEO:LUCZAK} tells us that
\[
\frac{2}{b + 1 + \veps}
\leq 
\dimh\left( E_{\infty}(\Phi)\right)
\leq
\frac{2}{b+1-\veps}.
\]
Letting $\veps\to 0$, we arrive at $\dimh\left( E_{\infty}(\Phi)\right) = \frac{2}{b+1}$.

To finish the proof, assume that $b= \infty$. Then, $\Phi(n)\geq e^{c^n}$ for all $c>0$ and every sufficiently large $n\in\Na$; that is, 
\[
E_{\infty}(\Phi)
\subseteq
\left\{ z\in\mfF: \| a_n(z)\|\geq e^{c^n} \text{ for infinitely many }n \in\Na \right\}.
\]
Theorem \ref{TEO:LUCZAK} yields $\dimh\left(E_{\infty}(\Phi)\right) \leq \frac{2}{1 + c}$ and, taking $c\to \infty$, we conclude that $\dimh\left(E_{\infty}(\Phi)\right) =0$. This completes the proof of Theorem \ref{TEO:MAIN}.
\end{enumerate}

\section{Concluding remarks and open problems}\label{SECTION:ConcludingRemarks}
In the recent literature, many new metrical results (including zero-one laws and computation of 
Hausdorff dimensions) on sets of real numbers $x$ defined in terms 
of their partial quotients $a_n (x)$ and/or of their convergents $p_n(x)/q_n(x)$ have been established. 
We survey some of them below and briefly discuss their complex analogs, where regular continued fractions of real numbers are replaced by Hurwitz continued fractions of complex numbers. 

As discussed in the introduction, stemming from Borel-Berstein and Wang-Wu theorems, the focus is on the growth properties of partial quotients. To this end, Sun and Wu \cite{SunWu14}  established that, for any $b > 0$, the Hausdorff dimension of the level set 
$$
\Bigl\{ x \in (0, 1) : \lim_{n \to + \infty} \frac{\log a_{n+1} (x)}{\log q_n (x)} = b \Bigr\}
$$
is equal to $1 / (b+2)$. A localized version (where $b$ is replaced by a continuous function of $x$) has subsequently 
been proved in \cite{ChSunZh}, see \cite{ZhSo21} for a further extension. 
Their complex analogs should hold. To establish it, 
one first needs to prove the complex analog of \cite[Lemma 3.2]{FLWW09}. Specifically, in \cite{FLWW09}, the authors calculated the Hausdorff dimension of the level set 
$$ 
E(\alpha)=\left\{x\in [0, 1): \lim_{n\to\infty}\frac{\log a_1(x)+ \cdots+\log a_n(x)}{n}=\alpha\right\}, 
$$ 
for any given $\alpha>0$. They proved, quite unexpectedly, that the graph of the function $\alpha \mapsto \dim E(\alpha)$ 
is neither convex nor concave. 
Later, Fan et al. \cite{FLWW13}  
considered the subsets of $E(\infty)$ defined by 
$$ 
E(\psi, \alpha)=\left\{x\in [0, 1): \lim_{n\to\infty}\frac{\log a_1(x)+\cdots+\log a_n(x)}{\psi(n)}=\alpha\right\} \text{ for all } \alpha>0, 
$$ 
where $\psi : \Na\to\Na$ is such that $\psi(n)/n$ tends to $\infty$ with $n$. 
They obtained the Hausdorff dimension of $E(\psi, \alpha)$. Subsequently, Liao and Rams \cite{LiRa16} 
computed the Hausdorff dimension of sets defined as $E(\psi, \alpha)$, but with 
the limit replaced by $\liminf$ (resp., by $\limsup$) 
and thereby exhibited a new phenomenon. 
It would be desirable to establish the analogs of all these results for Hurwitz continued fractions and to 
determine whether the unexpected phenomena also occur in this setting.

Very recently, He and Xiao \cite{HeXiao24} proved results concerning 
the sets of complex numbers defined by the
growth rate of products of consecutive partial quotients. Thereby they obtained complex analogs of results 
established for real continued fractions in \cite{HuWuXu20}. They observe that, in order to extend 
their Theorem 1.3 to functions that are not continuous, the auxiliary results on the geometry of irregular cylinders 
established in the present paper may be useful. 

It is worth mentioning that the study of consecutive partial quotients gives information regarding the improbability (or not) of Dirichlet's theorem, see \cite{BHS, HKWW2018, KleinbockWadleigh}. However, interestingly products of partial quotients arising from the arithmetic progressions give information regarding the set of exceptions to the convergence criterion of Host-Kra (2005) and Bourgain (1989) of multiple ergodic averages for the Gauss map \cite{HussainShulga2}. Establishing analogs of these results for Hcfs is possible by using the geometric properties proved in this paper, along with some other ideas.


There is a vast literature on the determination of the 
Hausdorff dimension of sets of real numbers whose partial quotients are all in a given finite (or infinite) set; see 
e.g. \cite{DFSU23, He92} and the references given therein. 
Getting a complex analog seems to be  a challenging problem because the symbolic space associated to Hcfs cannot be modelled as a countable Markov shift. 

Chevallier \cite{Che21} used flows on spaces of lattices in $\mathbb{C}^2$ to prove the analog of Dirichlet's theorem for 
approximation in $\QU({i})$. Let $c(z)$ denote the smallest real number such that, for every real number $Q>1$, there exist 
Gaussian integers $p, q$ such that 
$$
0 < |q| < Q, \quad |qz - p| \le \frac{c(z)}{Q}.
$$
He proved that $c(z) \le \sqrt{2} / (3 - \sqrt{3}) = 1.115\ldots$, with equality for almost all $z$ 
(in the sense of the Lebesgue measure). 
It would be interesting to find an alternative proof, in the spirit of
Davenport and Schmidt \cite{DaSc70}, which may also give an explicit 
formula to compute $c(z)$ in terms of the Hurwitz continued fraction of $z$. 
Another related question is the determination of the smallest possible constant $c(z)$. 
In the real case, this corresponds to the constant associated with the Golden Ratio \cite{DaSc70}.  

We finalize the paper with a simple question motivated by theorems \ref{TEO:LUCZAK:R} and \ref{TEO:LUC:USUALABS}. Let $\Phi:\Na\to\RE$ be an arbitrary function. We have defined
\[
E(\Phi)
\colon=
\left\{ z\in \mathbb C: |a_n(z)|\geq \Phi(n) \text{ for infinitely many }n\in\mathbb{N} \right\}.
\]
Recall that for any $x\in (0,1)\setminus\QU$ and $n\in\Na$ we denote $n$-th partial quotient of the regular continued fraction of $x$ as $A_n(x)$. Define the set
\[
E_{\RE}(\Phi)
\colon=
\left\{ x\in \mathbb{R}: |A_n(x)|\geq \Phi(n) \text{ for infinitely many }n\in\mathbb{N} \right\}.
\]
Theorems \ref{TEO:LUCZAK:R} and \ref{TEO:LUC:USUALABS} tell us that $\dimh(E(\Phi))=2\dimh(E_{\RE}(\Phi))$ holds whenever $\Phi$ has a subsequence of very fast growth. Whether it would hold for any $\Phi:\Na\to\RE$ is a natural open problem.


\bibliography{referencias}
\bibliographystyle{siam}

\end{document}